\documentclass[12pt]{amsart}
\usepackage{amsmath,amssymb,times,color}
\usepackage{graphicx}
\usepackage{fig4tex}
\usepackage{dsfont}
\definecolor{webred}{rgb}{0.75,0,0}
\definecolor{webgreen}{rgb}{0,0.75,0}
\usepackage[citecolor=webgreen,colorlinks=true,linkcolor=webred]{hyperref}
\usepackage{stmaryrd} 
\usepackage{mathrsfs}
\usepackage{fig4tex}

\newtheorem{thm}{Theorem}[section]

\newtheorem{lem}[thm]{Lemma}

\newtheorem{hyp}[thm]{Assumption}

\theoremstyle{definition}
\newtheorem{defn}[thm]{Definition}

\theoremstyle{remark}
\newtheorem{rem}[thm]{Remark}

\oddsidemargin=\evensidemargin

\numberwithin{equation}{section}

\renewcommand{\Re}{\operatorname{Re}}

\renewcommand{\Im}{\operatorname{Im}}

\newcommand{\Div}{\operatorname{\mathrm{div}}}

\newcommand{\rot}{\operatorname{\mathrm{curl}}}

\newcommand{\db}[1]{_{\raise-0.3ex\hbox{$\scriptstyle #1$}}}
\newcommand{\dd}[1]{_{\raise-1.5pt\hbox{$\scriptstyle #1$}}}

\newcommand{\di}{\displaystyle}
\newcommand{\te}{\textstyle}
\newcommand{\dr}{{\rm d}}
\newcommand{\m}{\mathsf{m}}

\newcommand{\iso}{_{\sf +}}
\newcommand{\con}{_{\sf -}}
\newcommand{\is}{^{\sf +}}

\newcommand  {\C}{{\mathbb C}}
\newcommand  {\N}{{\mathbb N}}
\renewcommand{\P}{{\mathbb P}}
\newcommand  {\Q}{{\mathbb Q}}
\newcommand  {\R}{{\mathbb R}}
\newcommand  {\T}{{\mathbb T}}
\newcommand  {\Z}{{\mathbb Z}}
\newcommand {\Id}{\mathbb {I}}

\newcommand{\BB}{\boldsymbol{\mathsf B}}
\newcommand  {\EE}{\boldsymbol{\mathsf E}}

\newcommand  {\GG}{\boldsymbol{\mathsf G}}
\newcommand  {\HH}{\boldsymbol{\mathsf H}}
\newcommand  {\KK}{\boldsymbol{\mathsf K}}
\newcommand  {\LL}{\boldsymbol{\mathsf L}}
\renewcommand  {\L}{{\mathrm L}}

\newcommand  {\RR}{\boldsymbol{\mathsf R}}
\newcommand  {\TT}{\boldsymbol{\mathsf T}}
\newcommand  {\ZZ}{\boldsymbol{\mathsf Z}}

\newcommand  {\VV}{\,\underline{\!{\boldsymbol{\mathfrak H}}\!}\,}

\newcommand  {\WW}{\,\underline{\!{\boldsymbol{\mathfrak E}}\!}\,}

\newcommand  {\XX}{\boldsymbol{\tau}}

\renewcommand{\aa}{{\boldsymbol{\mathsf a}}}

\newcommand  {\jj}{{\boldsymbol{\mathsf j}}}
\newcommand  {\nn}{\boldsymbol{\mathsf n}}

\newcommand  {\uu}{\boldsymbol{\mathsf u}}

\newcommand  {\xx}{\boldsymbol{\mathsf x}}
\newcommand  {\yy}{\boldsymbol{\mathsf y}}

\newcommand  {\cH}{H}

\newcommand  {\cL}{\mathcal{L}}
\newcommand  {\cO}{\mathcal{O}}
\newcommand  {\cU}{{\mathcal U}}

\newcommand  {\cV}{{\mathfrak H}}
\newcommand  {\cW}{{\mathfrak E}}

\newcommand  {\bH}{\mathbf{H}}
\newcommand  {\bL}{\mathbf{L}}

\newcommand  {\sH}{\mathsf{h}}
\newcommand  {\sG}{\mathsf{g}}
\newcommand  {\sj}{\mathsf{j}}

\newcommand  {\sv}{\mathfrak{h}}
\newcommand  {\fke}{\mathfrak{e}}

\newcommand  {\svV}{\mathfrak{h}^\theta}
\newcommand  {\gm}{\mathfrak{m}}

\newcommand  {\sw}{\mathsf{w}}
\newcommand  {\bsH}{\boldsymbol{\mathsf{h}}}

\newcommand  {\bs}{\underline \sigma}

\newcommand{\A}{\mathsf A}
\newcommand{\B}{\mathsf B}
\newcommand{\D}{\mathsf D}
\renewcommand{\P}{\mathsf P}
\renewcommand{\Q}{\mathsf Q}

\newcommand  {\MM}{{\mathfrak M}}

\newcommand  {\err}{\mathrm{err}}

\newcommand  {\brH}{\breve{\HH}}

\newcommand  {\tH}{\tilde{\mathsf{h}}}

\definecolor{mpurple}{rgb}{0.6,0,0.8}
\definecolor{myblue}{rgb}{0.,0.2,0.8}
\definecolor{mygreen}{rgb}{0,0.75,0.0}
\definecolor{mred}{rgb}{0.9,0,0}
\definecolor{mbrun}{rgb}{0.8,0.5,0}

\newcommand{\Bk}{\color{black}}

\begin{document}

\title{On the influence of the geometry on skin effect in electromagnetism}

\author{Gabriel Caloz, Monique Dauge, Erwan Faou, Victor P\'eron}

\begin{abstract}
We consider the equations of electromagnetism set on a domain made of a dielectric and a conductor subdomain in a regime where the conductivity is large. Assuming smoothness for the dielectric--conductor interface, relying on recent works we prove that the solution of the Maxwell equations admits a multiscale asymptotic expansion with profile terms rapidly decaying inside the conductor. This {\em skin effect} is measured by introducing a skin depth function that turns out to depend on the mean curvature of the boundary of the conductor. We then confirm these asymptotic results by numerical experiments in various axisymmetric configurations. We also investigate numerically the case of a nonsmooth interface, namely a cylindrical conductor. 
\end{abstract}

\maketitle

\tableofcontents

\section{Introduction}
Our interest lies in the influence of the geometry of a conducting body  on the skin effect in electromagnetism. This effect describes the rapid decay of electromagnetic fields with depth inside a metallic conductor. The skin effect reflects the flow of current near the surface of a conductor. After the early work \cite{Rytov40}, the mathematical analysis of the skin effect has been addressed more recently in several papers, \cite{S83,MS84,MS85,HJN08}.

The present work is motivated by recent studies \cite{HJN08,DFP09,CDP09,Pe09} in which authors analyze the behavior of the electromagnetic fields solution of the Maxwell equations  through an asymptotic expansion for large conductivity. In particular, uniform estimates for the electromagnetic field at high conductivity are proved in \cite{CDP09}, whereas in the note \cite{DFP09} a suitable skin depth function is introduced on the interface between a conductor and an insulator  to generalize the classical scalar quantity. An asymptotic expansion at high conductivity for this function shows the influence of the geometry of the interface : the skin depth is larger for high conductivity when the mean curvature of the conducting body surface is larger -- and here the sign of the curvature has a major influence, which means that the skin depth is larger in convex than in concave conductors.

\medskip
In this paper, our aim is twofold
\begin{enumerate}
\item Present elements of derivation for asymptotic expansions near the conductor-insulator interface,
\item Illustrate by numerical computations the theoretical behavior deduced from asymptotic analysis in \cite{DFP09}.
\end{enumerate} 

For our computations, we consider a special class of axisymmetric problems, which allows their reduction to two-dimensional scalar problems. This enables to measure the skin effect with reduced computational effort. We perform finite element computations for relevant benchmarks (cylinders and spheroids). The interface problem for the magnetic field is  solved with the library M\'{e}lina \cite{melina}. A postprocessing shows the accuracy of the asymptotic expansion, and exhibits the influence of the mean curvature of the conductor on the skin effect. Our computations also clearly display the effect of the edges of the cylindrical conductor: In this case, the decay is not exponential near edges while becoming exponential further away from edges.

The presentation of the paper proceeds as follows. In section \ref{S2}, we introduce the framework: we present the Maxwell equations and the result of \cite{CDP09} about existence of solution at high conductivity and the formulation for the magnetic field. In section \ref{S3}, we present the asymptotic expansion of the magnetic field in the smooth case, and the asymptotic behavior of the skin depth for high conductivity. In section \ref{S4}, we restrict our considerations to axisymmetric domains and orthoradial axisymmetric data, and we present the configurations A (cylindrical), B, and C (spheroidal), chosen for computations. In section \ref{S5}, we introduce a finite element discretization for the solution of the problem  in the meridian domain, and  we check the convergence of the  discretized problem. In section \ref{S6}, we present numerical simulations in all configurations to highlight the skin effect, and to exhibit the influence of the mean curvature of the conductor on the skin effect. In section \ref{S7}, we perform post-treatments of our numerical computations in configurations A (cylindrical) and B (spheroidal) in order to investigate the nature of the decay of the field inside the conductor: In smooth configurations (B) the rate of the exponential decay is very close to the expected theoretical one, while in corner configurations (A) the exponential decay shows up in a region which is not very close to the corner. Theoretical aspects in this latter case are going to be  studied in \cite{DaPePo10}. In Appendix A, we provide elements of proof for the multiscale expansion given in section \ref{S3}. We expand the Maxwell operators in power series, and we obtain the equations satisfied by the magnetic and electric profiles. In Appendix B, we derive the profiles of the orthoradial component of the magnetic field in relation with the 3D asymptotic expansion in axisymmetric configurations.

\section{Framework}
\label{S2}
Let $\Omega$ be a piecewise smooth Lipschitz domain. We denote by $\Gamma$ its boundary and by $\nn$ the outer unit normal field on $\Gamma$. We denote by $\bL^2(\Omega)$ the space of three-component fields square integrable on $\Omega$.

\begin{figure}[ht]
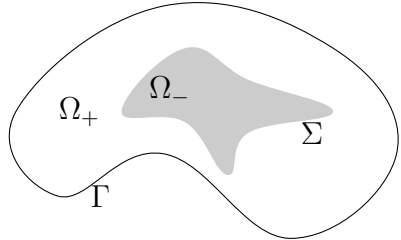

\figinit{0.8truept}
\figpt 1:(-100,20)
\figpt 2:(-30,70)\figpt 3:(50,50)
\figpt 4:(80,0)\figpt 5:(30,-40)
\figpt 6:(-30,0)\figpt 7:(-80,-20)
\figpt 8:(-50,20)\figpt 9:(-20,50)\figpt 10:(50,20)\figpt11:(10,35)
\figpt 12:(10,10)\figpt 13:(0,-10)
\figpt 14:(-20,15)%
\figpt 16:(-70,20)\figpt 15:(-10,30)
\figpt 17:(40,10) \figpt 18:(-60,-20)
\psbeginfig{}
\pscurve[1,2,3,4,5,6,7,1,2,3]
\pssetfillmode{yes}\pssetgray{0.8}
\pscurve[8,9,11,10,12,13,14,8,9,11]
\pssetfillmode{no}\pssetgray{0}
\psendfig
\figvisu{\figBoxA}{
}{
\figwritew 15: $\Omega_{-}$(6pt)
\figwritec [16]{$\Omega_{+}$}
\figwritec [17]{$\Sigma$}
\figwritec [18]{$\Gamma$}
}
\centerline{\box\figBoxA}
\vskip -6mm
\caption{The domain $\Omega$ and its subdomains $\Omega_+$ (dielectric) and $\Omega_-$ (conductor).}
\label{F1}
\end{figure}

We consider the Maxwell equations given by Faraday's and Amp\`ere's laws in $\Omega$:
\begin{equation}
    \rot \ \EE - i \omega\mu_0 \HH = 0 \quad \mbox{and}\quad
    \rot \ \HH + (i\omega\varepsilon_0 - \bs) \EE =   \jj
    \quad\mbox{in}\quad\Omega\, .
\label{MS}
\end{equation}
Here, $(\EE,\HH)$ represents the electromagnetic field, $\mu_{0}$ is the magnetic permeability, $\varepsilon_0$ the electric permittivity, $\omega$ the angular frequency, $\jj$ represents a current density and is supposed to belong to  
\[
   \bH_0(\Div,\Omega) =  \{\uu\in\bL^2(\Omega)\; |\ \Div\uu\in\L^2(\Omega),\ 
   \uu\cdot\nn=0 \ \mbox{on} \ \Gamma \} \ ,
\]
and $\bs$ is the electric conductivity. We assume that the domain $\Omega$ is made of two (connected) subdomains $\Omega_+$ and $\Omega_-$ in which the coefficient $\bs$ take two different values $(\sigma_+=0,\sigma_-\equiv\sigma)$. We denote by $\Sigma$ the interface between the subdomains $\Omega_+$ and $\Omega_-$. In the most part of our analysis\footnote{Except in one configuration in which the roles of the conductor and dielectric bodies are swapped.}, we will assume that the external boundary $\Gamma$ is contained in the dielectric region $\overline\Omega_+$, see Figure~\ref{F1}. 

To complement the Maxwell harmonic equations \eqref{MS}, we consider the perfectly insulating electric boundary conditions on the boundary $\Gamma$
\begin{equation}
\label{PIbc}
   \EE\cdot\nn=0 \quad\mbox{and}\quad \HH\times\nn=0\quad\mbox{on}\quad\Gamma\,.
\end{equation}
Note that we could also consider the perfectly conducting electric boundary conditions $\EE \times\nn=0$ and $\HH\cdot\nn=0$ on $\Gamma$.

\subsection{Existence of solutions}
Subsequently, we assume that the following condition on the limit problem in the dielectric part $\Omega_+$ is valid: 
\begin{hyp}
\label{H1}
The angular frequency $\omega$ is not an eigenfrequency of the problem
\begin{equation}
\label{PH1}
 \left\{
   \begin{array}{lll}
    \rot  \EE - i \omega\mu_0 \HH = 0 \quad \mbox{and}\quad
    \rot \HH + i\omega\varepsilon_0 \EE = 0 \quad&\mbox{in}\quad \Omega_{+}
\\[0.5ex]
\EE\times \nn  = 0 \quad \mbox{and} \quad
\HH\cdot\nn=0  \quad &\mbox{on}\quad \Sigma
\\[0.5ex]
\eqref{PIbc}  \quad &\mbox{on}\quad \Gamma.
   \end{array}
    \right.
\end{equation}
\end{hyp}
Hereafter, we denote by  $\|\cdot \|_{0,\mathcal{O}}$ the norm in $\bL^2(\mathcal{O})$. 
We quote from \cite[Th.2.3]{CDP09}:
\begin{thm}
\label{2T0}
If the interface $\Sigma$ is Lipschitz, under Assumption {\em\ref{H1}}, there are constants $\sigma_{0}$ and $C>0$, such that for all $\sigma\geqslant\sigma_0$, the Maxwell problem \eqref{MS} with boundary condition \eqref{PIbc} and data $\jj\in\bH_0(\Div,\Omega)$ has a unique solution $(\EE,\HH)$ in $\bL^2(\Omega)^2$, which satisfies:
\begin{equation}
   \|\EE\|_{0,\Omega} + \|\HH\|_{0,\Omega} 
   + \sqrt\sigma\, \|\EE\|_{0,\Omega_-} \leqslant C \|\jj\|_{\bH(\Div,\Omega)}.
\label{3E1}
\end{equation}
\end{thm}

For convenience, we introduce in the Maxwell equations \eqref{MS} the small parameter
\begin{equation}
\label{delta}
   \delta = \sqrt{{\omega\varepsilon_0} /{\sigma}}\,.
\end{equation}
Hence, $\delta$ tends to $0$ when $\sigma\to\infty$. For $\sigma\geqslant\sigma_0$, we denote by $(\EE_{(\delta)},\HH_{(\delta)})$ the solution of the system \eqref{MS} --\eqref{PIbc}.

\subsection{Magnetic formulation}
By a standard procedure we deduce from the Maxwell system \eqref{MS}-\eqref{PIbc} the following variational formulation for the magnetic field $\HH_{(\delta)}$. The variational space is $\bH_0(\rot,\Omega)$:
\begin{equation}
 \bH_0(\rot,\Omega) =\{\uu\in\bL^2(\Omega)\; |\ \rot\uu\in\bL^2(\Omega),\, 
 \uu\times\nn = 0  \ \mbox{on}\ \Gamma\}\ , 
\end{equation}
and the variational problem writes
\\[1ex]
{\it Find $\HH_{(\delta)} \in \bH_0(\rot,\Omega)$ such that for all $\KK \in \bH_0(\rot,\Omega)$}
\begin{equation}
\int_{\Omega} \big(\frac1{\varepsilon(\delta)}\rot\HH_{(\delta)} \cdot \rot\KK- \kappa^2\HH_{(\delta)}\cdot\KK\big)\dr\xx = \int_{\Omega} \big( \rot\frac{\jj}{\varepsilon(\delta)} \big) \cdot\KK\,\dr\xx\ ,
\label{FVH0}
\end{equation}
where we have set
\begin{equation}
\label{2Eeps}
   \varepsilon(\delta)={\mathbf{1}_{\Omega\iso}}  
   + ({1+\frac{i}{\delta^2}})\,{\mathbf{1}_{\Omega\con}}
   \quad\mbox{and}\quad
   \kappa:=\omega\sqrt{\varepsilon_{0}\mu_{0}}\,.
\end{equation}

\begin{hyp}
\label{H2}
We assume that the surfaces $\Sigma$ (interface) and $\Gamma$ (external boundary) are smooth.
\end{hyp}
Under Assumption \ref{H2}, the magnetic field $\HH_{(\delta)}$ is globally in $\bH^1(\Omega)$.

\section{Multiscale expansion}
\label{S3}

Several works are devoted to asymptotic expansions at high conductivity of the electromagnetic field in a domain made of two subdomains {\em when the interface is smooth}: see \cite{S83,MS84,MS85} for plane interface and eddy current approximation, and \cite{HJN08,CDP09,DFP09,Pe09} for a three-dimensional model of skin effect in electromagnetism.

In this section we recall the results presented in \cite{DFP09} on the behavior of the skin depth function. This relies on asymptotic expansions for the electromagnetic field analyzed in the PhD thesis \cite{Pe09}. For the sake of completeness we give in Appendix \ref{AppA} elements of proofs.

\begin{figure}[ht]
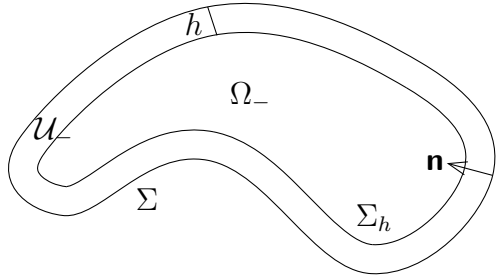

\begin{center}
\def\EpsTube{-10}
\figinit{1.1pt}
\pssetupdate{yes}
\figpt 8:(-48,0)
\figpt 22:(-92,15.5)\figpt 23:(-27.2,60.8)
\figpt 24:(46.3,41.5)\figpt 25:(70.8,2.8)
\figpt 26:(29,-31)\figpt 27:(-27.2,8)\figpt 28:(-77,-11)
\figpt 2:(-24,51)
\figpt 3:(50,50)
\figpt 16:(-70,20)\figpt 15:(0,30)
\figpt 17:(40,10) \figpt 18:(-60,-20)
\figpt 21:(50,65) \figpt 31:(30,-18)  \figpt 34:(-32,55) \figpt 36:(-200,30)
 \figpt 35:(-80,18) 
\figpt 19:(61,65) \figpt 20:(3,27)  \figpt 30:(42.6,33)  \figpt 32:(55,7) \figpt 33:(61.6,6.1) 
\psbeginfig{}
\pssetfillmode{no}\pssetgray{0}
\pscurve[22,23,24,25,26,27,28,22,23,24]
\figptscontrolcurve 40,\NbC[22,23,24,25,26,27,28,22,23,24]
\psEpsLayer \EpsTube,\NbC[40,41,42,43,44,45,46,47,48,49,50,51,52,53,54,55,56,57,58,59,60,61]
\psline[23,2]
\psarrow[25,32]
\pssetfillmode{yes}\pssetgray{0.5}
\psendfig
\figvisu{\figBoxA}{}
{
\figwritew 15: $\Omega\con$(6pt)
\figwritec [34]{$h$}
\figwrites 8: $\Sigma$ (2pt)
\figwriten 31: $\Sigma_{h}$ (2pt)
\figwritew 32: $\nn$(1pt)
\figwritec [35]{$\cU_{-}$}
\figsetmark{$\figBullet$}
}
\centerline{\box\figBoxA}
 \caption{A tubular neighbourhood of the surface $\Sigma$}
\label{Fig1}	
\end{center}
\end{figure}

Here, we assume that $\jj$ is smooth and that Assumptions \ref{H1} and \ref{H2} hold. 
Let $\cU_-$ be a tubular neighborhood of the surface $\Sigma$ in the conductor part $\Omega_-$, see Figure~\ref{Fig1}. We denote by $(y_{\alpha},y_{3})$ a local {\em normal coordinate system} to the surface $\Sigma$ in  $\cU_-$: Here, $y_\alpha$, $\alpha=1,2$, are tangential coordinates on $\Sigma$ and $y_3$ is the normal coordinate to $\Sigma$, cf.\ \cite{Na63,Fa02}.
\input contribF4T.tex

In view of our numerical computations, we concentrate on the magnetic field $\HH_{(\delta)}$ solution of \eqref{FVH0}. It is denoted by $\HH^+_{(\delta)}$ in the dielectric part $\Omega_{+}$, and by $\HH^-_{(\delta)}$ in the conducting part $\Omega_{-}$. Both parts  exhibit series expansions in powers of $\delta$: 
\begin{gather}
\label{6E4a}
     \HH^+_{(\delta)}(\xx) = \HH^+_0(\xx) + \delta\HH^+_{1}(\xx)
     + \mathcal{O}(\delta^2)\, ,\quad \xx\in\Omega_+
\\
\label{6E4b}
   \HH^-_{(\delta)}(\xx) = \HH^-_0 (\xx;\delta)  + \delta\HH^-_{1} (\xx;\delta) 
   + \mathcal{O}(\delta^{2})\, ,\quad\xx\in\Omega_-
  \\
  \label{6E4c}
   \quad\mbox{with}\quad   
   \HH^-_j(\xx;\delta) = \chi(y_3) \,\VV_j(y_\beta,\frac{y_3}{\delta})\, .
\end{gather}
In \eqref{6E4a}-\eqref{6E4b}, the symbol $\mathcal{O}(\delta^2)$ means that the remainder is uniformly bounded by $\delta^2$, and in \eqref{6E4c}, the function $\yy\mapsto\chi(y_3)$ is a smooth cut-off with support in $\overline\cU_-$ and equal to $1$ in a smaller tubular neighborhood of $\Sigma$. The vector fields 
$\VV_{j}: (y_\alpha,Y_3)\mapsto\VV_{j}(y_\alpha,Y_3)$ are \textit{profiles} defined on $\Sigma\times\R^+$: They are exponentially decreasing with respect to $Y_{3}$ and are smooth in all variables.

\subsection{First terms of asymptotics in the conductor region}
Hereafter, we present the construction of the first profiles $\VV_{j}=( {\cV^{\alpha}_{j}},\sv_{j})$ and of the first terms $\HH^+_j$. The normal component $\sv_0$ of the first profile in the conductor is zero:
\begin{equation*}
\sv_{0}=0\, .
\end{equation*}
Then, the first term of the magnetic field in the dielectric region solves Maxwell equations with perfectly conducting conditions on $\Sigma$:
\begin{equation}
\label{H0p}
 \left\{
   \begin{array}{lll}
    \rot\rot \HH^+_0 - \kappa^2\HH^+_0 = \rot\jj  \quad&\mbox{in}\quad \Omega_{+}
\\[0.5ex]
   \HH^+_0\cdot\nn=0  \quad\mbox{and}\quad\!
   \rot\HH^+_0\times\nn=0 \quad &\mbox{on}\quad \Sigma
\\[0.5ex]
   \HH^+_0\times\nn=0 \quad\mbox{and}\quad\!\Div\HH^+_0=0   \quad &\mbox{on}\quad \Gamma.
   \end{array}
    \right.
\end{equation}
Thus the trace $\bsH_0$ of $\HH^+_0$ on the interface $\Sigma$ is \emph{tangential}.

The first profile in the conductor region is exponential with the complex rate $\lambda$ such that $\lambda^2=-i\kappa^2$:
\begin{equation}
\label{V0cd}
 \quad \VV_{0}(y_{\beta},Y_{3}) =
 \bsH_{0} (y_{\beta})\,\mathrm{e}^{-\lambda Y_{3}} 
 \quad\mbox{with}\quad
 \lambda=\kappa\, \mathrm{e}^{-i\pi/4}\,  .
\end{equation}
Note that, if $\bsH_0$ is not identically $0$, there exists $C_0>1$ independent of $\delta$ such that
\begin{equation}
\label{Eestim}
   C_0^{-1}\sqrt \delta \leqslant
   \|\HH^-_0(\,\cdot\,;\delta) \|_{0,\Omega_{-}}\leqslant C_0\sqrt \delta \ .
\end{equation}    

The next term which is determined in the asymptotics is the normal component $\sv_1$ of the profile $\VV_1$:
\begin{equation}
\label{v1cda}
\sv_{1}(y_{\beta},Y_{3}) =\lambda^{-1} D_{\alpha}\sH^{\alpha}_{0}(y_{\beta}) \;\mathrm{e}^{-\lambda Y_{3}}\,.
\end{equation}
Here $D_{\alpha}$ is the covariant derivative on $\Sigma$ and we use the summation convention of repeated indices.
The next term in the dielectric region solves:
\begin{equation}
\label{H1p}
 \left\{
   \begin{array}{lll}
    \rot\rot \HH^+_1 - \kappa^2\HH^+_1 = 0  \quad&\mbox{in}\quad \Omega_{+}
\\[0.5ex]
   \HH^+_1\cdot\nn= \sv_1
    \quad\mbox{and}\quad\!
   \rot\HH^+_1\times\nn= i \lambda \bsH_0 \quad &\mbox{on}\quad \Sigma
\\[0.5ex]
   \HH^+_1\times\nn=0 \quad\mbox{and}\quad\!\Div\HH^+_1=0   \quad &\mbox{on}\quad \Gamma.
   \end{array}
    \right.
\end{equation}
Like above, $\bsH_1$ is the trace of $\HH^+_1$ on the interface $\Sigma$, and $\sH^\alpha_1$ denote its tangential components.
The tangential components ${\cV}^{\alpha}_1$ of the profile $\VV_{1}$ are given by
\begin{equation}
\label{V1cd}
   {\cV^{\alpha}_{1}}(y_{\beta},Y_{3}) =
   \Big[ \sH^{\alpha}_{1}+Y_{3} \big(\cH\,\sH^{\alpha}_{0}
   +b_{\sigma}^{\alpha} \sH^{\sigma}_{0}\big)\Big](y_{\beta}) \;\mathrm{e}^{-\lambda Y_{3}},
   \quad\alpha=1,2\,. 
\end{equation}    
Here $b_{\alpha}^{\gamma}=a^{\gamma\beta}b_{\beta\alpha}$, where $a^{\gamma\beta}$ is the inverse of the metric tensor $a_{\alpha\beta}$ in $\Sigma$, and $b_{\alpha\beta}$ is the curvature tensor in $\Sigma$ and
\[
   \cH=\tfrac12\, b_{\alpha}^{\alpha}
\]
is the \textit{mean curvature} of the surface $\Sigma$. In particular, the sign of $\cH$ depends on the orientation of the surface $\Sigma$. As a convention, the unit normal vector $\nn$ on the surface $\Sigma$ is inwardly oriented to $\Omega_{-}$, see Figure~\ref{Fig1}.

\subsection{Asymptotic behavior of the skin depth}
 In a one-dimensional model, when the conductor $\Omega_{-}$ is a half-space, the classical {\em skin depth} parameter is given by
\begin{equation}
 \label{El}
 \ell(\sigma)=\sqrt{\frac{2}{\omega\mu_{0}\sigma}}\ .
\end{equation}
This length corresponds to the distance from the surface of the conductor where the field has decreased of a rate $\mathrm{e}$. In our situation, following \cite[def.\,4.1]{DFP09}, we extend this definition to curved interfaces.
For a data $\jj$, let us define  
$$
\VV_{(\delta)}(y_\alpha,y_3) := \HH^-_{(\delta)}(\xx),\quad 
y_\alpha\in\Sigma,\quad 0\leqslant y_3<h_0\ ,
$$
for $h_{0}$ small enough. Hereafter for any $\ZZ=(z_1,z_2,z_3)\in\C^3$, $|\ZZ|$ denotes the vector-norm $(|z_1|^2+|z_2|^2+|z_3|^2)^{1/2}$ in $\C^3$ and $\langle\cdot\,,\cdot\rangle$ the corresponding hermitian scalar product.

\begin{defn}
\label{DefL}
Let $\Sigma$ be a smooth surface, and $\jj$ a data of problem \eqref{MS}. such that for all $y_{\alpha}$ in $\Sigma$, $\VV_{(\delta)}(y_{\alpha},0) \neq 0$. The skin depth is the length $\cL(\sigma,y_{\alpha})$ defined on $\Sigma$ and taking the smallest positive value such that
\begin{equation}
\label{epaisspeau}
   |\VV_{(\delta)}\big(y_{\alpha},\cL(\sigma,y_{\alpha})\big) |=
   |\VV_{(\delta)}(y_{\alpha},0) | \ {\mathrm{e}}^{-1}\,.
\end{equation}
\end{defn}

Thus the length $\cL(\sigma,y_{\alpha})$ is the distance from the interface where the field has decreased of a fixed rate. It depends on the conductivity $\sigma$ and of each point $y_{\alpha}$ in the interface $\Sigma$. A priori it also depends on the data $\jj$.

As a consequence of \eqref{V0cd} and \eqref{V1cd}, there holds
\begin{equation}
\label{EV}
\left\{ \begin{array}{l}
   |\VV_{(\delta)}(y_{\alpha},y_{3}) |^2 = 
   |\bsH_{0}(y_{\alpha}) |^2 \ \gm(y_{\alpha},y_3;\delta) \ 
   \mathrm{e}^{-2 y_{3} \Re(\lambda) / \delta} \ , 
   \quad\mbox{with}\\[1ex]
   \gm(y_{\alpha},y_3;\delta) :=  \\ \displaystyle \qquad\qquad
   1 + 2 y_{3} \cH(y_{\alpha}) +
   2 \delta \ \frac{\Re \big\langle  \bsH_{0}(y_{\alpha}) , \bsH_{1}(y_{\alpha}) \big\rangle}
   { | \bsH_{0}(y_{\alpha}) |^{2}} 
   + \cO\big((\delta+y_{3})^2\big)
    \, . 
\end{array}\right.
\end{equation}
Relying on this formula, one can exhibit the asymptotic behavior of the skin depth $\cL(\sigma,y_{\alpha})$ for high conductivity $\sigma$, cf.\ \cite[Th.\,4.2]{DFP09}:
\begin{thm}
Let $\Sigma$ be a regular surface with mean curvature $\cH$. 
Recall that $\ell(\sigma)$ is defined by \eqref{El}.
We assume that $\bsH_{0}(y_{\alpha}) \neq 0$. The skin depth $\cL(\sigma,y_{\alpha})$ has the following behavior for high conductivity:
\begin{equation}
\label{daepaisspeau}
   \cL(\sigma,y_{\alpha})=
   \ell(\sigma)\Big(1+\cH(y_{\alpha})\, \ell(\sigma)+ \cO(\sigma^{-1})\Big),
   \quad\sigma\to\infty\,.
\end{equation}
\end{thm}

\begin{rem}
The higher order terms $\cO(\sigma^{-1})$ in equation \eqref{daepaisspeau} do depend on the data $\jj$ of problem \eqref{MS}.
\end{rem}

\section{Axisymmetric domains}
\label{S4}
\begin{figure}[ht]
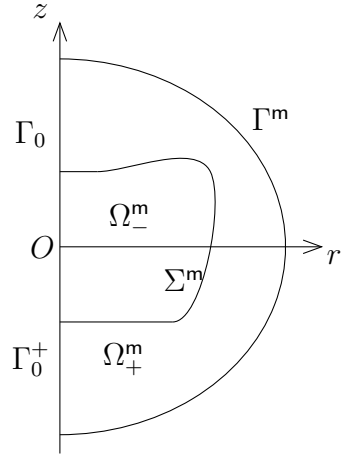

\begin{center}
\figinit{0.5cm}
\def\rx{6} \def\ry{5} 
\figpt 6: $O$ (0,0)
\figpt 7: $r$ (7,0)
\figpt 8: (0,-5.5)
\figpt 9: $z$ (0,6)
\figpt 1: (3,-2)
\figpt 2: (4,2)
\figpt 3: (1,2)
\figpt 4: (0,2)
\figpt 16: (3,2) 
\figpt 17: (0,2)
\figpt 18: (0,-2)
\figpt 10: (0,-8) \figpt 11: (0,8)
\figpt 12: (1.2,0.2) \figpt 13: (1,-3)
\figpt 14: (5,3)
\figpt 15: (4,-0.9)
\figpt 16: (0,3)
\psbeginfig{}

\psarrow[6,7] \psarrow[8,9]
\psset (width=\defaultwidth)
\psline [18,1]
\psline [3,17]
\psset (width=\defaultwidth)
\pscurve[17,1,2,3,17]

\psarcell 6 ; \rx, \ry (-90,90,0)
\psendfig

\figvisu{\figBoxA}{}{
\figwritew 6:$O$(2pt)
\figwritese 7:$r$ (2pt) \figwritenw 9:$z$ (2pt)
\figsetmark{}
\figwritene 12:${\Omega\con^\m}$(2pt)
\figwritee 13:$\Omega\iso^\m$(2pt)
\figwritew 13:$\Gamma^+_0$(18pt)
\figwritene 14:$\Gamma^\m$ (2pt)
\figwritew 16:$\Gamma_{0}$ (2pt)
\figwritew 15:$\Sigma^\m$(2pt)
}
\centerline{\box\figBoxA}
 \caption{The meridian domain $\Omega^\m=\Omega^\m\con\cup\Omega^\m\iso\cup\Sigma^\m$ with boundary $\partial\Omega^\m=\Gamma^\m\cup\Gamma_{0}$}
\label{Fmer}
\end{center}
\end{figure}

In order to perform scalar two dimensional computations which could represent correctly the features of a three-dimensional problem, we choose to consider an \emph{axisymmetric configuration} in which $\Omega_+$ and $\Omega_-$ are axisymmetric domains with the same axis $\Xi_0$: in \emph{cylindrical coordinates} $(r,\theta,z)$ associated with this axis, there exists bi-dimensional ``meridian'' domains $\Omega^\m$ and $\Omega^\m_\pm$ such that
\begin{gather*}
   \Omega = \{ \xx\in\R^3\ | \ \  (r,z)\in\Omega^\m,\ \theta\in\T\}, \\
   \Omega_\pm = \{ \xx\in\R^3\ | \ \  (r,z)\in\Omega^\m_\pm,\ \theta\in\T\}.
\end{gather*}
Here $\T=\R/(2\pi\Z)$ is the one-dimensional torus. We denote by $\Gamma^\m$ and $\Sigma^\m$ the meridian curves corresponding to $\Gamma$ and $\Sigma$, respectively, and by $\Gamma_0$, $\Gamma^+_0$ the following subsets of the rotation axis $\Xi_0$, see Fig.\ \ref{Fmer}
\[
   \Gamma_0 = \Xi_0\cap\overline\Omega{}^\m \quad\mbox{and}\quad
   \Gamma^+_0 = \Xi_0\cap\overline\Omega{}^\m_+. 
\]

On such an axisymmetric configuration we consider a modification of problem \eqref{MS}-\eqref{PIbc}: We take $\jj\equiv0$ and impose instead of \eqref{PIbc} non-homogeneous magnetic boundary conditions
\begin{equation}
\label{EDbc}
   \EE\cdot\nn=0 \quad\mbox{and}\quad \HH\times\nn=\GG \times\nn \quad\mbox{on}\quad\Gamma\,,
\end{equation}
for a given data $\GG\in\bH(\rot,\Omega):=\{ \uu\in\bL^2(\Omega)\ | \ \rot \uu \in \bL^2(\Omega) \}$.
Then, the variational problem for the magnetic field $\HH_{(\delta)}$ solution of the Maxwell equations \eqref{MS}-\eqref{EDbc} writes
\\[1ex]
{\it Find $\HH_{(\delta)}\in  \bH_0(\rot,\Omega)+\GG$\ , such that for all $\KK\in \bH_0(\rot,\Omega)$,}
\begin{equation}
\int_{\Omega} \big(\frac1{\varepsilon(\delta)}\rot\HH_{(\delta)} \cdot \rot\KK- \kappa^2\HH_{(\delta)}\cdot\KK\big)
\,\dr\xx =0  \ .
\label{FVH}
\end{equation}

\subsection{Formulation in cylindrical components}
For a vector field $\HH=(H_{1},H_{2},H_{3})$, we introduce its \textit{cylindrical components} $(H_{r},H_{\theta},H_{z})$ according to
$$
\left\{ \begin{array}{rcl}
   H_{r}(r,\theta,z) &=& H_{1}(\xx)\cos\theta+H_{2}(\xx) \sin\theta\ ,\\
   H_{\theta}(r,\theta,z) &=& -H_{1}(\xx)\sin\theta+H_{2}(\xx) \cos\theta\ ,\\
   H_{z}(r,\theta,z) &=& H_{3}(\xx)\ ,
\end{array}\right.
$$
and we set
\[
   \brH(r,\theta,z) = \big(H_{r}(r,\theta,z),H_{\theta}(r,\theta,z),H_{z}(r,\theta,z)\big).
\]
We say that $\HH$ is axisymmetric if $\brH$ does not depend on the angular variable $\theta$.

The Maxwell problem \eqref{FVH} is axisymmetric, which means that, expressed in cylindrical variables $(r,\theta,z)$ and components $\brH$, its coefficients do not depend on $\theta$ \cite{HiLe05,B-D-M99,Nk05}.
Recall that for a vector field $\HH=(H_{1},H_{2},H_{3})$ the \textit{cylindrical components} of its curl write
\begin{equation}
\label{Ecurl}
\left\{\begin{aligned}
   &\te (\rot\HH)_r = \frac1r\partial_\theta H_z - \partial_z H_\theta\ ,\\
   &(\rot\HH)_\theta = \partial_z H_r - \partial_r H_z\ ,\\
   &\te (\rot\HH)_z =  \frac1r \big( \partial_r(r H_\theta) - \partial_\theta H_r\big)\ .
\end{aligned}\right.
\end{equation}
As a consequence, if the right-hand side $\GG$ is axisymmetric, and if \eqref{FVH} has a unique solution, then this solution is axisymmetric. According to \eqref{Ecurl}, when $\HH$ is axisymmetric its curl reduces to
\begin{equation}
\label{EcurlAxi}
\left\{\begin{aligned}
   &\te (\rot\HH)_r = - \partial_z H_\theta\ ,\\
   &(\rot\HH)_\theta = \partial_z H_r - \partial_r H_z\ ,\\
   &\te (\rot\HH)_z = \frac1r \, \partial_r(r H_\theta) \ .
\end{aligned}\right.
\end{equation}

\subsection{Axisymmetric orthoradial problem}
We say that $\HH$ is orthoradial if its components $H_r$ and $H_z$ are zero.

We assume that $\GG$ is axisymmetric and orthoradial, i.e.,
\[
   \breve\GG(r,\theta,z) = (0,\sG(r,z),0).
\] 
The components $(H_r,H_z)$ and $H_\theta$ being uncoupled in \eqref{EcurlAxi}, and the solution of problem \eqref{FVH} being unique, we obtain that this solution is orthoradial 
$$
   \brH_{(\delta)}(r,\theta,z) = (0,\sH_{(\delta)}(r,z),0).
$$ 

\subsubsection{Variational formulation}
In this framework, the change from Cartesian to cylindrical coordinates requires the modification of the  solution spaces $\bH_0(\rot,\Omega)$ used in problem \eqref{FVH}. Precisely, the weighted space characterizing the \textit{orthoradial component} $\sH_{(\delta)}(r,z)$ is 
$$
\mathrm{V}_{1,\Gamma^\m}^1(\Omega^\m)=\{v\in \mathrm{H}_{1}^1(\Omega^\m) \,|\,  v\in\mathrm{L}_{-1}^2(\Omega^\m)\quad \mbox{and} \quad v=0 \quad \mbox{on}\quad\Gamma^\m\}.
$$
Here,
$$ \mathrm{H}_{1}^1(\Omega^\m)=\{v\in \mathrm{L}_{1}^2(\Omega^\m) \,|\ \ 
\partial_{r}^j\partial_{z}^{1-j} v\in\mathrm{L}_{1}^2(\Omega^\m), \  
j=0,1 \},
$$ 
and for all $\alpha\in\R$, the space $\mathrm{L}_{\alpha}^2(\Omega^\m)$ is the set of measurable functions $v(r,z)$ such that  
\begin{equation*}
\|v \|^2_{\mathrm{L}_{\alpha}^2(\Omega^\m)} = \int_{\Omega^\m}  |v|^2 \,r^\alpha drdz <+\infty\ .
\end{equation*}

\begin{rem}
The space $\mathrm{V}_{1,\Gamma^\m}^1(\Omega^\m)$ incorporates essential boundary conditions, in particular on $\Gamma_{0}$, where $v=0$, see \cite[Remark\, \textrm{II.1.1}]{B-D-M99}. 
\end{rem}

This leads us to solve the following two-dimensional scalar problem set in $\Omega^\m$.
\\[1ex] 
{\it Find $\sH_{(\delta)} \in \mathrm{V}_{1,\Gamma^\m}^1(\Omega^\m)+\sG$ such that for all 
$\ \sw\in  \mathrm{V}_{1,\Gamma^\m}^1(\Omega^\m)$,}
\begin{equation}
a^{\delta}(\sH_{(\delta)} ,\sw) 
= 0\ ,
\label{FVHtheta}
\end{equation}
where 
$$
   a^{\delta}(\sH ,\sw):=
   \int_{\Omega^\m} \frac1{\varepsilon(\delta)}\Big( \partial_{z} \sH \, \partial_{z}  \sw
   +\frac1{r} \partial_{r} (r \sH )\, \frac1{r}\partial_{r}(r  \sw) \Big) \,rdrdz  
   -  \kappa^2 \int_{\Omega^\m}  \sH\,  \sw \,rdrdz \ .
$$

\subsubsection{Asymptotic expansion}
\label{AEh}
 Let $\left(r(\xi),z(\xi)\right)=\XX(\xi)$, $\xi \in (0,L)$, be an \textit{arc-length coordinate} on the interface $\Sigma^\m$. Here $\xi \mapsto \XX(\xi)$ is a $\mathcal{C}^{\infty}$ function, and $L$ is the length of the curve $\Sigma^\m$. Let $(\xi,y_3)$ be the associate normal coordinate system in a tubular neighborhood of $\Sigma^\m$ inside $\Omega^\m_-$. Then the normal vector $\nn(\xi)$ at the point $\XX(\xi)$ can be written as (Frenet frame) 
\begin{equation}
\label{nzp}
   \nn(\xi)=\big( -z'(\xi),r'(\xi)\big)
   \quad\mbox{with}\quad
   z'(\xi) = \frac{d z}{d\xi} \quad\mbox{and}\quad r'(\xi) = \frac{d r}{d\xi}.
\end{equation}
Finally we denote by $k(\xi)$ the curvature of $\Sigma^\m$ in $\XX(\xi)$.

In accordance with \eqref{6E4a}-\eqref{6E4b}-\eqref{6E4c}, we can exhibit series expansions in powers of $\delta$ for the magnetic field $\sH_{(\delta)}$ which we denote by $\sH^+_{(\delta)}$ in the dielectric part $\Omega_{+}$, and by $\sH^-_{(\delta)}$ in the conducting part $\Omega_{-}$:  
\begin{gather*}
   \sH^+_{(\delta)}(r,z) = \sH^+_0(r,z) + \delta\sH^+_1(r,z) + \mathcal{O}(\delta^2)\, ,\\
   \sH^-_{(\delta)}(r,z) = \sH^-_0(r,z;\delta) + \delta\sH^-_1(r,z;\delta) 
   + \mathcal{O}(\delta^2),
   \quad   
   \sH^-_j(r,z;\delta) = \chi(y_3) \,\svV_j(\xi,\frac{y_3}{\delta})\, .
\end{gather*}
Here the profiles $\svV_{j}$ are defined on $\Sigma^\m\times\R^+$. Hereafter, we focus on the first terms $\sH^+_0$, $\svV_{0}$, $\sH^+_1$, and $\svV_{1}$. We introduce the interior and boundary operators 
\begin{multline}
\label{EOpDB}
\D(r,z;\partial_r,\partial_z) = \partial^2_r + \frac1{r} \partial_r  + \partial^2_z -\frac1{r^2} \quad \\ \mbox{and} \quad 
\B(\xi;\partial_r,\partial_z) = -z'(\xi)(\partial_r + \frac1{r}) + r'(\xi)\partial_z \, .
\end{multline}
In appendix \ref{AppB1}, we give the expansion of these operators in power series of $\delta$ inside the domain $\Omega_{-}$ and on the interface $\Sigma^\m$. The terms $\sH^+_0$, $\svV_{0}$, $\sH^+_1$, and $\svV_{1}$ satisfy the following problems coupled by their boundary conditions on the interface $\Sigma^\m$ (corresponding to $Y=0$) -- compare with \eqref{H0p}-\eqref{V0cd}, \eqref{H1p}-\eqref{V1cd}:
\begin{equation}
\label{Eh0+}
\left\{
    \begin{array}{lllll}
 \D \sH^+_0  +\kappa ^2  \sH^+_0 &=& 0 & \mbox{in}\quad \Omega^\m_{+} \ ,  
\\[0.5ex] 
 \B \sH^+_0 &=& 0 & \mbox{on}\quad \Sigma^\m \ ,
\\[0.5ex]
 \sH^+_0 &=&  \sG  & \mbox{on}\quad \Gamma^\m\cup \Gamma_{0}^+ \ ,
    \end{array}
\right.
\end{equation}
\begin{equation}
\label{Ev0}
\left\{
   \begin{array}{lllll}
 (\partial^2_{Y} -\lambda^2)\svV_0
 &=& 0 & \mbox{for}\quad 0<Y<+\infty \ ,
\\[0.5ex]
 \svV_0& =&  \sH_0^{+} & \mbox{for}\quad Y=0\ ,
   \end{array}
\right.
\end{equation}
with $\lambda$ defined by \eqref{V0cd}, from which we deduce:
\begin{equation}
\label{v0}
\svV_{0}(\xi,Y)=\mathrm{e}^{-\lambda Y} \sH^+_0\big(\XX(\xi)\big)\ .
\end{equation}
The next problem in the dielectric part is
\begin{equation}
\label{Eh1+}
\left\{
   \begin{array}{lllll}
 \D \sH^{+}_1 + \kappa ^2  \sH^{+}_1 &=& 0& \mbox{in}\quad\Omega^\m_{+}\ ,
\\[0.5ex]
\B \sH^{+}_1&=&-i  \partial_{Y}\svV_{0}  & \mbox{on}\quad \Sigma^\m\ ,
\\[0.5ex]
 \sH^{+}_1& =& 0 & \mbox{on}\quad \Gamma^\m\cup \Gamma_{0}^+\ ,
   \end{array}
\right.
\end{equation}
Thus \eqref{v0} yields that $\sH^+_1$ satisfies  $\B \sH^+_1=i \lambda  \sH^+_0\big|_{\Sigma^\m} $ on the interface.

The next problem in the conductor part is
\begin{equation}
\label{Ev1}
\left\{
   \begin{array}{lllll}
 (\partial^2_{Y} - \lambda^2)\svV_1 &=& -\A_{1} \svV_{0} &  \mbox{for}\quad 0<Y<+\infty \ ,
\\[0.5ex]
 \svV_1& =&  \sH_1^{+} &  \mbox{for}\quad Y=0\ .
   \end{array}
\right.
\end{equation}
 Here, $\A_{1}\svV_{0}= - \big(k+\frac{z'}{r}\big)(\xi) \partial_{Y} \svV_0$. From \eqref{v0}, we infer $\A_{1}\svV_{0}=\mathrm{e}^{-\lambda Y} \lambda  \big(k+\frac{z'}{r}\big) (\xi) \sH^+_0\big|_{\Sigma^\m}$. Then, from equation \eqref{Ev1}, we obtain 
\begin{equation}
\label{v1}
\svV_{1}(\xi,Y)=\mathrm{e}^{-\lambda Y} \Big[\sH^+_1(\XX(\xi))
+ \frac{Y}{2}  \Big(k+\frac{z'}{r}\Big)(\xi)\,\sH^+_0(\XX(\xi))
\Big]\ .
\end{equation}
Note that we can also deduce this profile $\svV_{1}$ from equations \eqref{V0cd} and  \eqref{V1cd}, but this is not obvious because the cylindrical coordinates are not a normal coordinate system, see Appendix \ref{AppB2}.

\begin{rem}
\label{Disc}
Subsequently, we assume that the data $\sG$ is a real valued function. Thus, the right hand side of the boundary value problem \eqref{Eh0+} is real. Hence $\sH^+_0$ is a real valued function. 
Recall that $\B \sH^+_1= \kappa\, \mathrm{e}^{i\pi/4} \sH^+_0\big|_{\Sigma}$. From the boundary value problem \eqref{Eh1+}, we infer: $\Re \sH^+_1=\Im \sH^+_1$. We will exploit this relationship in the numerical simulations of skin effect, see \S 6. 
\end{rem}

\subsection{Configurations chosen for computations}
We consider three classes of geometric configurations: one cylindrical configuration (A) and two spheroidal configurations (B and C).

\subsubsection*{Configuration A: cylindrical geometry}
\begin{figure}[ht]
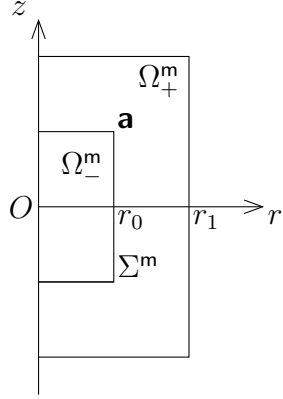

\begin{center}
\figinit{0.5cm}
\figpt 6: $O$ (0,0)
\figpt 7: $r$ (6,0)
\figpt 8: (0,-5)
\figpt 9: $z$ (0,5)
\figpt 1:(0,2)
\figpt 2:(2,2)
\figpt 3:(2,-2)
\figpt 4:(0,-2)
\figpt 5:(2,0)
\figpt 20:(4,0)
\figpt 16:(0,4)
\figpt 17:(4,4)
\figpt 18:(4,-4)
\figpt 19:(0,-4)
\figpt 10: (0,-8) \figpt 11: (0,8)
\figpt 12: (0.5,0.5) \figpt 13: (3,-2.5)
\figpt 14: (3,2)
\figpt 15: (2,-1.2)
\psbeginfig{}

\psarrow[6,7] \psarrow[8,9]
\psset (width=\defaultwidth)
\psline [1,2,3,4]
\psline [3,4]
\psset (width=\defaultwidth)
\psline[16,17,18,19]
\psendfig

\figvisu{\figBoxA}{}{
\figwritew 6:$O$(2pt)
\figwritese 5:$r_{0}$ (2pt)
\figwritese 20:$r_{1}$ (2pt)
\figwritese 7:$r$ (2pt) \figwritenw 9:$z$ (2pt)
\figsetmark{}
\figwritene 12:$\Omega\con^\m$(2pt)
\figwritesw 17:$\Omega\iso^\m$(5pt)
\figwritese 15:$\Sigma^\m$(2pt)
\figwritene 2:$\aa$(2pt)
}
\centerline{\box\figBoxA}
 \caption{The meridian domain $\Omega^\m$ in configuration A}
\label{FCylGeo}
\end{center}
\end{figure}

We assume that $\Omega$ is a circular cylinder of radius $r_{1}$ and length $\ell_{1}$, and $\Omega\con$ is a coaxial cylinder of radius $r_{0}$ and  length $\ell_{0}$. Hence, $\Omega^\m$ is a rectangle of width $r_{1}$ and length $\ell_{1}$, and $\Omega^\m\con$ is a coaxial rectangle of width $r_{0}$ and length $\ell_{0}$, see Figure~\ref{FCylGeo}.
We choose the parameters $r_{0}=1$, $\ell_{0}=2$, $r_{1}=2$, $\ell_{1}=4$ in computations.

\subsubsection*{Configuration B: spheroidal geometry}
\begin{figure}[ht]
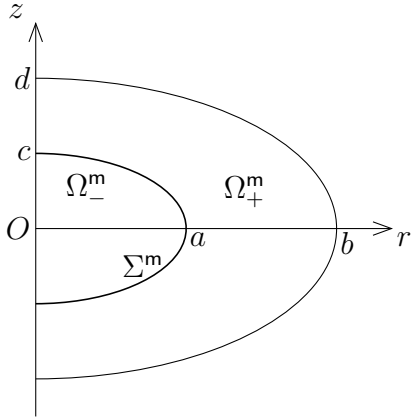

\begin{center}
\figinit{0.5cm}
\def\rx{4} \def\ry{2} \def\Rx{8} \def\Ry{4}
\figpt 6: $O$ (0,0)
\figpt 7: $r$ (9.5,0)
\figpt 8: (0,-5)
\figpt 9: $z$ (0,5.5)
\figpt 1:(3,-2)
\figpt 2:(4,2)
\figpt 3:(2,3)
\figpt 4: (0,3)
\figptorthoprojline 5:$A_5$=1/6,4/
\figpt 16:(4,0) 
\figpt 17:(8,0) 
\figpt 20:(0,2)
\figpt 21:(0,4)
\figpt 10: (0,-8) \figpt 11: (0,8)
\figpt 12: (0.7,0.5) \figpt 13: (5,1.4)
\figpt 15: (2.2,-0.6)
\psbeginfig{}
\psarrow[6,7] \psarrow[8,9]
\psset (width=0.6)
\psarcell 6 ; \rx, \ry (-90, 90,0) 
\psset (width=\defaultwidth)
\psarcell 6 ; \Rx, \Ry (-90,90,0)
\psendfig
\figvisu{\figBoxA}{}{
\figwritew 6:$O$(2pt)
\figwritew 20:$c$(2pt)
\figwritew 21:$d$(2pt)
\figwritese 7:$r$ (2pt) \figwritenw 9:$z$ (2pt) \figwritese 16:$a$ (2pt)
 \figwritese 17:$b$ (2pt)
\figsetmark{}
\figwritene 12:$\Omega\con^\m$(2pt)
\figwritese 13:$\Omega\iso^\m$(0pt)
\figwritese 15:$\Sigma^\m$(2pt)
}
\centerline{\box\figBoxA}
 \caption{The meridian domain $\Omega^\m$ in configuration B1}
\label{FEllGeo}
\end{center}
\end{figure}

We assume that $\Omega$ is a spheroid, and $\Omega\con$ is a coaxial spheroid. We denote by  $a$ and $c$, respectively $b$ and $d$, the semimajor and semiminor axis of $\Omega\con$, respectively $\Omega$. Hence, $\Omega^\m$ is a semi-ellipse with semimajor axis $b$ and semiminor axis $d$, and $\Omega^\m\con$ is a coaxial semi-ellipse with axis lengths $a$ and $c$, see Figure~\ref{FEllGeo}.
In computations, we consider the following parameters 
\begin{equation}
\label{EConfB}
   \begin{array}{ll}
\mbox{Configuration B1:} \quad a=2,b=4,c=1,d=2\ ,  \ \ 
\mbox{(oblate spheroid)}
\\
\mbox{Configuration B2:} \quad a=4,b=8,c=1,d=2\ ,   \ \ 
\mbox{(more oblate spheroid)}
.\\
  \end{array}
\end{equation}
Note that all the numerical parameters are given in SI units, in particular, $a,b,c$, and $d$ are defined in meters, and the conductivity $\sigma$ is defined in siemens per meter $S.m^{-1}$.

\subsubsection*{Configuration C: spheroidal geometry}
\begin{figure}[ht]
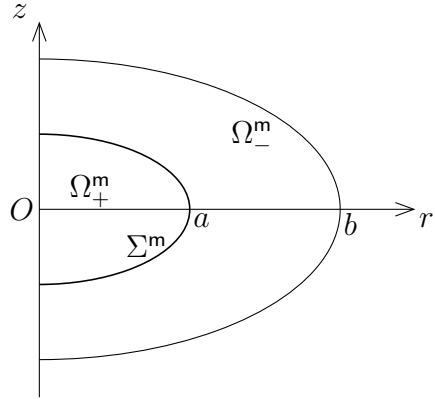

\figinit{0.5cm}
\def\rx{4} \def\ry{2} \def\Rx{8} \def\Ry{4}
\figpt 6: $O$ (0,0)
\figpt 7: $r$ (10,0)
\figpt 8: (0,-5)
\figpt 9: $z$ (0,5)
\figpt 1:(3,-2)
\figpt 2:(4,2)
\figpt 3:(2,3)
\figpt 4: (0,3)
\figptorthoprojline 5:$A_5$=1/6,4/
\figpt 16:(4,0) 
\figpt 17:(8,0) 
\figpt 10: (0,-8) \figpt 11: (0,8)
\figpt 12: (0.7,1) \figpt 13: (5,1.4)
\figpt 15: (2.2,-0.6)
\psbeginfig{}
\psarrow[6,7] \psarrow[8,9]
\psset (width=0.6)
\psarcell 6 ; \rx, \ry (-90, 90,0) 
\psset (width=\defaultwidth)
\psarcell 6 ; \Rx, \Ry (-90,90,0)
\psendfig
\figvisu{\figBoxA}{}{
\figwritew 6:$O$(2pt)
\figwritese 7:$r$ (2pt) \figwritenw 9:$z$ (2pt) \figwritese 16:$a$ (2pt)
 \figwritese 17:$b$ (2pt)
\figsetmark{}
\figwritene 13:$\Omega\con^\m$(2pt)
\figwritese 12:$\Omega\iso^\m$(2pt)
\figwritese 15:$\Sigma^\m$(2pt)
}
\centerline{\box\figBoxA}
 \caption{The meridian domain $\Omega^\m$ in configuration C1}
\label{FEllGeoC}
\end{figure}

We introduce a configuration C, switching the roles of the subdomains $\Omega\con^\m$ and $\Omega\iso^\m$ in the configuration B, see Figure~\ref{FEllGeoC}.
We denote by C1 and C2 the configuration C, which correspond to the choice of parameters $a,b,c,d$, see \eqref{EConfB}.

\subsubsection*{Right hand sides of problems}
In configurations A and B, we take the data $\sG=r$. Hence, $\sH_{(\delta)}$ satisfies the following inhomogeneous Dirichlet boundary condition 
$$
\sH_{(\delta)}(r,z)=r  \quad \mbox{on}   \quad  \Gamma^\m \ .
$$ 
For configuration C1, we take $\sG=0$ and an interior data $f:=(\rot\jj)_\theta$ with support inside $\Omega\iso^\m$: $f=10^2$ if $r^2 /4 +z^2 \leqslant 0.8 $ and $f=0$ otherwise.

\section{Finite element discretizations and computations}
\label{S5}

\subsection{Finite element method}
In this section, we consider two benchmarks for the computational domain: configuration A, see Figure~\ref{FCylGeo}, and configuration B, see Figure~\ref{FEllGeo}. We use high order elements available in the finite element library M\'{e}lina, see \cite{melina}, and quadrangular meshes in the meridian domain. We discretize the variational problem \eqref{FVHtheta}. The script used to solve the problem \eqref{FVHtheta} is adapted from \cite{Bo05}. In the computations, we fix the angular frequency $\omega=3.10^7$. We denote by $\sH_{(\delta)}^{p,\MM}$ the computed solution of the discretized problem \eqref{FVHtheta} with an interpolation degree $p$ and a mesh $\MM$. We define
$$
A_{\sigma}^{p,\MM}:=\|\sH_{(\delta)}^{p,\MM} \|_{\mathrm{L}_{1}^2(\Omega\con^\m)} \quad\mbox{with}\quad
\sigma = \omega\varepsilon_0\delta^{-2},\ \ \mbox{cf. \eqref{delta}}\ .
$$

\subsection{Interpolation degree}
We first check the convergence when the interpolation degree of the finite elements increases.

\subsubsection{Configuration B1}
\begin{figure}[ht]
\begin{center}
\includegraphics[keepaspectratio=true,width=4.cm]{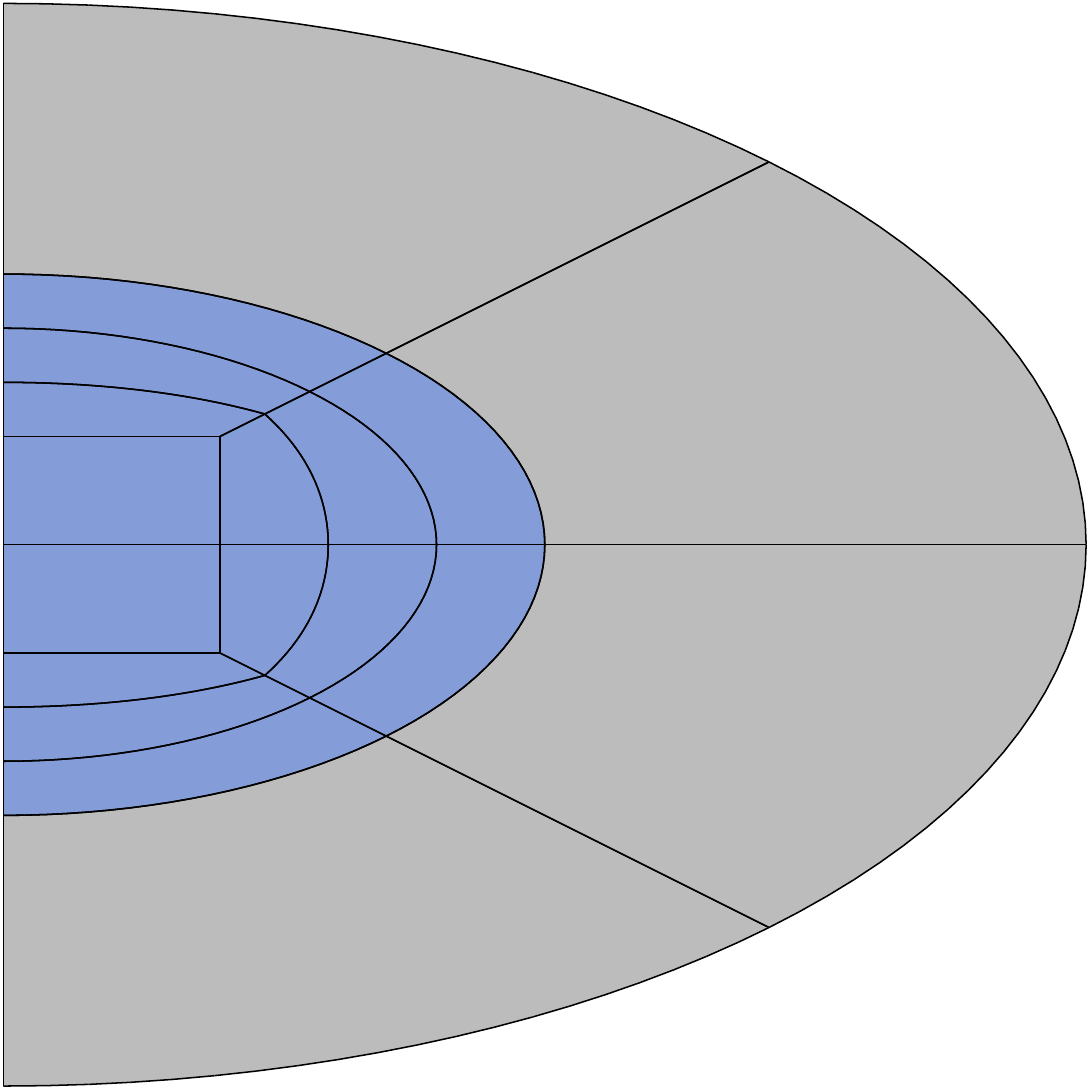} \quad
\includegraphics[keepaspectratio=true,width=4.cm]{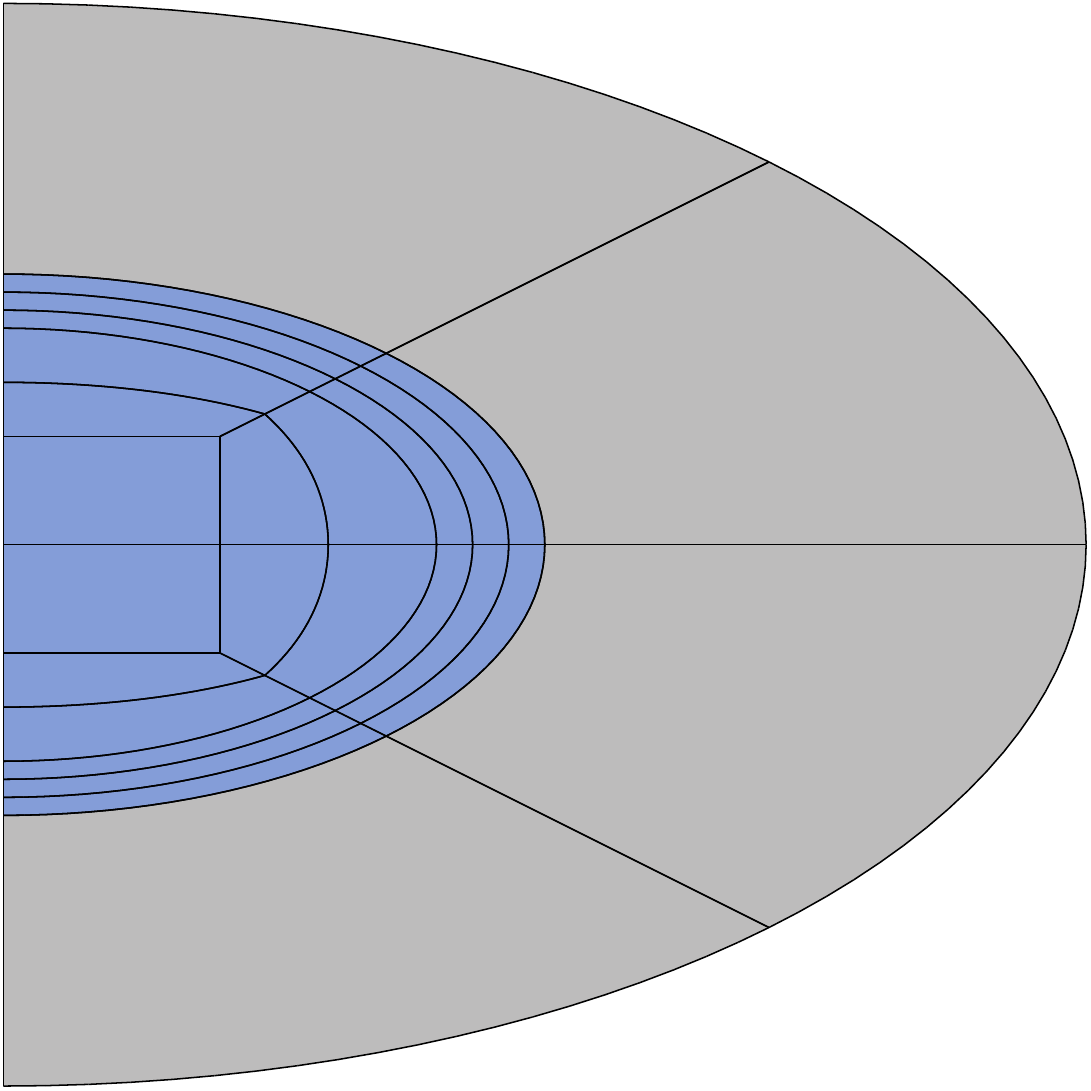} \quad
\includegraphics[keepaspectratio=true,width=4.cm]{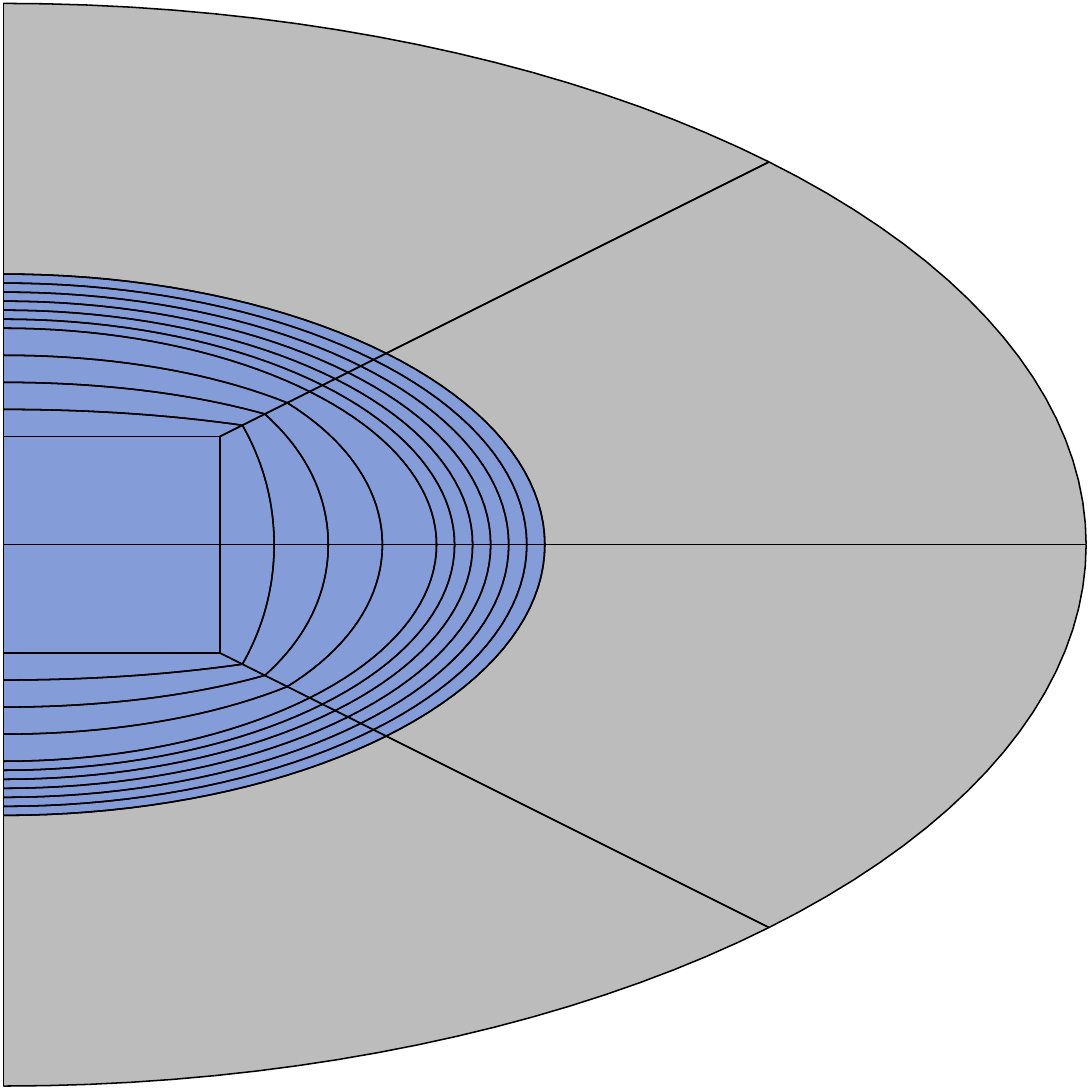}
 \caption{The meshes $\MM_1$, $\MM_3$ and $\MM_6$ for configuration B1}
\label{F0}
\end{center} 
\end{figure}

We consider the discretized problem with different degrees: $Q_p$\footnote{ Recall that $Q_p$ is the vector space of polynomials of two variables and partial degree $p$ defined on the reference element $\widehat{K}:=[0,1]\times [0,1]$}, for all $p=1,\cdots, 20$, and with three different meshes $\MM_1$, $\MM_3$ and $\MM_{6}$ with $1$, $3$ or $6$ layers of elongated elements in the skin region of the conductor $\Omega\con^\m$, see Figure~\ref{F0}.

\begin{figure}[ht]
\begin{center}
\includegraphics[keepaspectratio=true,width=7cm]{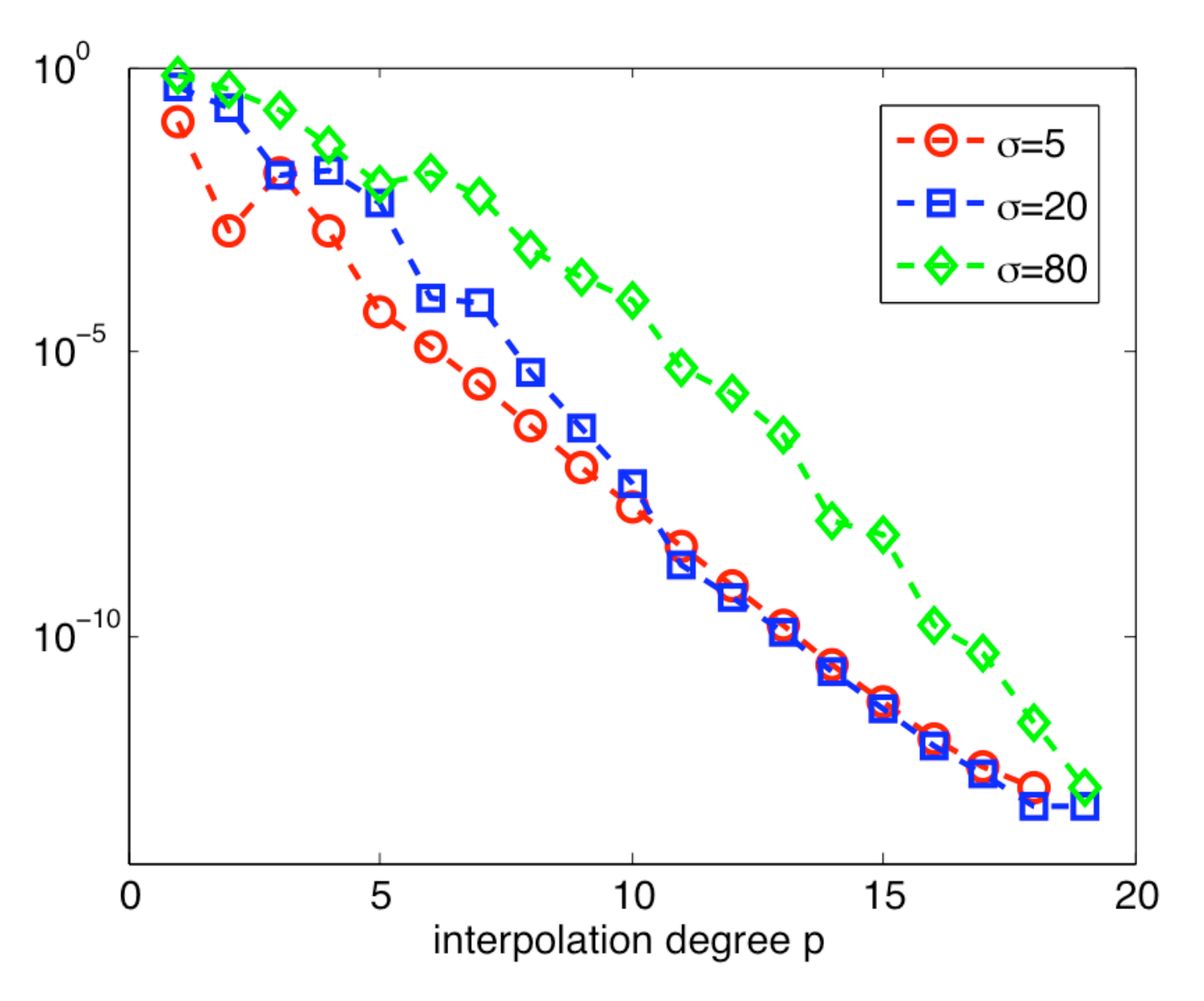}
\caption{Graph of $\big|A_{\sigma}^{p,\MM_1}-A_{\sigma}^{20,\MM_6}\big|$ with respect to $p=1,\cdots,19$ in semi-log coordinates, for $\sigma\in\{5,20,80\}$ and configuration B1}
\label{NormL2a2K1p20K6}
\end{center} 
\end{figure}

\begin{figure}[ht]
\begin{center}
\includegraphics[keepaspectratio=true,width=7cm]{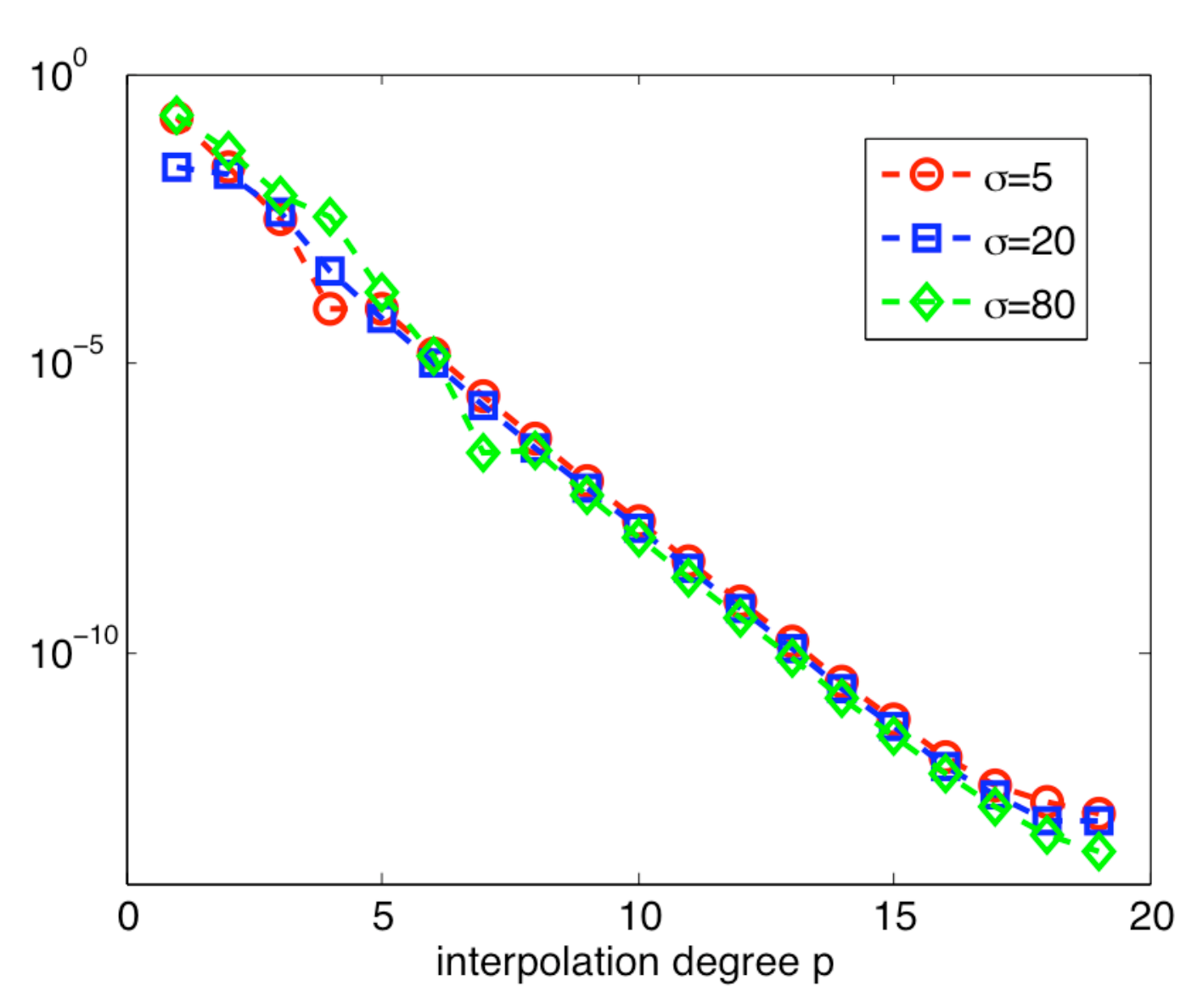}
\caption{Graph of $\big|A_{\sigma}^{p,\MM_3}-A_{\sigma}^{20,\MM_6}\big|$ with respect to $p=1,\cdots,19$ in semi-log coordinates, for $\sigma\in\{5,20,80\}$ and configuration B1}
\label{NormL2a2K3p20K6}
\end{center} 
\end{figure}

We represent  in Figure~\ref{NormL2a2K1p20K6} the absolute value of the difference between the weighted norms $A_{\sigma}^{p,\MM_1}$ and $A_{\sigma}^{20,\MM_6}$, versus $p$ in semilogarithmic coordinates, and in each case: $\sigma=5$ with circles, $\sigma=20$ with squares, and $\sigma=80$ with diamonds.  We refer to \cite{ScSu96} for theoretical results of convergence for the $p$-version in presence of an exponentially decreasing boundary layer. In  Figure~\ref{NormL2a2K3p20K6} we use mesh $\MM_3$ instead of $\MM_1$.

\begin{figure}[ht]
\begin{center}
\includegraphics[keepaspectratio=true,width=7cm]{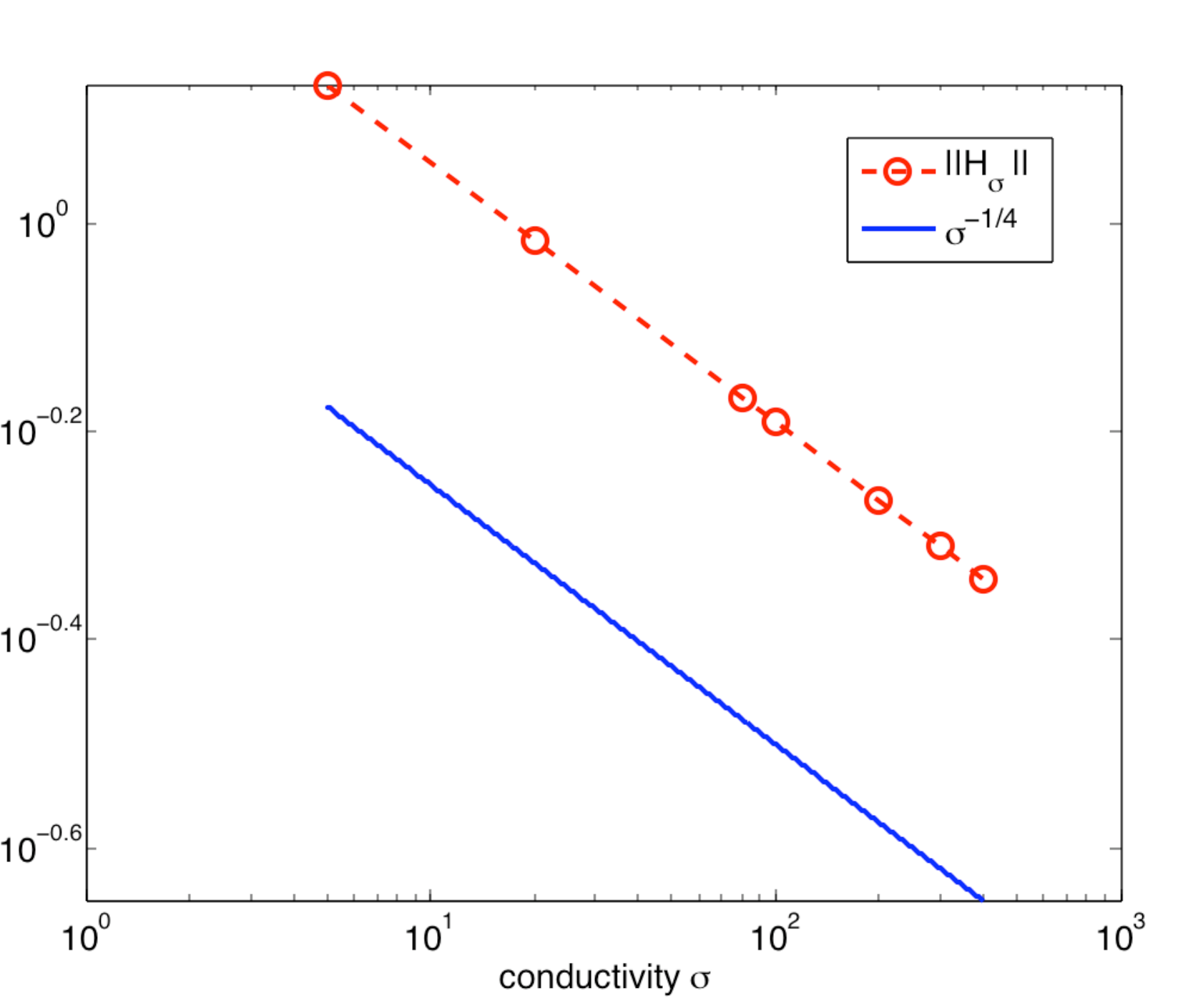}
  \caption{In circle: $A_{\sigma}^{16,\MM_3}$ for $\sigma=5,20,80,100,200,300,400$. In solid line: the  graph of the function $\sigma \mapsto \sigma^{-1/4}$ in log-log coordinates for configuration B1}
\label{FnormL2bis}
\end{center} 
\end{figure}

In Figure~\ref{FnormL2bis}, we plot in log-log coordinates the weighted norm $A_{\sigma}^{16,\MM_3}$ with respect to $\sigma=5,20,80,100,200,300,400$  with circles, and the graph of $\sigma \mapsto \sigma^{-1/4}$ by a solid line. The figure shows that $A_{\sigma}^{16,\MM_3}$ behaves like $\sigma^{-1/4}$ when $\sigma\to\infty$. This behavior is consistent with the asymptotic expansion \eqref{6E4b} and the estimate \eqref{Eestim}, (recall formula \eqref{delta}).

\clearpage

\subsubsection{Configuration A}
\begin{figure}[ht]
\begin{center}
\includegraphics[keepaspectratio=true,width=2.cm]{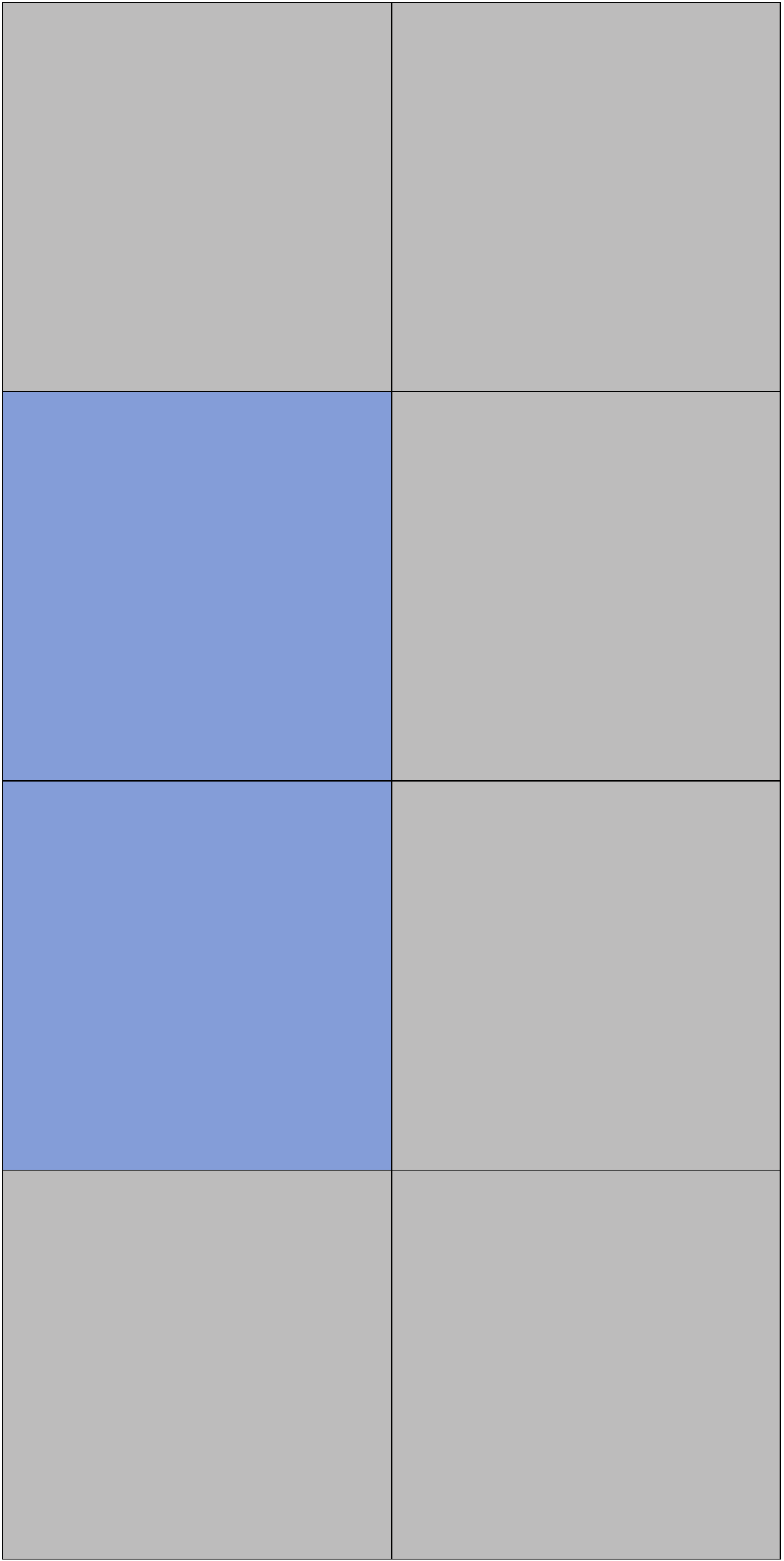}  \quad
\includegraphics[keepaspectratio=true,width=2.cm]{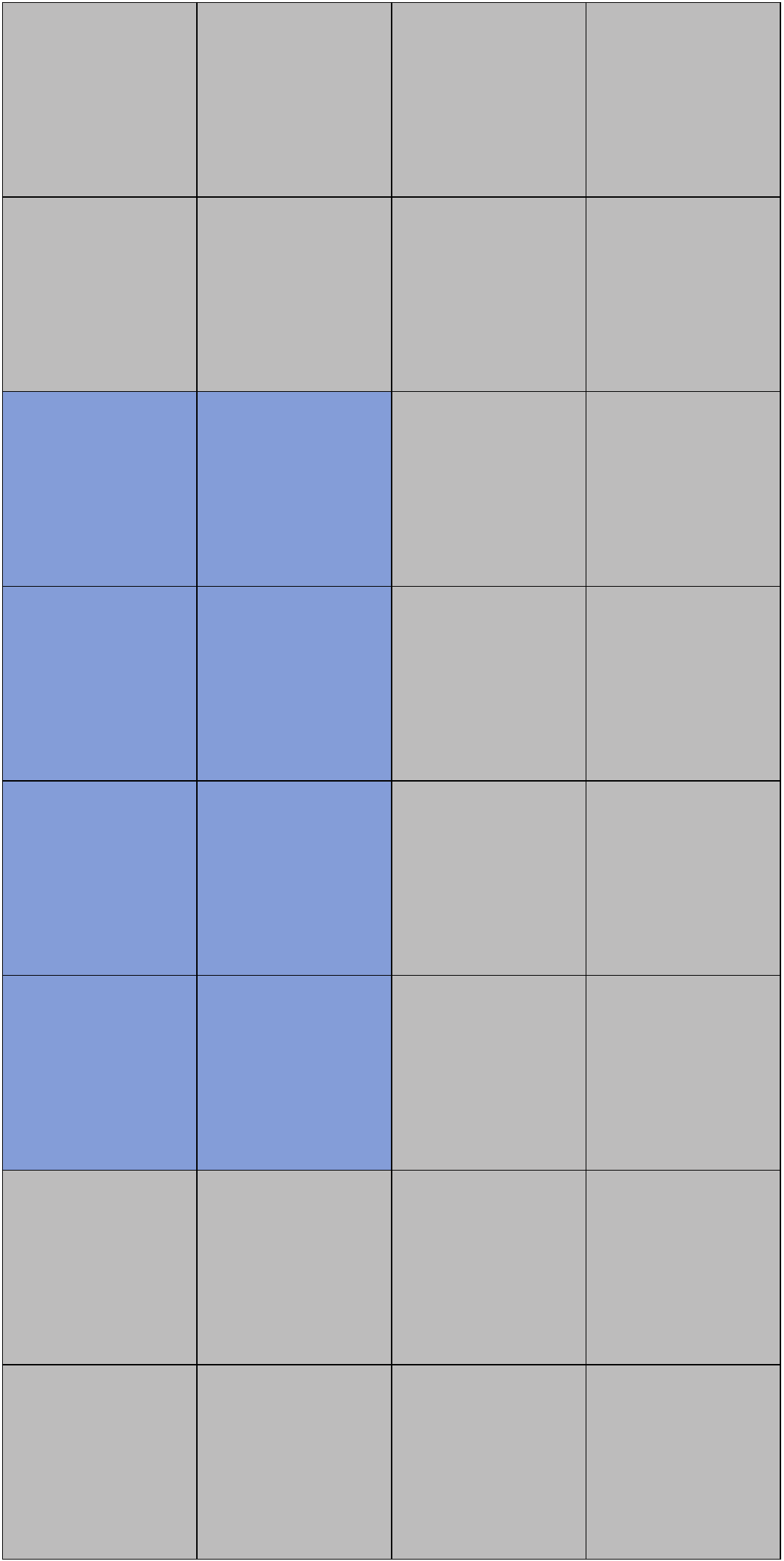} \quad
\includegraphics[keepaspectratio=true,width=2.cm]{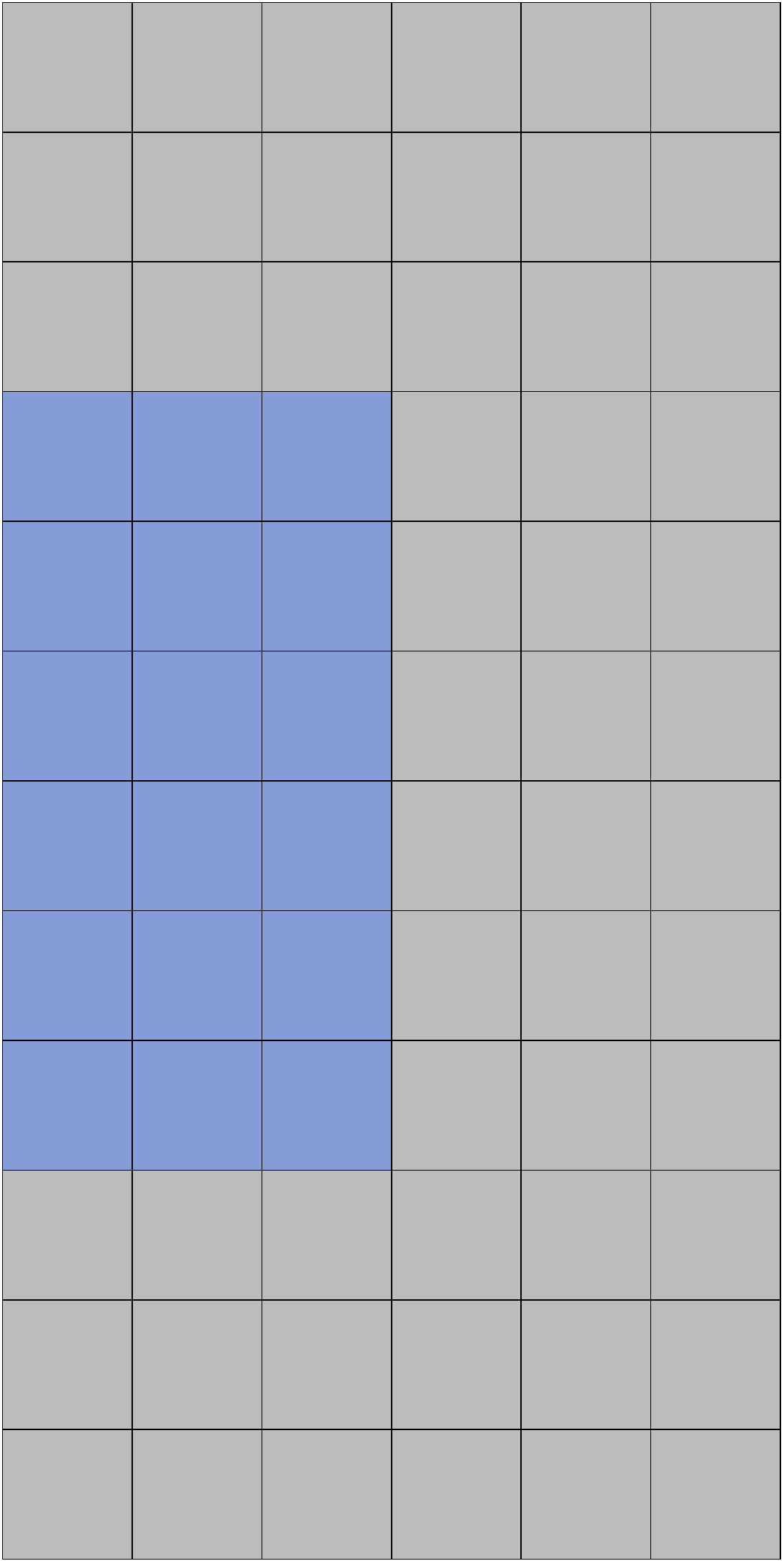} \quad
\includegraphics[keepaspectratio=true,width=2.cm]{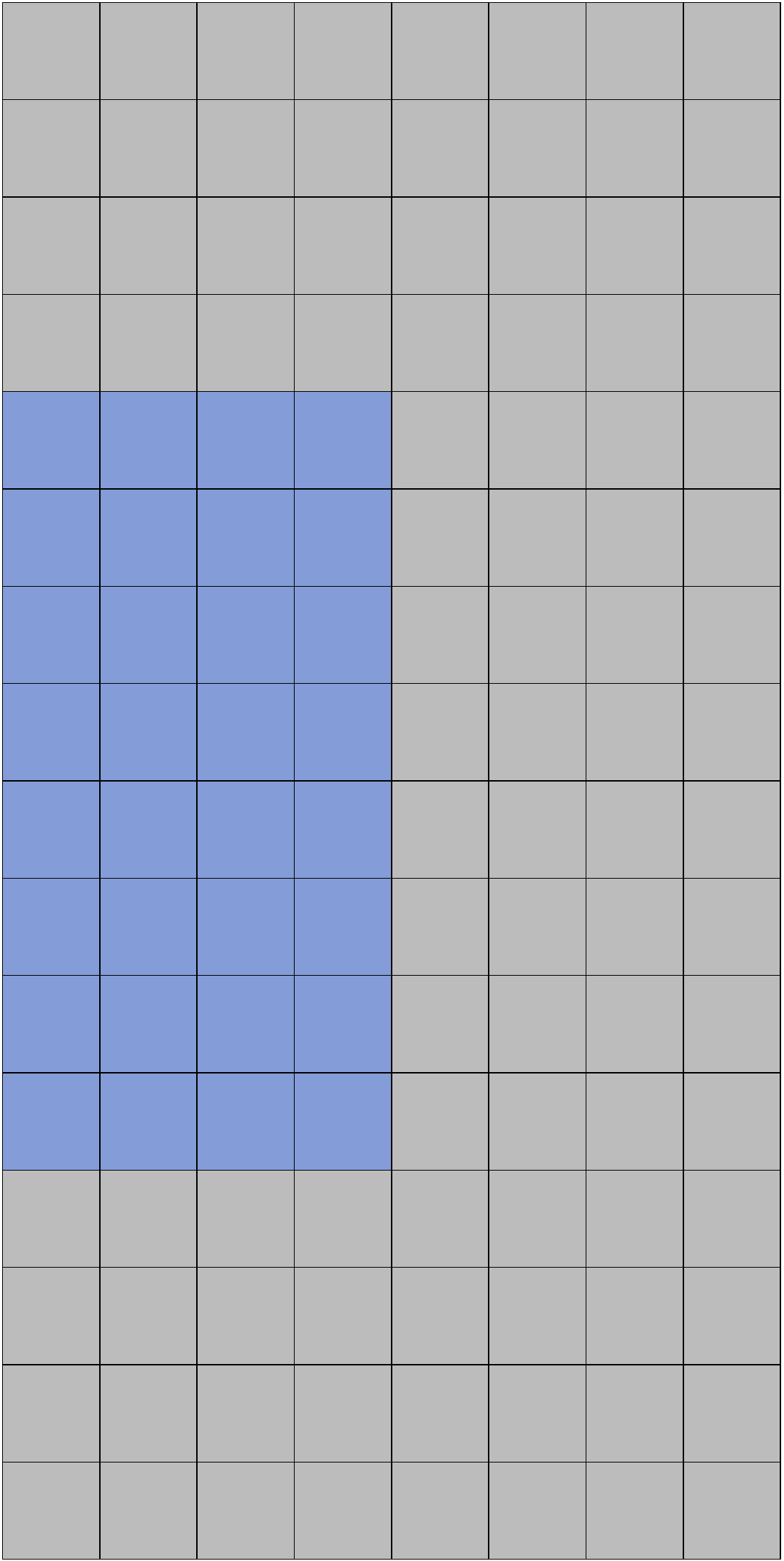}
\caption{Meshes $\MM_1$, $\MM_2$, $\MM_3$, and $\MM_4$ for configuration A}
\label{MeshRect0}
\end{center} 
\end{figure}

\begin{figure}[ht]
\begin{center}
\includegraphics[keepaspectratio=true,width=7cm]{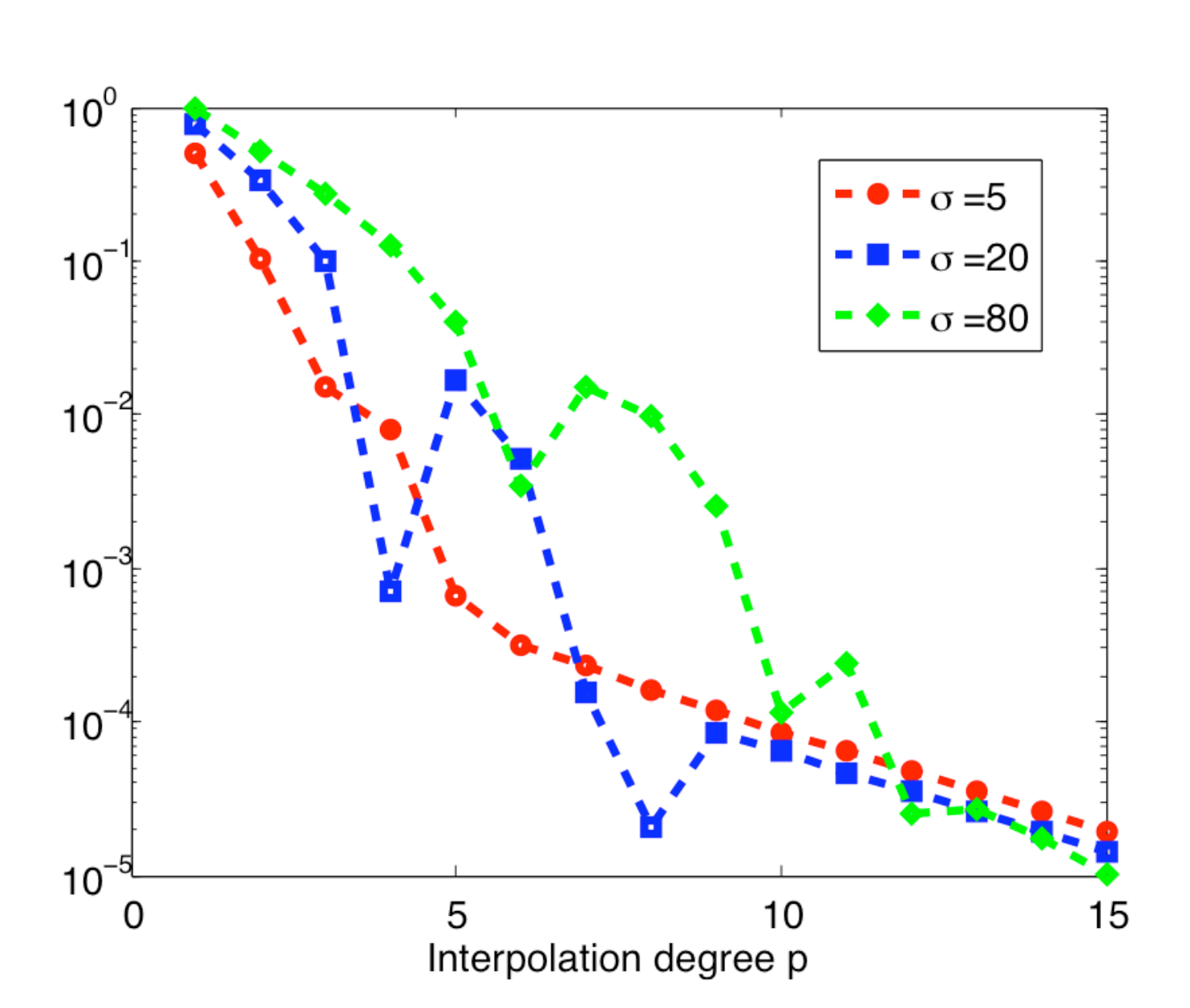}
\caption{Graph of $\big|A_{\sigma}^{p,\MM_2}-A_{\sigma}^{16,\MM_3}\big|$ with respect to $p=1,\cdots,15$ in semi-log coordinates, for $\sigma\in\{5,20,80\}$ for configuration A} 
\label{NormL2cd}
\end{center} 
\end{figure}

We consider a family of eight meshes with square elements $\MM_k$, $k=1,\ldots,8$ with size $h=1/k$, see Figure~\ref{MeshRect0}.
We represent in Figure~\ref{NormL2cd} the absolute value of the difference between $A_{\sigma}^{p,\MM_2}$ and $A_{\sigma}^{16,\MM_3}$, versus $p$ in semilogarithmic coordinates, and in each case: $\sigma=5$ with circles, $\sigma=20$ with squares, and $\sigma=80$ with diamonds. The figure shows that $A_{\sigma}^{p,\MM_2}$ approximates $A_{\sigma}^{16,\MM_3}$ better than $10^{-4}$ when $p\geqslant 12$. 

\clearpage

\subsection{Stability of the h-version}
\begin{figure}[ht]
\begin{center}
\includegraphics[keepaspectratio=true,width=7cm]{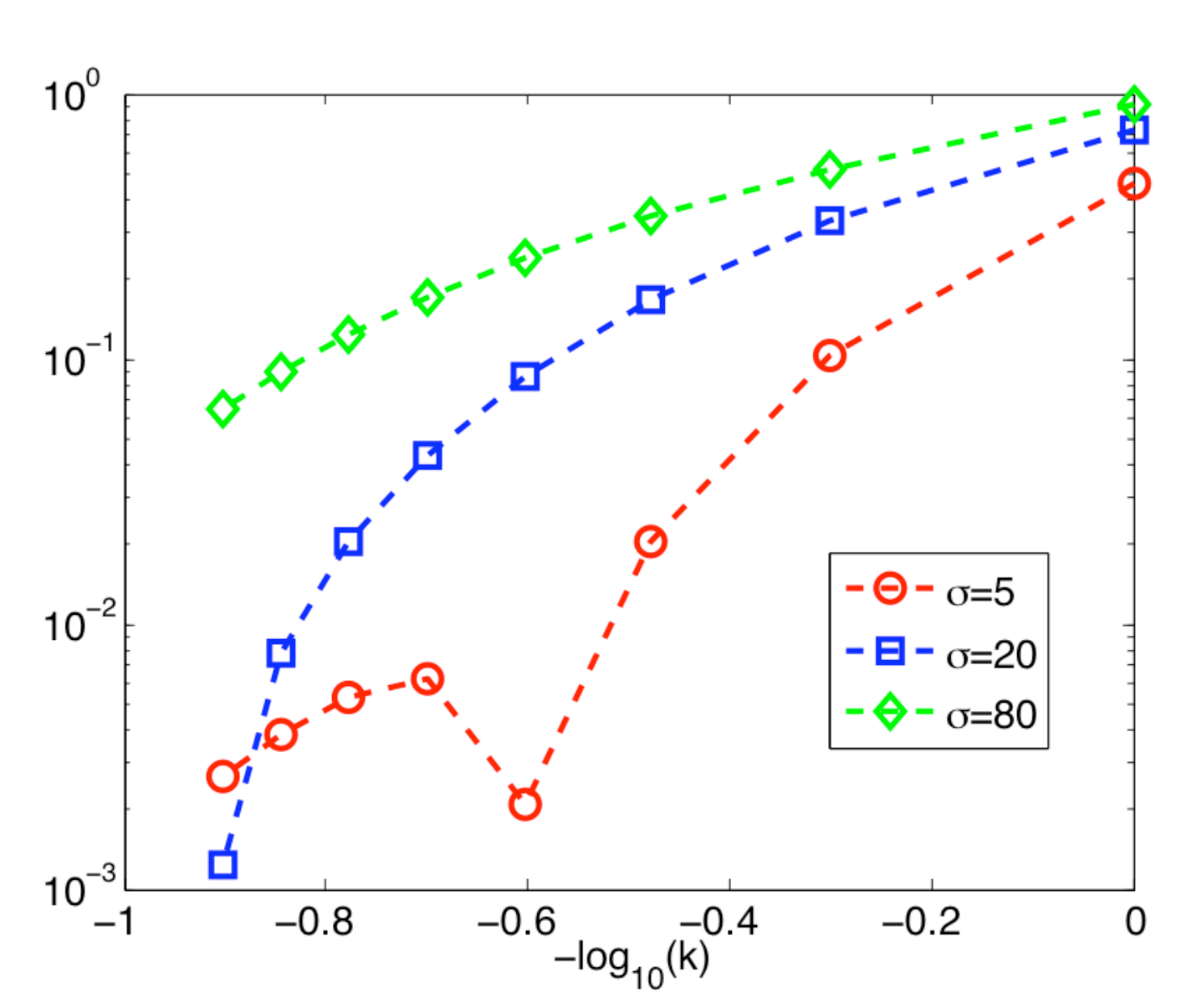}   
\caption{Graph of $\big|A_{\sigma}^{2,\MM_k}-A_{\sigma}^{6,\MM_8}\big|$ with respect to $-\log_{10}(k)$ in semi-log coordinates, when $k=1,..,8$ and $\sigma=5,20,80$ for configuration A}
\label{FQ12rect}
\end{center} 
\end{figure}

In this subsection, we check the convergence of the discretized problem in configuration A for several meshes $\MM_k$ of the computational domain when $k$ increases. We fix the interpolation degree $p=2$ of the finite elements and use the square meshes $\MM_k$.
We plot in Figure~\ref{FQ12rect} the absolute value of the difference between $A_{\sigma}^{2,\MM_k}$ and $A_{\sigma}^{6,\MM_8}$ with respect to $k=1,..,8$, in log-log coordinates. 

 
\section{Numerical simulations of skin effect}
\label{S6} 
Recall that the asymptotic expansion of $\sH_{(\delta)}$ is described in \S \ref{AEh}. According to Remark \ref{Disc} the first term $\sH_{0}^+$ is a real valued function. Hence the imaginary part of the magnetic field $\sH_{(\delta)}$ is small in the dielectric $\Omega_{+}^\m$, because
$$
|\Im \sH_{(\delta)}^+|= \mathcal{O}(\delta) \ .
$$
Thus, the imaginary part of the computed field is located in the conductor $\Omega\con^\m$. We display this imaginary part to highlight the boundary layer near the surface of the conductor, see Figures~\ref{CImagin}, \ref{EModul} bottom, \ref{MCSEa2} and \ref{MCSEa4}.  In this section, $\tH_{\sigma}:=\sH_{(\delta)}^{p,\MM_k}$ denote the computed solution associated with the numerical parameters considered in each subsection.

\subsection{Skin effect in configuration A} 
\label{S6.1}
\begin{figure}[htbp]
\begin{center}
\begin{tabular}{cc}
\includegraphics[keepaspectratio=true,width=6cm]{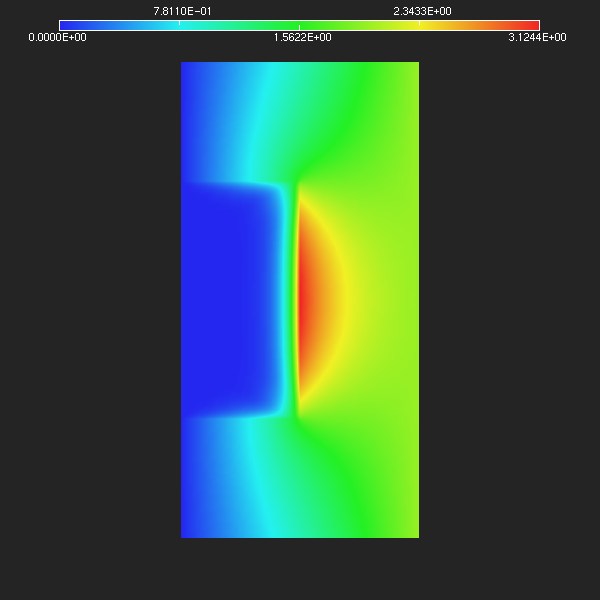}  &
\includegraphics[keepaspectratio=true,width=6cm]{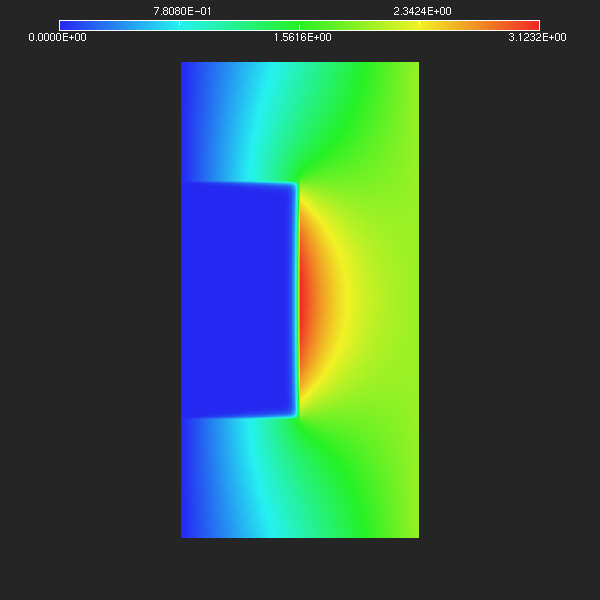} 
\end{tabular}
 \caption{Configuration A. On the left, $| \tH_{\sigma}|$ when  $\sigma=5$. On the right, $| \tH_{\sigma}|$ when  $\sigma=80$}
\label{CModul}
\end{center} 
\end{figure}

\begin{figure}[ht]
\begin{center}
\begin{tabular}{ccc}
 \includegraphics[keepaspectratio=true,width=6cm]{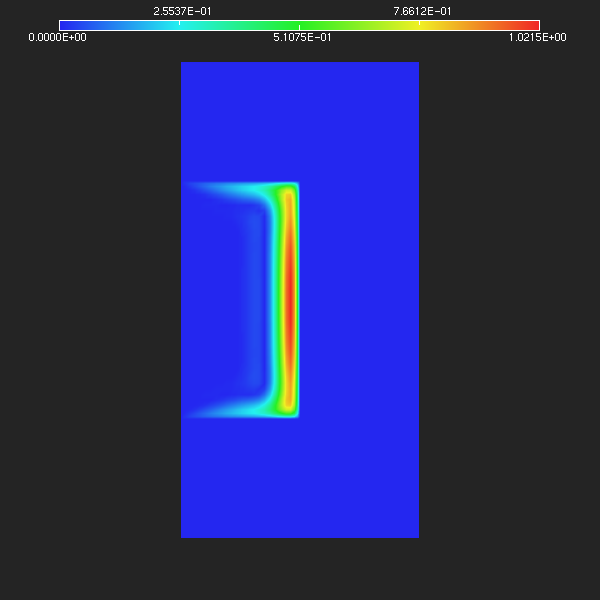}  
 &
 \includegraphics[keepaspectratio=true,width=6cm]{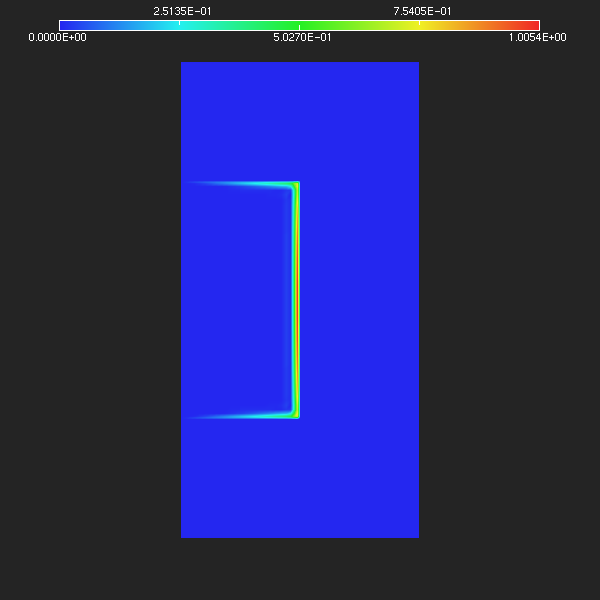}  
\end{tabular}
 \caption{Configuration A. On the left, $|\Im \tH_{\sigma}|$ when  $\sigma=5$. On the right, $|\Im \tH_{\sigma}|$ when $\sigma=80$}
\label{CImagin}
\end{center} 
\end{figure}

\begin{figure}[ht]
\begin{center}
\begin{tabular}{ccc}
\includegraphics[keepaspectratio=true,width=6cm]{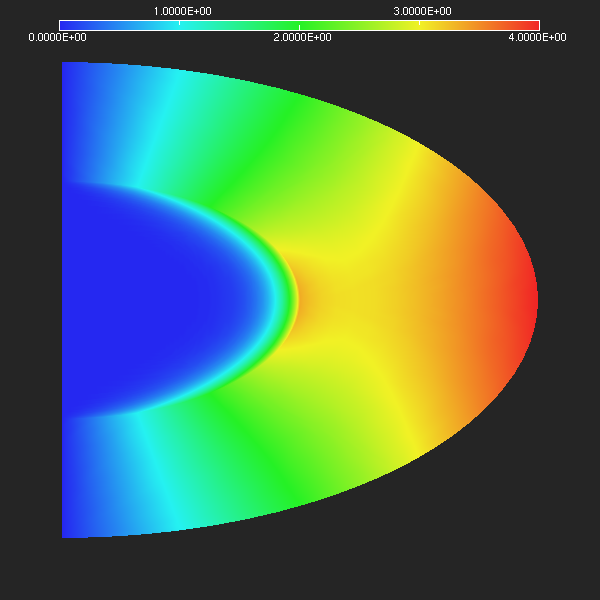}  &
\includegraphics[keepaspectratio=true,width=6cm]{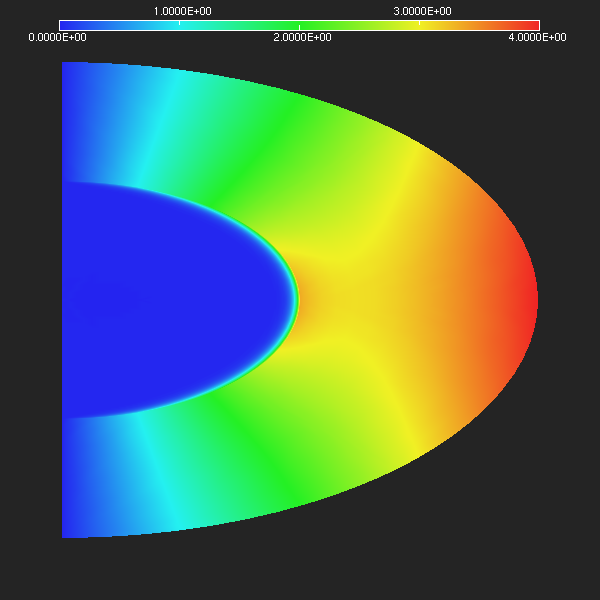} 
\\
\includegraphics[keepaspectratio=true,width=6cm]{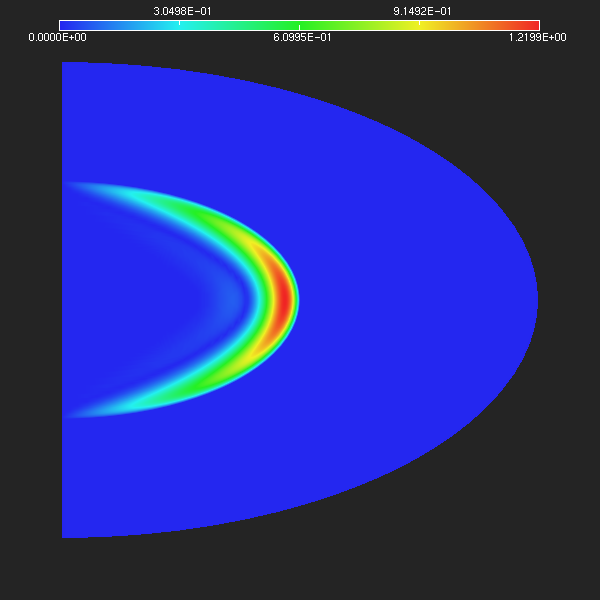}  &
\includegraphics[keepaspectratio=true,width=6cm]{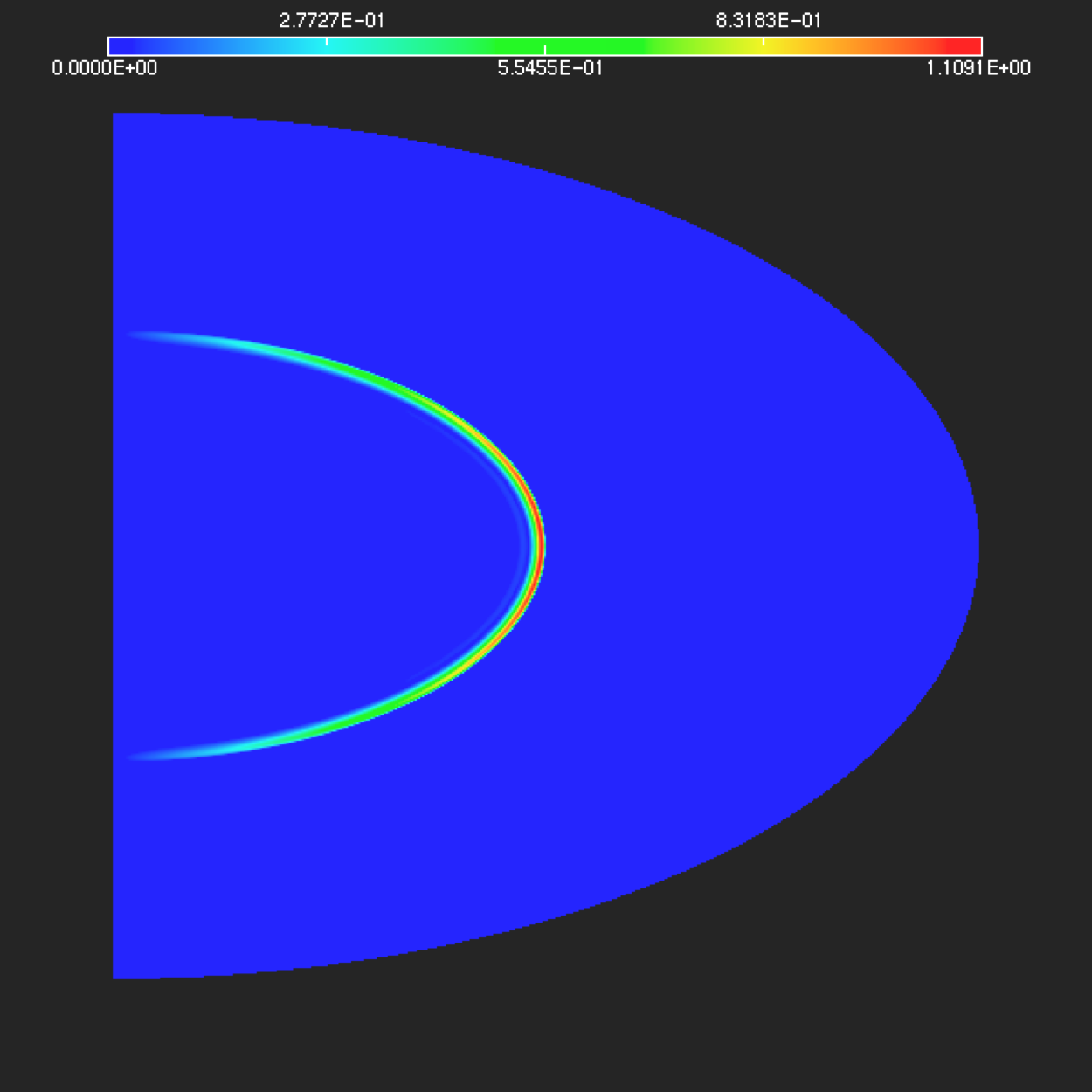}
\end{tabular}
 \caption{Configuration B1. At the top, $|\tH_{\sigma}|$ when  $\sigma=5$ (on the left), and $\sigma=80$ (on the right). At the bottom, $|\Im \tH_{\sigma}|$ when  $\sigma=5$ and $\sigma=80$}
\label{EModul}
\end{center} 
\end{figure}

We fix the mesh $\MM_2$, see Figure~\ref{MeshRect0}, and the interpolation degree of the finite elements : $p=16$. Here $\tH_{\sigma}=\sH_{(\delta)}^{16,\MM_2}$ is computed in configuration A for several values of $\sigma$, see Figure~\ref{CModul}. Similarly, we compute $|\Im \tH_{\sigma}|$, see Figure~\ref{CImagin}.


\subsection{Skin effect in configuration B}
\label{S6.2}
Here $\tH_{\sigma}=\sH_{(\delta)}^{p,\MM_{3}}$ is the computed solution for several values of $\sigma$, with a fixed mesh $\MM_{3}$ (Figure~\ref{F0}). We represent $|\tH_{\sigma}|$ and $|\Im \tH_{\sigma}|$ in configuration B1, with an interpolation degree $p=16$, see Figure~

\subsection{Skin effect in configurations B and C} 
\begin{figure}[ht]
\begin{center}
\begin{tabular}{ccc}
 \includegraphics[keepaspectratio=true,width=6cm]{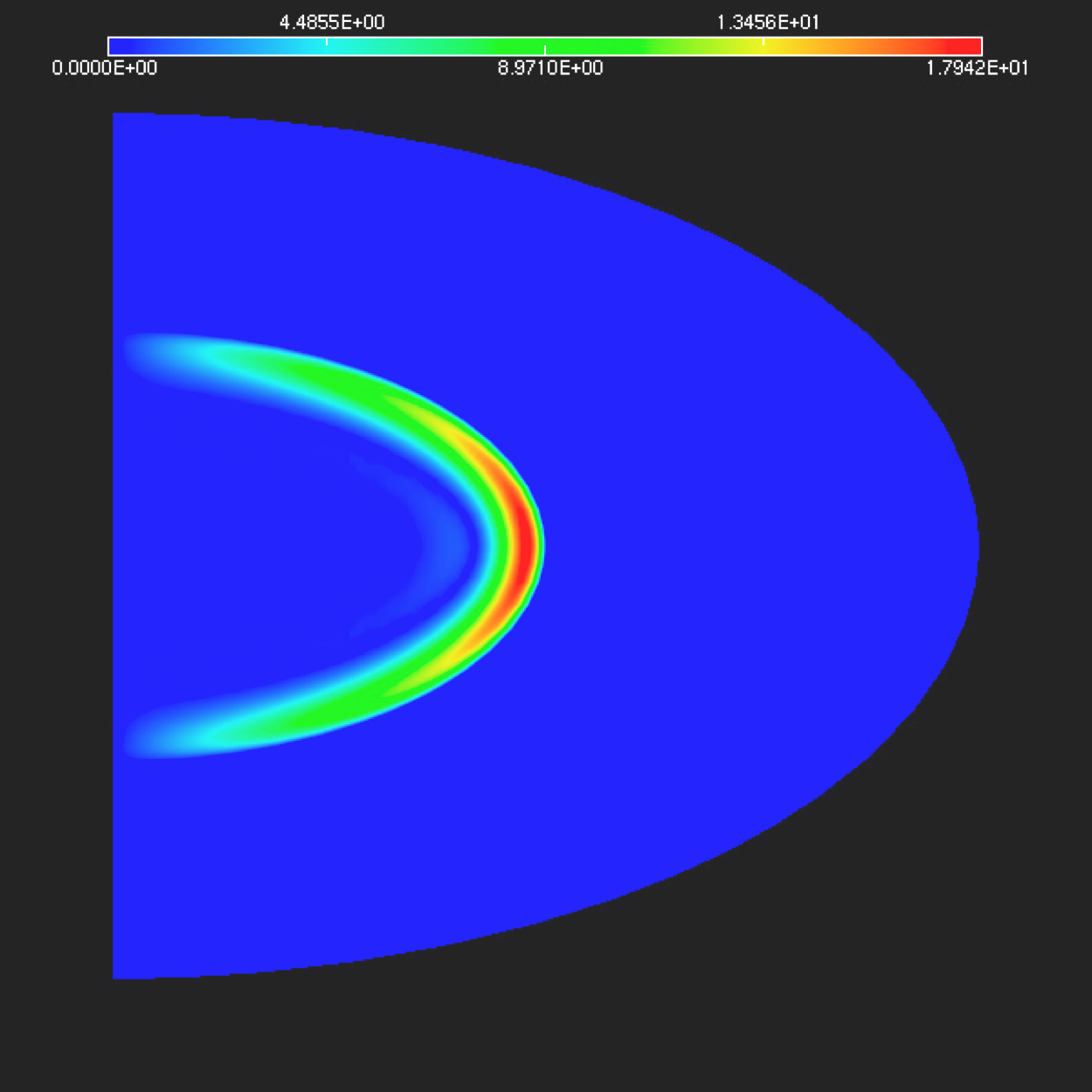}  &
\includegraphics[keepaspectratio=true,width=6cm]{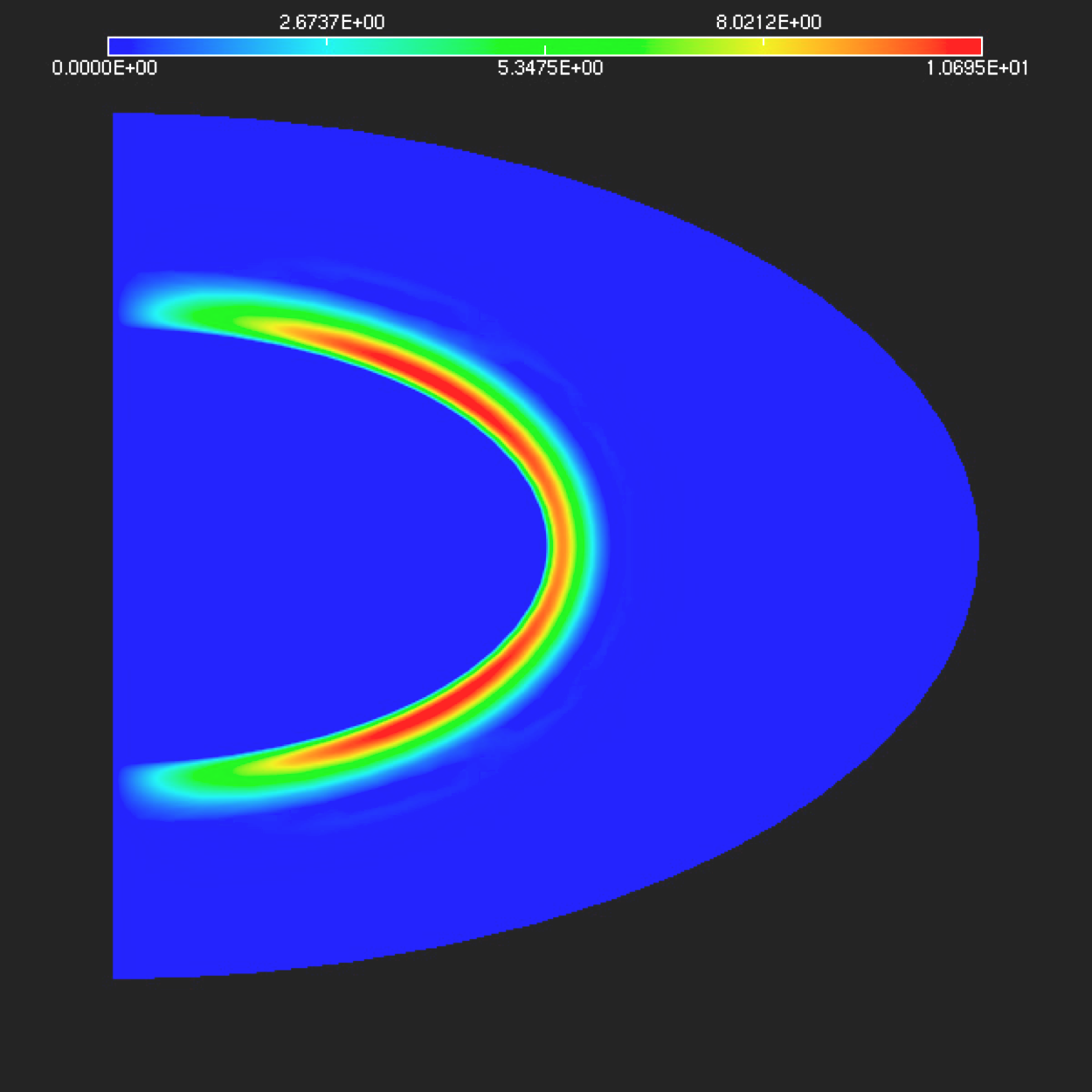}  
\end{tabular}
\caption{ On the left, $|\Im \tH_{\sigma}|$ in configuration B1 ($\cH>0$). On the right, $ |\Im \tH_{\sigma}|$  in configuration C1 ($\cH<0$). $\sigma=5$}
\label{MCSEa2}
\end{center}
\end{figure}

\begin{figure}[ht]
\begin{center}
\begin{tabular}{cc}
 \includegraphics[keepaspectratio=true,width=6cm]{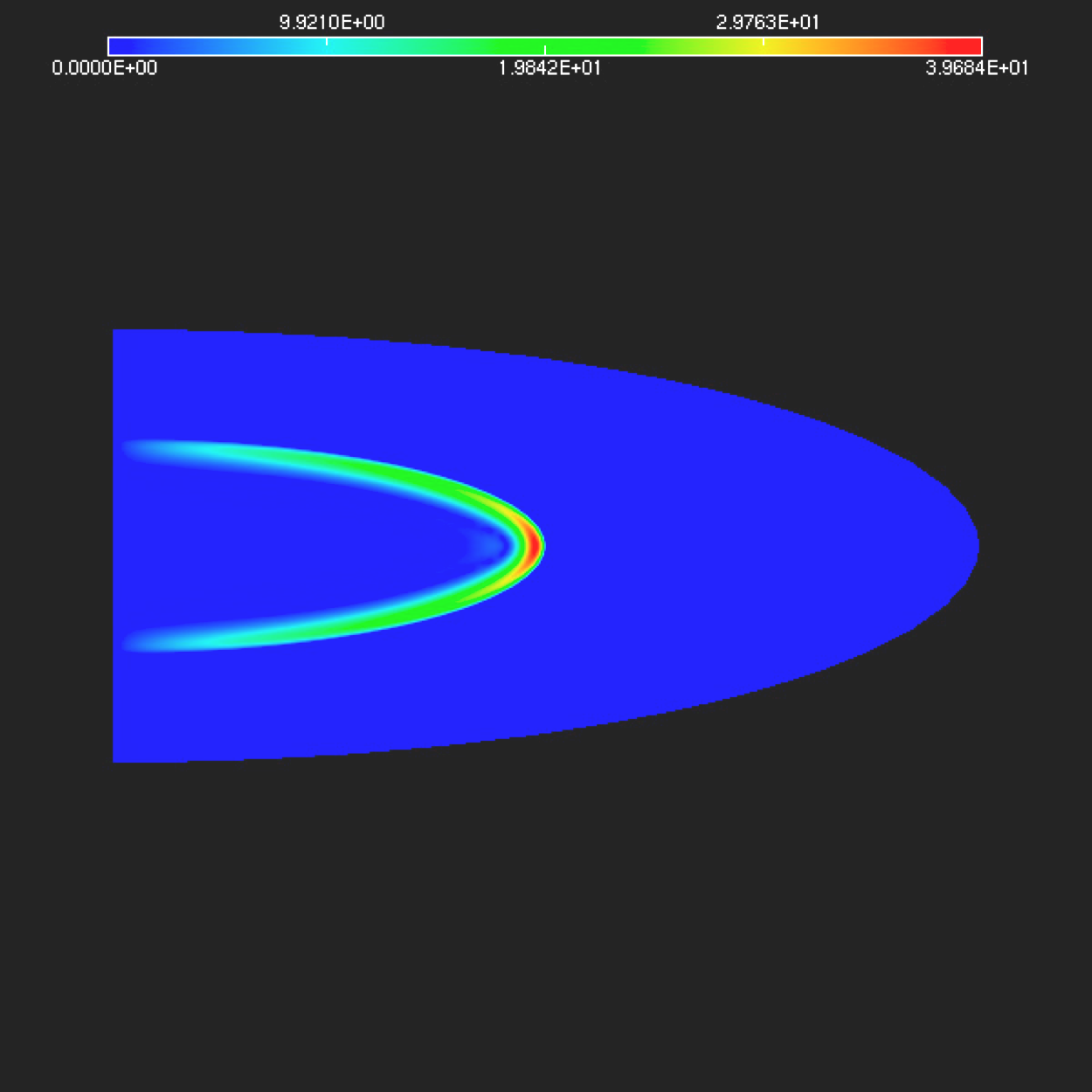}  &
\includegraphics[keepaspectratio=true,width=6cm]{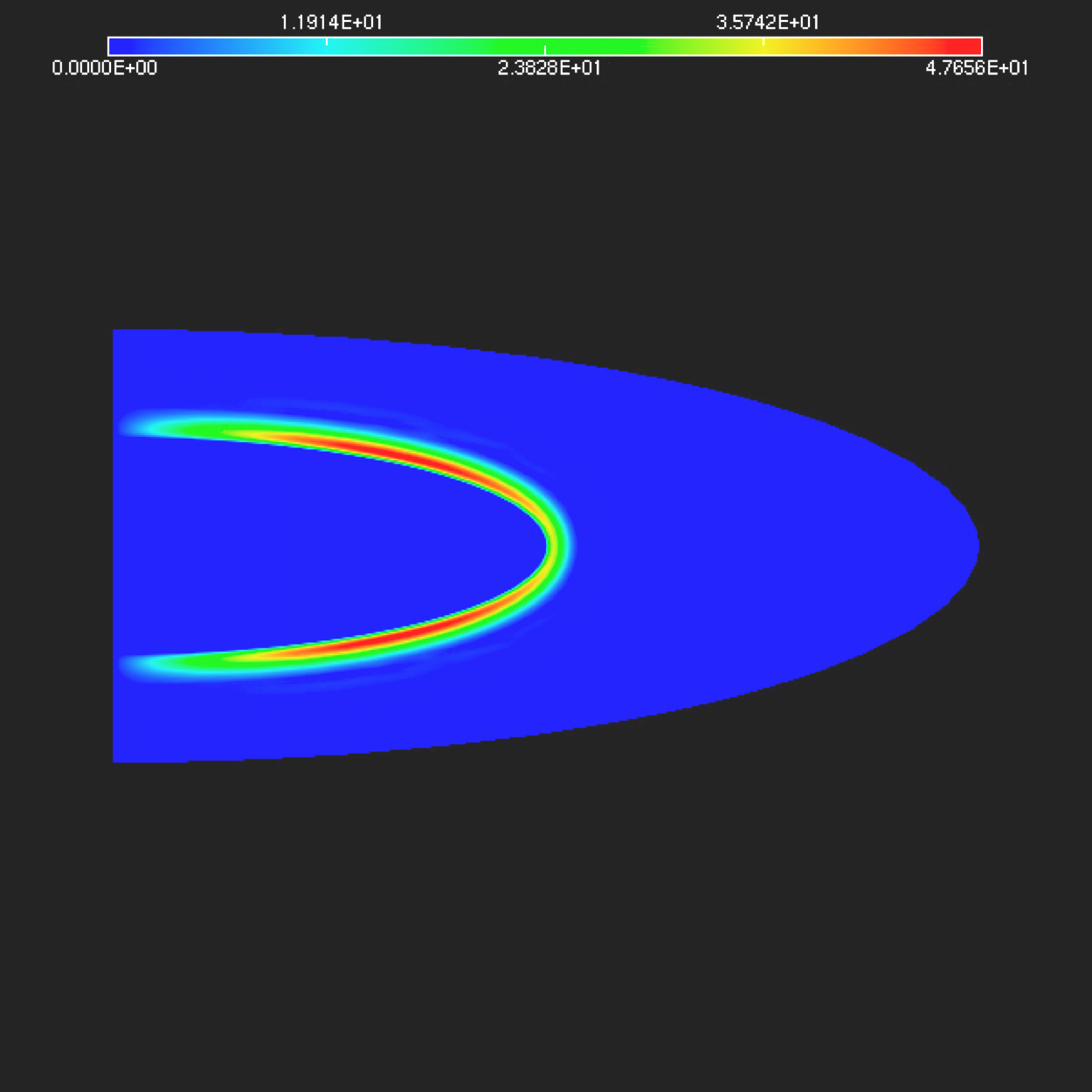}  
\end{tabular} 
\caption{On the left, $|\Im \tH_{\sigma}|$ in configuration B2. On the right, $|\Im \tH_{\sigma}|$ in configuration C2. $\sigma=5$}
\label{MCSEa4}
\end{center}
\end{figure}

In order to exhibit the influence of the sign of the mean curvature $\cH$ of the interface $\Sigma$ on the skin effect, we fix in this subsection the conductivity $\sigma=5$. We perform computations in configurations B and C. We note that $\cH>0$ in configuration B, and $\cH<0$ in configuration C. We compute $|\Im \tH_{\sigma} |$  in configurations B1 and C1, see Figure~\ref{MCSEa2}. 
We then compute $| \Im \tH_{\sigma} |$ in  configurations B2 and C2: $\cH>0$, and $\cH<0$, see Figure~\ref{MCSEa4}. 

Figures \ref{MCSEa2} and \ref{MCSEa4} show that the skin depth is larger for a fixed conductivity when the mean curvature of the conducting body surface is larger. Moreover, the sign of the curvature has an influence on the skin depth. This length is larger in convex (Figures \ref{MCSEa2} and \ref{MCSEa4} on the left) than in concave conductors (same figures on the right).

\clearpage

\section{Postprocessing}
\label{S7}
In this section, we perform numerical treatments from computations in configuration B1, see subsection \ref{S6.2}, and configuration A, see subsection \ref{S6.1}, in order to investigate whether solutions are exponentially decreasing inside the conductor and with which rate. Let us recall that the standard skin depth $\ell(\sigma)$ is given by \eqref{El}.

\subsection{Configuration B}
The mesh of the computational domain is the mesh $\MM_3$, see Figure~\ref{F0}. We extract values of $|\tH_{\sigma}|$ in $\Omega\con^\m$ along edges of the mesh $\MM_3$ for $z=0$: in this configuration, the normal coordinate writes $y_{3}:=2-r$. 

\begin{figure}[ht]
\begin{center}
\begin{tabular}{cc}
\includegraphics[keepaspectratio=true,width=7cm]{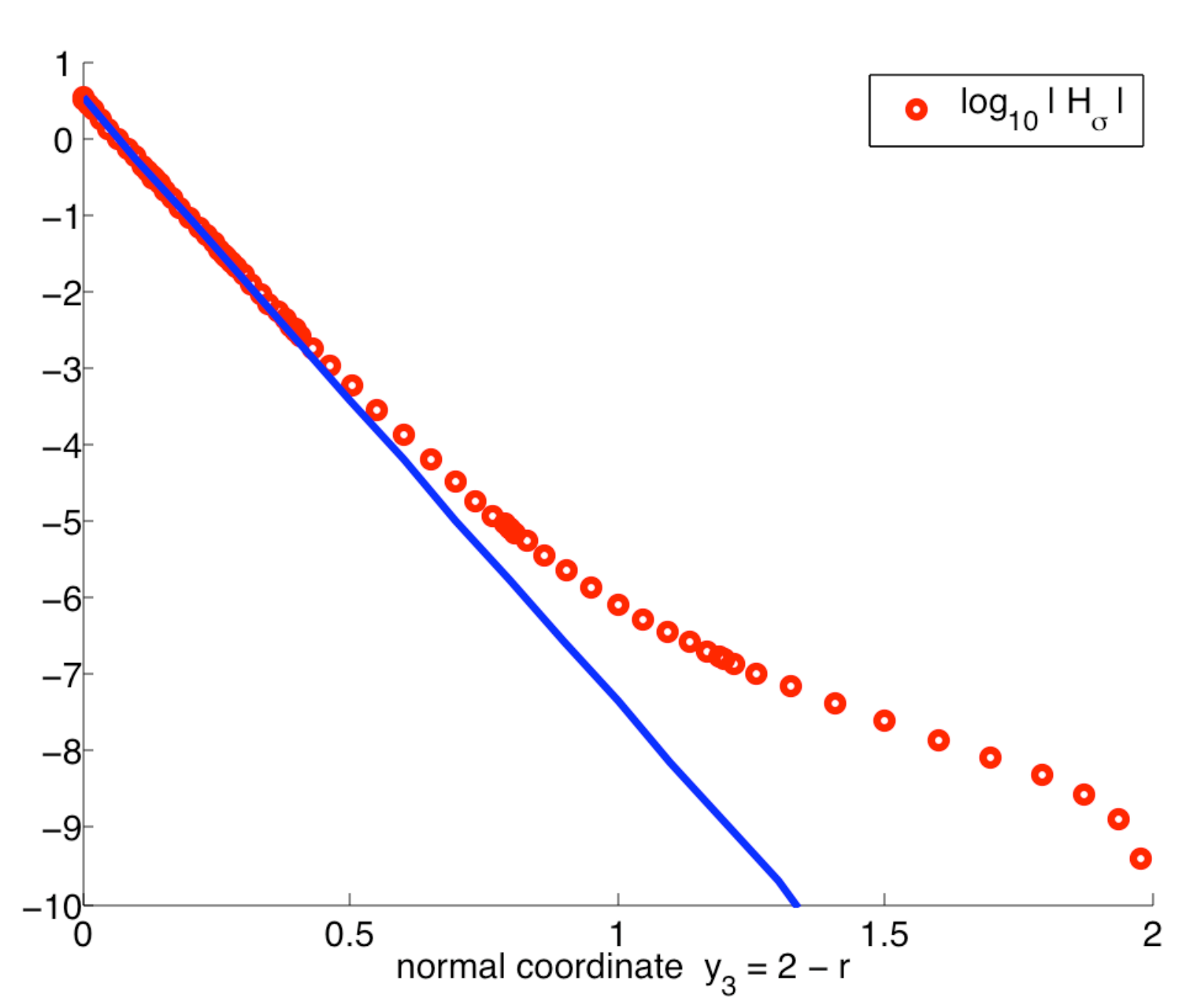}
&
\includegraphics[keepaspectratio=true,width=7cm]{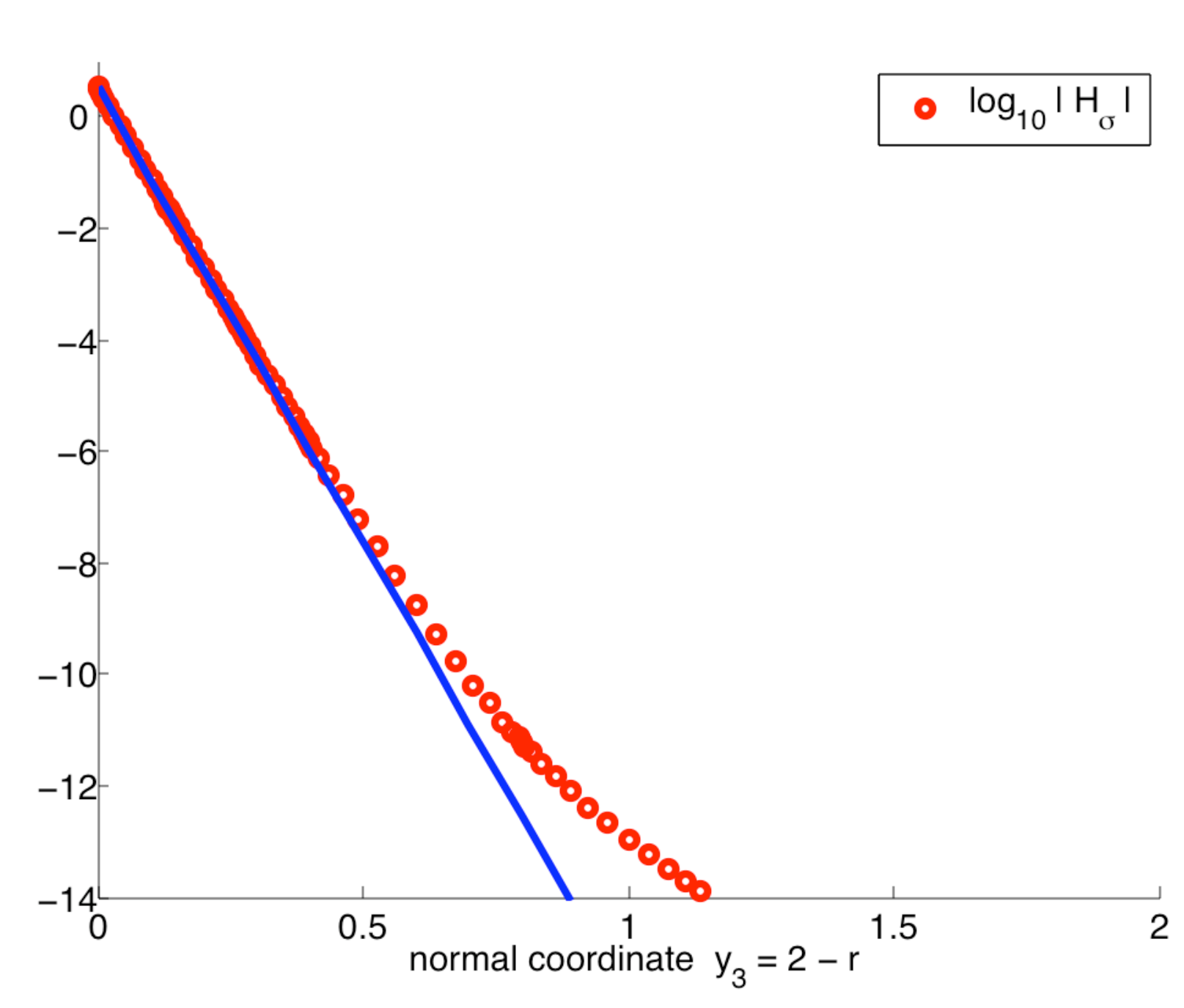}  
\end{tabular}
\caption{In circles, extracted values of $\log_{10} |\tH_{\sigma}{(y_{3})}|$. In solid line, linear regression of $\log_{10} |\tH_{\sigma}{(y_{3})}|$. On the left  $\sigma=20$ (interpolation degree $p=12$). On the right, $\sigma=80$ (degree $p=16$).  $n(20)=6$\,, $\tilde s(20)=7.86906$ and $n(80)=5$\,, $\tilde s(80)=16.29634$. Configuration B1}
 \label{Fextr}
\end{center}
\end{figure}

Then, we perform a linear regression from values of $\log_{10}|\tH_{\sigma}(y_{3})|$ in the skin depth $\ell(\sigma)$, see Figure~\ref{Fextr}. We denote  by $n(\sigma)$ the number of extracted values on the axis $z=0$ in the skin depth $\ell(\sigma)$. From the linear regression, we derive a numerical slope $\tilde s(\sigma)$ such that
$$
\log_{10} | \tH_{\sigma}(y_{3})| = - \tilde s(\sigma) y_{3} + b \ , \quad  b\in\R  \ .
$$

\subsubsection*{Accuracy of asymptotics}
Relying on formula \eqref{EV}, we can derive a Taylor expansion of  
$\log_{10}  |\VV_{(\delta)}(y_{\alpha},y_{3}) |$ with respect to $y_{3}$: for $y_{3}$ small enough
\begin{multline}
\label{ElogV}
\log_{10}  |\VV_{(\delta)}(y_{\alpha},y_{3}) | =  
\log_{10} |\bsH_{0}(y_{\alpha}) |  
\\
- s(y_{\alpha},\sigma) y_{3}
+ \frac{\delta}{\ln 10}  | \bsH_{0}(y_{\alpha}) |^{-2} 
   \Re \big\langle  \bsH_{0}(y_{\alpha}) , \bsH_{1}(y_{\alpha}) \big\rangle
+ \cO\big((\delta+y_{3})^2\big) \ .
\end{multline}
Here the function
$$
s(y_{\alpha},\sigma):= \frac1{\ln 10}\Big(\frac1{\ell(\sigma)}- \cH(y_{\alpha}) \Big)  
$$ 
depends on the skin depth $\ell(\sigma)$ and the mean curvature $\cH(y_{\alpha})$ at the point $y_{\alpha}$ of the interface $\Sigma$. We note that the mean curvature of the surface $\Sigma$ is constant when $z=0$. Hence, we introduce hereafter the theoretical slope 
$$
 s(\sigma):=s(y_{\alpha},\sigma) \quad \mbox{when}\quad z=0 \ .
$$
In configuration B1, the principal curvatures at a point of the interface $\Sigma$ when $z=0$ are the constants $\kappa_{1}=2$ and $\kappa_{2}=\frac1{2}$. Hence, $\cH(y_{\alpha})=\frac{5}{4}$. We infer 
$$
s(\sigma)= \frac1{\ln 10}\Big(\frac1{\ell(\sigma)}- \frac{5}{4}\Big)\, .
$$
The accuracy of the asymptotic expansion is tested by representing the relative error between numerical and theoretical slopes: 
$$
\err(\sigma):=\Big|\frac{{s(\sigma)-\tilde s(\sigma)} }{s(\sigma)}\Big| \ ,
$$
see the table~\ref{a2}. 
In order to make clear wether the influence of the curvature is visible in computations, we also display the theoretical curvature ratio
\[
   {\rm curv\_ratio}(\sigma) := \frac{\frac54}{\frac1{\ell(\sigma)}- \frac{5}{4}}.
\] 
\Bk

\begin{table}[ht]
\begin{center}
\begin{tabular}{|l|c|c|c|}
\hline
$\sigma$ & $5$ & $20$ & $80$ \\
\hline
\hline
$\ell(\sigma)$  & $0.103$ & $0.0515$ & $0.0258$ \\
\hline
$s(\sigma)$ & $3.67332$ & $7.88951$ & $16.32188$ \\
\hline
${\rm curv\_ratio}(\sigma)$  & $0.148$ & $0.069$ & $0.033$ \\
\hline
\hline
degree $p$  & $10$ & $12$ & $16$ \\
\hline
$n(\sigma)$ & $7$ & $6$ & $5$ \\
\hline 
$\tilde s(\sigma)$ & $3.64686$ & $7.87347$ & $16.308279$ \\
\hline
$\err(\sigma)$ & $0.0072$ & $0.002$ & $0.0008$ \\
\hline
\end{tabular}
\bigskip

\caption{Postprocessing in configuration B1 with the mesh $\MM_3$}
\label{a2}
\end{center}
\end{table}

The relative error decreases more than a half when the conductivity is multiplied by $4$. 
Anyway, this relative error is much smaller than ${\rm curv\_ratio}(\sigma)$, which exhibits numerically the influence of the curvature on the skin depth.
We perform similar computations with the mesh  $\MM_6$ represented in Figure~\ref{F0}, see Table \ref{tK6}, and obtain still better results.

\begin{table}[ht]
\begin{center}
\begin{tabular}{|l|c|c|c|}
\hline
$\sigma$ & $5$ & $20$ & $80$ \\
\hline
\hline
degree $p$  & $8$ & $12$ & $16$ \\
\hline
$n(\sigma)$ & $13$ & $9$ & $7$ \\
\hline
$\tilde s(\sigma)$ & $3.64239$ & $7.88170$ & $16.33051$ \\
\hline
$\err(\sigma)$ & $0.0084$ & $0.001$ & $0.0005$ \\
\hline
\end{tabular}
\bigskip

\caption{Postprocessing in configuration B1 with the mesh $\MM_6$ (see Table \ref{a2} for the theoretical values $\ell(\sigma)$, $s(\sigma)$, and ${\rm curv\_ratio}(\sigma)$)}
\label{tK6}
\end{center}
\end{table}

\begin{figure}[ht]
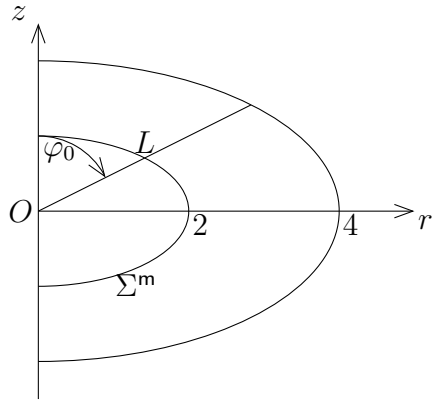

\begin{center}
\figinit{1.cm}
\def\rx{2} \def\ry{1} \def\Rx{4} \def\Ry{2}
\figpt 6: $O$ (0,0)
\figpt 7: $r$ (5,0)
\figpt 8: (0,-2.5)
\figpt 9: $z$ (0,2.5)
\figpt 1:(0,2)
\figpt 2:(2,2)
\figpt 3:(2,-2)
\figpt 4:(0,-2)
\figpt 5:(2,0)
\figpt 20:(4,0)
\figpt 16:(0,4)
\figpt 17:(5,4)
\figpt 18:(5,-4)
\figpt 19:(0,-2)
\figpt  25:(1.4142,0.7071) 
\figpt  26:(2.8284,1.4142)
\figpt  27:(0,3.1623)
\figptrot 29:= 27/ 6, -63.4377/
\figpt 10: (0,-3) \figpt 11: (0,3)
\figpt 12: (0.5,0.5) \figpt 13: (2,-1.5)
\figpt 14: (3,2.2)
\figpt 15: (1.3,-1.2)
\psbeginfig{}

\psarrow[6,7] \psarrow[8,9]
\psset (width=\defaultwidth)
\psline[6,29]
\def\rap{.2}
\figpthom 22:= 1 /6,\rap/
\figpthom 23:= 29 /6,\rap/
\figptbary 24:[23,22;1,1]

\pssetwidth{\defaultwidth}
\psarrowcircP 6 ; -1 [22,23]

\psarcell 6 ; \Rx, \Ry (-90,90,0)
\psarcell 6 ; \rx, \ry (-90, 90,0) 
\psendfig

\figvisu{\figBoxA}{}{
\figwritew 6:$O$(2pt)
\figwriten 25:$L$(2pt)
\figwritese 5:$2$ (2pt)
\figwritese 20:$4$ (2pt)
\figwritese 7:$r$ (2pt) \figwritenw 9:$z$ (2pt)
\figsetmark{}
\figwriten 15:$\Sigma^\m$(2pt)
\figwriten 24:${\varphi}_{0}$(9pt)
}
\centerline{\box\figBoxA}
\caption{In the meridian domain in configuration B1, $L$ is the point on $\Sigma^\m$  with a colatitude $\varphi_{0}=\pi/2-\arctan{(1/2)}$}
\label{F3}
\end{center}
\end{figure}

\begin{rem}
A similar postprocessing along the $r$-axis was performed for configuration A, see \cite[Ch.\,8, \S 8.3.2]{Pe09}. Relative errors  $\err(\sigma)$ are still consistent with the expansion \eqref{ElogV}. Non-radial postprocessings was performed along the segment $OL$, when $L$ is the point with a colatitude $\varphi_{0}$ in the interface $\Sigma^\m$, see Figure~\ref{F3}: when $\varphi_{0}=\pi/2-\arctan{(1/2)}$ (i.e. $r_{L}=\sqrt{2}, z_{L}=1/\sqrt{2}$), see \cite[Ch.\,8, \S 8.5.1]{Pe09}. Relative errors are again very small.
\end{rem}

\subsection{Configuration A}

\begin{figure}[htp]
\begin{center}
\begin{tabular}{cc}
\includegraphics[keepaspectratio=true,width=7cm]{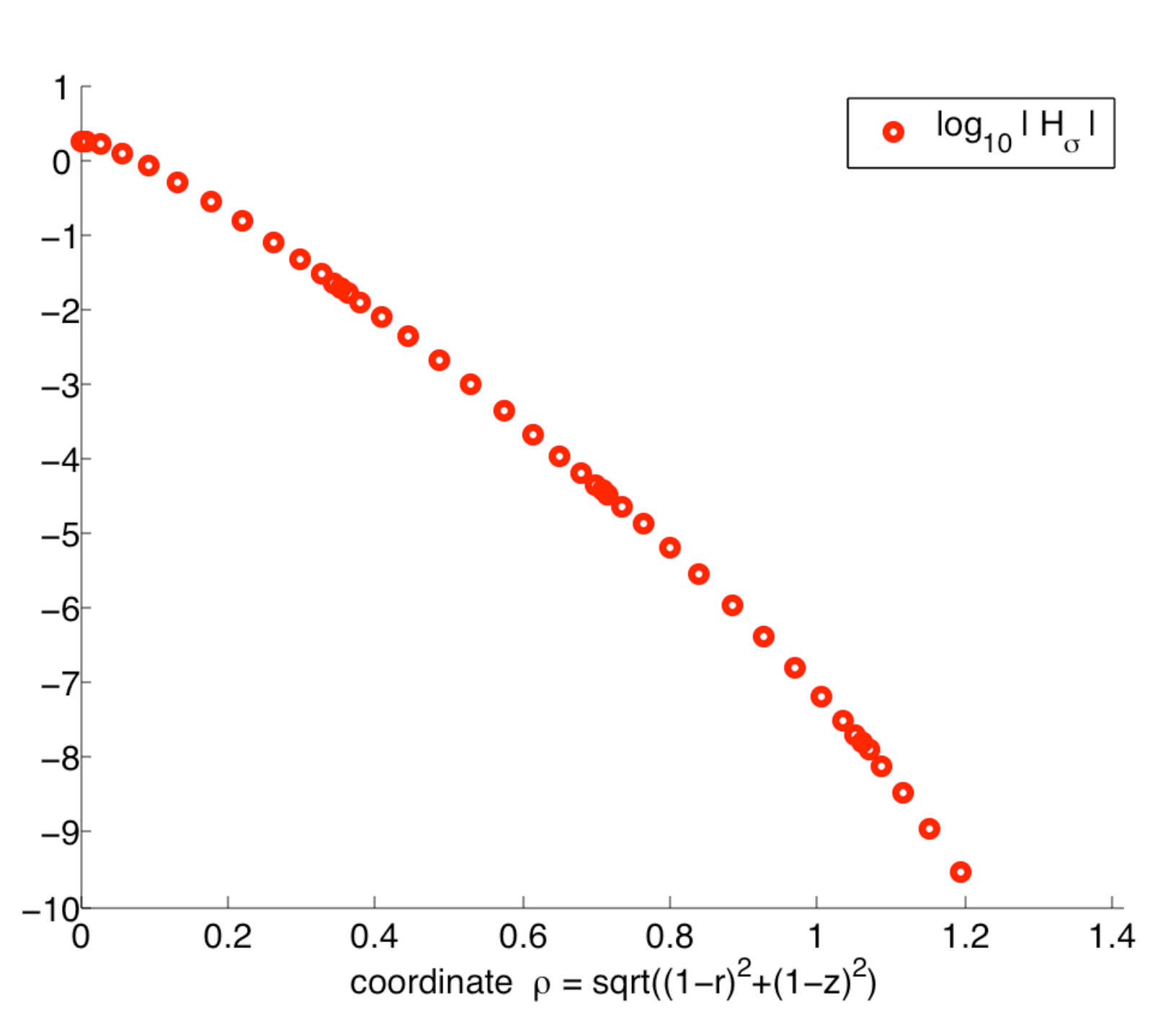}
&
\includegraphics[keepaspectratio=true,width=7cm]{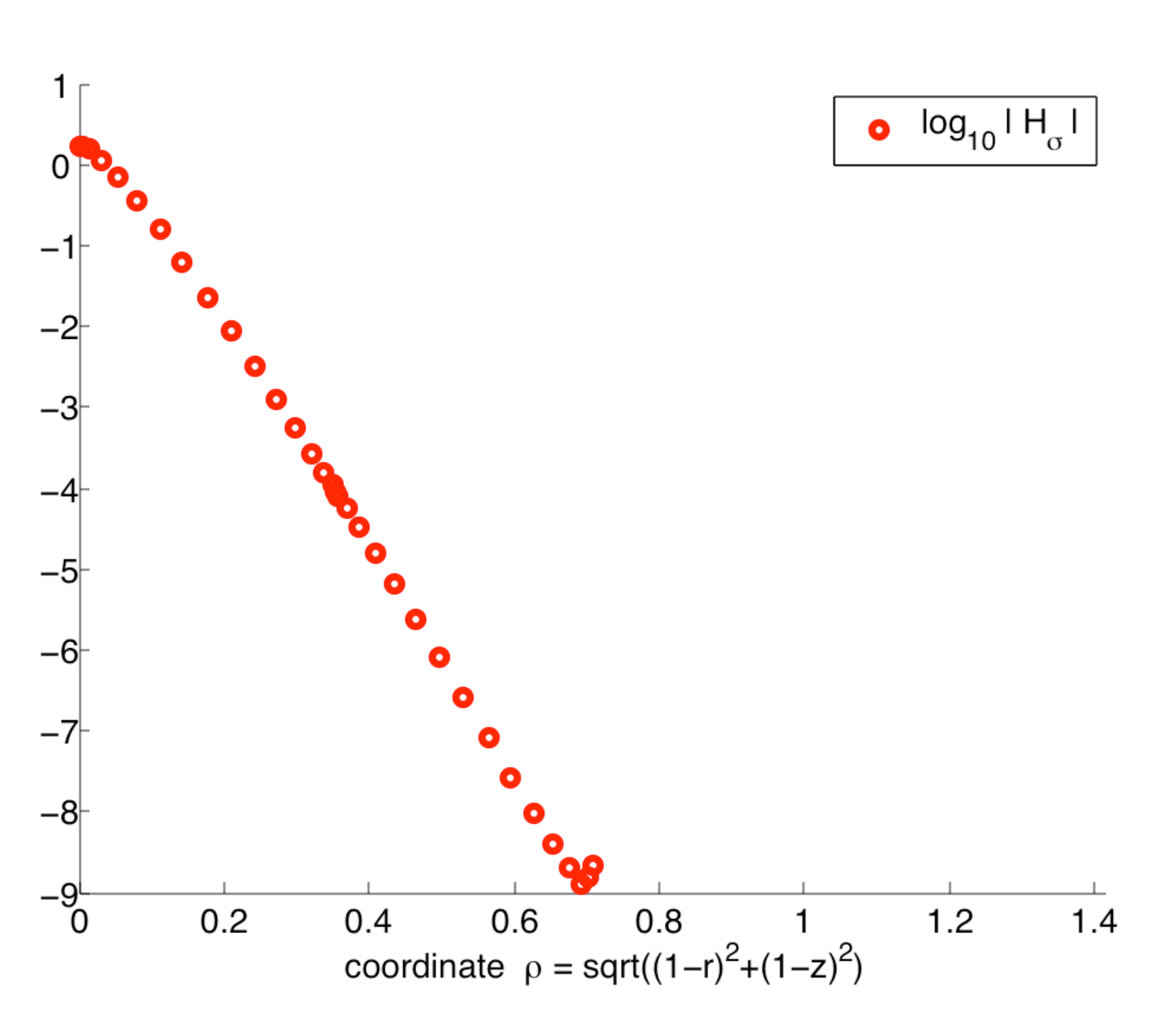}  
\end{tabular}
\caption{In circles, extracted values of $\log_{10} |\tH_{\sigma}{(\rho)}|$, $\rho$ is the distance to the corner. On the left  $\sigma=20$ (interpolation degree $p=12$). On the right, $\sigma=80$ (degree $p=16$). Configuration A}
 \label{extrLogH}
\end{center}
\end{figure}

The mesh of the computational domain is the mesh $\MM_4$ represented in Figure~\ref{MeshRect0}. We extract values of $\log_{10}|\tH_{\sigma}|$ in $\Omega\con^\m$ along the diagonal axis $r=z$, see Figure~\ref{extrLogH}: here we denote by 
$$
\rho:=\sqrt{(1-r)^2+(1-z)^2}\ ,
$$ 
the distance to the corner point $\aa$ with coordinates $(r=1,z=1)$, cf.\ Figure~\ref{FCylGeo}.  

\begin{figure}[ht]
\begin{center}
\includegraphics[keepaspectratio=true,width=8cm]{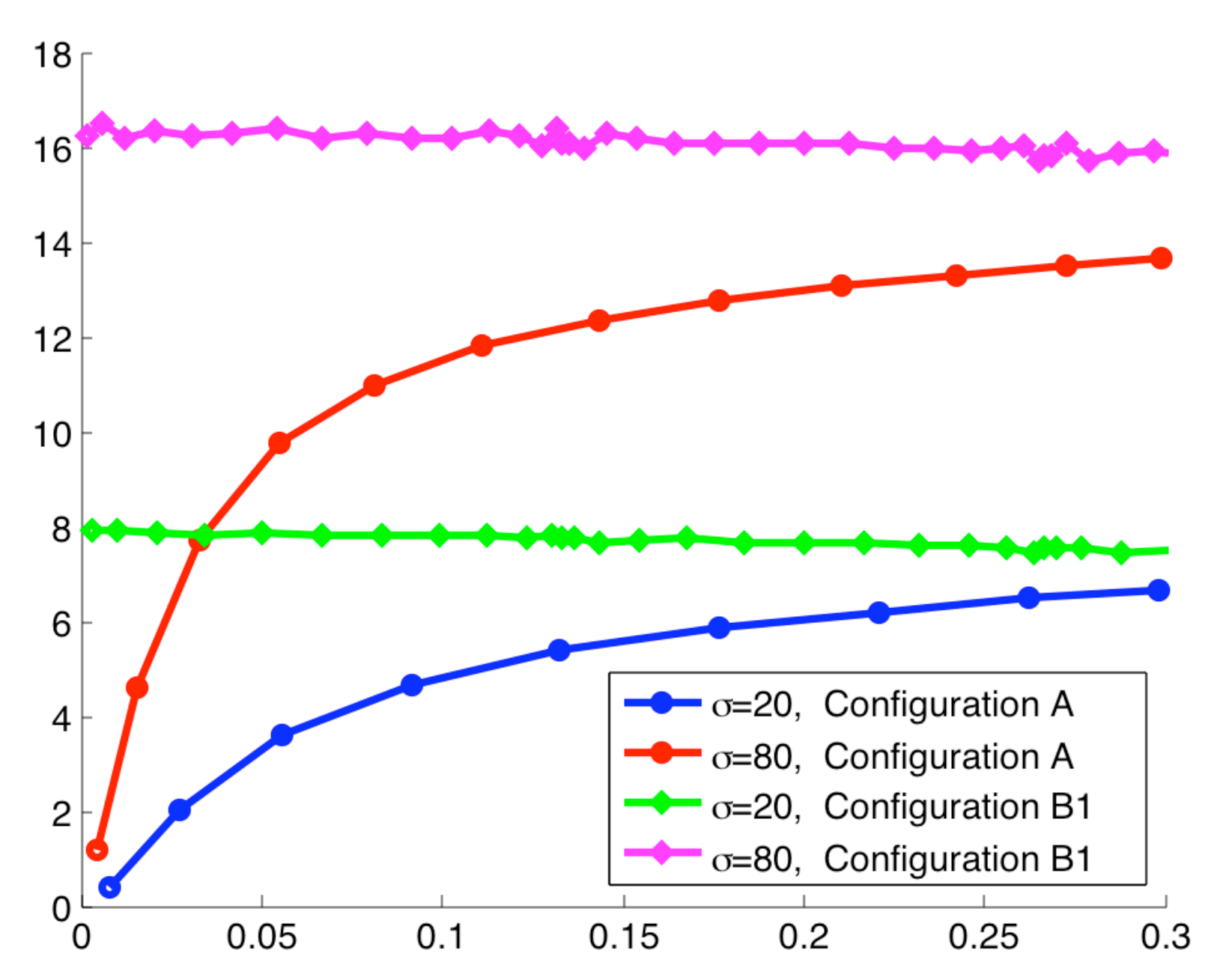}
\caption{The graphs of the slopes $\tilde s_{i}(\sigma)$. Configuration A and B1, and $\sigma=20$ (interpolation degree $p=12$), $\sigma=80$ (degree $p=16$)}
\label{slope}
\end{center}
\end{figure}

 When compared with Figure~\ref{Fextr}, we see that the curves do not exactly behave like lines, which means that the exponential decay is not obvious.
In order to measure a possible exponential decay, we define the slopes $\tilde s_{i}(\sigma)$ of extracted values 
$\log_{10} | \tH_{\sigma}(\rho)|$ by:
$$
\tilde s_{i}(\sigma):=\frac{\log_{10} | \tH_{\sigma}(r_i,z_i)|-\log_{10} | \tH_{\sigma}(z_{i+1},r_{i+1})|}{\rho_{i+1}-\rho_{i}}\ .
$$
Here, $\rho_{i}$ is the distance to the corner point $\aa$ defined by 
$\rho_{i}:=\sqrt{(1-r_{i})^2+(1-z_{i})^2}$ with $(r_i,z_i)$ the extraction points. 
We present in the Figure~\ref{slope} the graph of the slopes $\tilde s_{i}(\sigma)$  for each curve on the Figure~\ref{extrLogH}  ($\sigma=20,80$). For the sake of comparison, we also represent on the same figure the slopes of Figure~\ref{Fextr} corresponding to configuration B1.  

Whereas in configuration B1, the slopes clearly converge to a positive limit value as $\rho=y_3$ tends to $0$, in configuration A the slopes tend to $0$, which means that, {\em stricto sensu}, there is no exponential decay near the corner. Nevertheless we notice that in a region which is further away from the corner, a sort of exponential convergence is restored. This phenomenon is due to the fact that the principal asymptotic contribution inside the conductor is a \emph{profile globally defined on an infinite sector} $\mathcal{S}$ (of opening $\frac\pi2$ in the present case) solving, instead the 1D problem \eqref{Ev0} the model Dirichlet problem 
\begin{equation}
\label{Ev0A}
\left\{
  \begin{array}{lllll}
  (\partial^2_{X} +\partial^2_{Y})\svV_0 - \lambda^2\svV_0 &=& 0 & 
  \mbox{in}\quad \mathcal{S} \ ,
\\ [4pt]
 \svV_0& =&  \sH_{0}^+(\aa) & \mbox{on}\quad \partial\mathcal{S}\ .
   \end{array}
\right.
\end{equation}

\section{Conclusion}
Even though addressing axisymmetric configurations for which the Maxwell system can be reduced to one scalar equation, our numerical experiments are in significative accordance with our theoretical results concerning the decay of solutions inside the conductor and their structure  in the skin layer. Our asymptotics provide an {\em a priori} knowledge on solutions, which can be used for the design of meshes in view of a good quality finite element approximation: The mesh should fit the boundary of the conductor and can be coarse far from its boundary inside the conductor -- depending on the skin depth. In this perspective, it is interesting to compare with \cite{PardoDTP06} where an adaptive {\em a posteriori} approach based on the $hp$ method has been used for an industrial axisymmetric problem in electromagnetism.

\appendix
\section{Elements of proof for the multiscale expansion}
\label{AppA}

	Subsequently, we assume Assumption \ref{H1} on $\omega$ and Assumption \ref{H2} on the domains. In this framework, Theorem \ref{2T0} gives the existence of $\delta_{0}$ such that for all $\delta\leqslant\delta_{0}$, the  problem \eqref{MS}-\eqref{PIbc} has a unique solution $(\EE_{(\delta)}, \,  \HH_{(\delta)})$ which is denoted by $(\EE^+_{(\delta)}, \,  \HH^+_{(\delta)})$ in the dielectric part $\Omega_{+}$, and $(\EE^-_{(\delta)}, \,  \HH^-_{(\delta)})$ in  the conducting part  $\Omega_{-}$.
Furthermore, we suppose that the right hand side $\jj\in\bH_{0}(\Div,\Omega)$ is smooth and its support does not meet the conductor domain $\Omega_{-}$.

Recall that $(y_{\alpha},y_{3})$ is a local {\em normal coordinate system} to the surface $\Sigma$ in  $\cU_-$, see Figure~\ref{Fig1}. The function $\yy\mapsto\chi(y_3)$ is a smooth cut-off with support in $\overline\cU_-$ and equal to $1$ in a smaller tubular neighborhood of $\Sigma$.  
\begin{thm}
Under the above assumptions, the solution $(\EE_{(\delta)}, \,  \HH_{(\delta)})$  possesses the asymptotic expansion (see subsection {\rm\ref{AConv}} below for precise estimates):  
\begin{gather}
\label{EAa}
   \EE^+_{(\delta)}(\xx) \approx \sum_{j\geqslant0} \delta^j\EE^+_j(\xx) 
    \quad\mbox{and}\quad
     \HH^+_{(\delta)}(\xx) \approx \sum_{j\geqslant0} \delta^j\HH^+_j(\xx) \, ,
\\  
\label{EAb}
   \EE^-_{(\delta)}(\xx) \approx \sum_{j\geqslant0} \delta^j\EE^-_j(\xx;\delta) 
   \quad\mbox{with}\quad   
   \EE^-_j(\xx;\delta) = \chi(y_3) \,\WW_j(y_\beta,\frac{y_3}{\delta})\, ,
\\
\label{EAc}
   \HH^-_{(\delta)}(\xx) \approx \sum_{j\geqslant0} \delta^j\HH^-_j(\xx;\delta) 
   \quad\mbox{with}\quad   
   \HH^-_j(\xx;\delta) = \chi(y_3) \,\VV_j(y_\beta,\frac{y_3}{\delta})\, ,
\end{gather}
where $\WW_j(y_\beta,\frac{y_3}{\delta})\rightarrow 0$ and $\VV_j(y_\beta,\frac{y_3}{\delta})\rightarrow 0$ when $\frac{y_3}{\delta}\rightarrow \infty$. Moreover, for any $j\in\N$, there holds
\begin{equation}
\label{EAd}
   \EE^+_j,\HH^+_j\in\bH(\rot,\Omega_+) \quad\mbox{and}\quad
   \WW_j,\VV_{j}\in\bH(\rot,\Sigma\times\R_+).
\end{equation}
\end{thm}

Hereafter, we present elements of proof of this theorem and details about the terms in asymptotics \eqref{EAb}--\eqref{EAc}. In \S\ref{AH1}, we expand the ``magnetic'' Maxwell operators in power series of $\delta$ inside the boundary layer $\cU_-$. 
We deduce in \S\ref{AH2} the equations satisfied by the magnetic profiles, and derive explicitly the first ones in \S\ref{AH3}. Then in \S\ref{A2} and \S\ref{A3}, we do the same for the electric profiles. As an alternative, we show how to deduce directly the magnetic profiles from the electrical ones in \S\ref{A4}. In \S \ref{AConv}, we conclude to the validation of the asymptotic expansion with a convergence result.

\subsection{Expansion of the operators}
\label{AH1}
Integrating by parts in the magnetic variational formulation \eqref{FVH0}, we find the following Maxwell transmission problem
\begin{equation}
\label{EMH}
 \left\{
   \begin{array}{lll}
   \rot\rot  \HH^+_{(\delta)} - \kappa^2 \HH^+_{(\delta)} = \rot\jj   \quad&\mbox{in}\quad \Omega_{+}
\\[0.5ex]
\rot\rot \HH^-_{(\delta)} -\kappa^2({1+\frac{i}{\delta^2}}) \HH^-_{(\delta)} =0
 \quad&\mbox{in}\quad \Omega_{-}
\\[0.5ex]
  \rot\HH^+_{(\delta)}\times\nn=  ({1+\frac{i}{\delta^2}})^{-1}\rot\HH^-_{(\delta)}\times\nn \quad &\mbox{on} \quad \Sigma
\\[0.5ex]
\HH^+_{(\delta)}\times\nn=\HH^-_{(\delta)}\times\nn  \quad &\mbox{on}\quad \Sigma
\\[0.5ex]
  \HH^+_{(\delta)} \times\nn= 0 
  \quad &\mbox{on}\quad \Gamma .
   \end{array}
    \right.
\end{equation}
It is important to notice that, since $\kappa\neq0$, it is a consequence of the above equations that 
\[
   \Div\HH_{(\delta)} = 0 \quad\mbox{in}\quad\Omega.
\]
Therefore, we have in particular the extra transmission condition
\begin{equation}
\label{AHn}
 \HH^+_{(\delta)}\cdot\nn=\HH^-_{(\delta)}\cdot\nn  \quad \mbox{on}\quad \Sigma.
\end{equation}

We denote by $\LL(y_{\alpha}, h;D_{\alpha}, \partial_{3}^h)$ the 2d order Maxwell operator $\rot\rot-\kappa^2({1+\frac{i}{\delta^2}})\Id$ set in $\cU_{-}$ in a {\em normal coordinate system}. Here $D_{\alpha}$ is the covariant derivative on the interface $\Sigma$, and $\partial_{3}^h$ is  the partial derivative with respect to the normal coordinate $y_{3}=h$. 

Let $a_{\alpha\beta}(h)$ be the metric tensor of the manifold $\Sigma_{h}$, which is the surface contained in $\Omega_-$ at a distance $h$ of $\Sigma$, see Figure~\ref{Fig1}.
According to {\rm\cite{DFP09,Fa02}}, the metric tensor in such a coordinate system writes
\begin{equation}
\label{Emett}
a_{\alpha\beta}(h)=a_{\alpha\beta}-2 b_{\alpha\beta} h + b_{\alpha}^{\gamma}b_{\gamma\beta} h^2 \, ,
\end{equation}
and its inverse expands in power series of $h$
\begin{equation*}
a^{\alpha\beta}(h)=a^{\alpha\beta}+2 b^{\alpha\beta} h + \mathcal{O}( h^2) \, .
\end{equation*}
With this metric, a three-dimensional vector field $\VV$ can be split into its normal component $\sv$ and its tangential component that can be alternatively viewed as a vector field $\cV^{\alpha}$ or a one-form field $\cV_{\alpha}$ with the relation
\begin{equation}
\label{Emet}
\cV^{\alpha} = a^{\alpha \beta}(h)\cV_{\beta}. 
\end{equation}\Bk
Subsequently, we use a property of the covariant derivative, that it acts on functions like the partial derivative:  $D_{\alpha} \sv=\partial_{\alpha} \sv$.  

We denote by $\TT(y_{\alpha} ,h; D_{\alpha}, \partial_{3}^h)$ the tangent trace operator $({1+\frac{i}{\delta^2}})^{-1}\rot\cdot\times\nn$ on $\Sigma$. If $\VV=(\cV_{\alpha},\sv)$, then 
\begin{equation}
\label{ETTT}
\TT(y_{\alpha}, h;D_{\alpha} , \partial_{3}^h)\VV=({1+\frac{i}{\delta^2}})^{-1}(\partial_{3}^h\cV_{\alpha} - D_{\alpha} \sv) dy^{\alpha} \ ,
\end{equation}
see \cite[Ch. 6, eq. (6.87)]{Pe09}. The operators $\LL$ and $\TT$ expand in power series of $h$ with intrinsic coefficients with respect to $\Sigma$, see \cite{Fa02} for the 3D elasticity operator on a thin shell. We make the scaling $Y_{3}=\delta^{-1} h$ to describe the boundary layer with respect to $\delta$. Then, the three-dimensional harmonic Maxwell operators in $\cU_{-}$ are written $\LL[\delta]$ and $\TT[\delta]$. These operators expand in power series of $\delta$ with coefficients intrinsic operators :
\begin{equation*}
 \LL[\delta]={\delta}^{-2}\di\sum_{n=0}^{\infty} \delta^n\LL^{n} \quad \mbox{and} \quad  \TT[\delta]=\di\sum_{n=1}^{\infty} \delta^n\TT^{n} \ .
\end{equation*}
We denote by $L_{\alpha}^n$ and $T_{\alpha}^n$ the surface components of $\LL^{n}$ and $\TT^{n}$. With the summation convention of repeated two dimensional indices (represented by greek letters), there holds
\begin{equation}
L_{\alpha}^0(\VV)=-\partial_{3}^2 \cV_{\alpha}-i\kappa^2 \cV_{\alpha} 
  \ \  \mbox{and}  \ \ 
L_{\alpha}^1(\VV)=-2b_{\alpha}^\beta\partial_{3} \cV_{\beta} 
+ \partial_{3} D_{\alpha} \sv 
+  b_{\beta}^{\beta}\partial_{3}\cV_{\alpha}\ ,
\label{EHsL0L1}
\end{equation}
and
\begin{equation}
T_{\alpha}^n(\VV)=
\left\{
   \begin{array}{lll}
(-i)^p\partial_{3} \cV_{\alpha}  \quad &\mbox{if} \quad n=2p-1
\\[5pt]
i^p\partial_{\alpha}\sv   \quad &\mbox{if} \quad n=2p \ .
  \end{array}
   \right.
\label{EHTn}
\end{equation}
Here, $\partial_{3}$ is the partial derivative with respect to $Y_{3}$. We denote by $L_{3}^n$ the transverse components of $\LL^{n}$. There holds
\begin{equation}
L_{3}^0(\VV)=-i\kappa^2 \sv 
  \quad \mbox{and}  \quad
L_{3}^1(\VV)=\gamma_{\alpha}^{\alpha}(\partial_{3}\VV)+  b_{\beta}^{\beta}\partial_{3} \sv\ ,
\label{EHtL0L1}
\end{equation}
where $\gamma_{\alpha\beta}(\VV)=\frac12 (D_{\alpha}\cV_\beta + D_{\beta}\cV_\alpha)-b_{\alpha\beta}\sv$ is the change of metric tensor.

\subsection{Equations for the coefficients of the magnetic field}
\label{AH2}

The profiles $\VV_{j}$ and the coefficients $\HH_{j}^+$ of the magnetic field satisfy the following system
\begin{gather}
\label{ELV}
 \LL[\delta] \sum_{j\geqslant0} \delta^j\VV_j(y_{\alpha},Y_{3})=0  \quad \mbox{in}\quad \Sigma \times I,
 \\
\label{ETV}
 \TT[\delta] \sum_{j\geqslant0} \delta^j\VV_j(y_{\alpha},0)= \sum_{j\geqslant0} \delta^j \rot \HH_{j}^+ \times\nn \quad \mbox{on} \quad \Sigma \ ,
\end{gather}
where $I = (0,+\infty)$.
We perform in \eqref{ELV}-\eqref{ETV} the identification of terms with the same power in $\delta$.
The components of equation \eqref{ELV} are the collections of equations
\begin{equation}
\label{EsLVj}
\LL^0(\VV_{0})=0\ , \quad \LL^0(\VV_{1})+\LL^1(\VV_{0})=0 \ , \quad \mbox{and} \quad \di\sum_{l=0}^{n} \LL^{n-l}(\VV_{l})=0  \ ,
\end{equation}
for all $n\geqslant2$. Similarly, the surface components of the equation \eqref{ETV} write 
\begin{equation}
\label{EsTVj}
\rot \HH_{0}^+ \times\nn  =0\ ,  \quad \mbox{and} \quad 
\di\sum_{k=1}^n \;  \TT^{k}\VV_{n-k} = \rot\HH_{n}^+\times\nn  \ ,
\end{equation}
for all $n\geqslant 1$. 
Using the expression of the operator $\LL^0$, and expanding $\HH^+_{(\delta)}$ in $\Omega_+$, we thus see that, according to the system \eqref{EMH}, the profiles 
$\VV_n = (\cV_{n},\sv_n)$ and the terms $\HH^+_{n}$ have to satisfy, for all $n \geq 0$, 
\begin{equation}
\label{EH}
 \left\{
   \begin{array}{clll}
   (i) & \quad &
   -\lambda^2 \sv_n =  \di \sum_{j = 0}^{n-1} L_3^{n-j}( \VV_j) 
 \quad&\mbox{in}\quad \Sigma \times I
\\[0.5ex]
   (ii) & \quad &
   \rot\rot  \HH^+_{n} - \kappa^2 \HH^+_{n} = \delta_n^0 \rot\jj   \quad&\mbox{in}\quad \Omega_{+}
\\[0.5ex]
   (iii) & \quad &
    \rot\HH_{n}^+\times\nn  = \di\sum_{j = 0}^{n-1} \;  \TT^{n-j}\VV_{j}\quad &\mbox{on} \quad \Sigma
\\[0.5ex]
   (iv) & \quad &
  \HH^+_{n} \times\nn= 0 \quad 
  &\mbox{on}\quad \Gamma 
\\[0.5ex]
   (v) & \quad &
\partial_{3}^2 \cV_{n,\alpha}- \lambda^2 \cV_{n,\alpha} =  \di \sum_{j = 0}^{n-1} L_\alpha^{n-j}( \VV_j) 
 \quad&\mbox{in}\quad \Sigma \times I
 \\[0.5ex]
   (vi) & \quad &
 \cV_{n} =\nn\times \bsH_n\times\nn &\mbox{on} \quad \Sigma   .
   \end{array}
    \right.
\end{equation}
where $\lambda=\kappa\, \mathrm{e}^{-i\pi/4}$, cf \eqref{V0cd} (so that  $-\lambda^2=i\kappa^2$) and $\bsH_n$ denotes the trace of $\HH^+_{n}$ on $\Sigma$. In \eqref{EH}, we use the convention that the sums are $0$ when $n=0$.  The transmission condition \eqref{AHn} implies the extra continuity condition 
\begin{equation}
\label{EH1}
\HH^+_{n} \cdot  \nn = \sv_n \quad \mbox{on}\quad \Sigma. 
\end{equation}
The set of equations \eqref{EH}--\eqref{EH1} allows to determine $\VV_n $ and $\HH^+_{n}$ by induction.

\subsection{First terms of the magnetic field asymptotics}
\label{AH3}
According to equation $(i)$ in \eqref{EH}, the normal component $\sv_0$ of the first profile in the conductor vanishes:
\begin{equation}
\label{Esv0}
\sv_{0}=0\ , 
\end{equation}
because $\kappa\neq 0$, thus $\lambda\neq0$.

Hence, according to \eqref{EH} $(ii)$-$(iv)$ and \eqref{EH1}, the first term of the magnetic field in the dielectric region solves Maxwell equations with perfectly conducting conditions on $\Sigma$:
\begin{equation}
\label{H0pbis}
 \left\{
   \begin{array}{lll}
    \rot\rot \HH^+_0 - \kappa^2\HH^+_0 = \rot\jj  \quad&\mbox{in}\quad \Omega_{+}
\\[0.5ex]
   \HH^+_0\cdot\nn=0  \quad\mbox{and}\quad\!
   \rot\HH^+_0\times\nn=0 \quad &\mbox{on}\quad \Sigma
\\[0.5ex]
   \HH^+_0\times\nn=0 
   \quad &\mbox{on}\quad \Gamma.
   \end{array}
    \right.
\end{equation}
Thus the trace $\bsH_0$ of $\HH^+_0$ on the interface $\Sigma$ is \emph{tangential}.

According to equations $(v)$-$(vi)$ in  \eqref{EH}, $\cV_{0}$ satisfies the following ODE 
\begin{equation}
\left\{
   \begin{array}{lll}
\partial_{3}^2\cV_{0}(.,Y_{3}) -\lambda^2 \cV_{0}(.,Y_{3})&=0 
&\mbox{for} \quad Y_{3}\in I=(0,\infty) 
\\[2pt]
\cV_{0}(.,0)& = (\nn\times\bsH_0)\times\nn \ .
\end{array}
  \right.
\label{EV0}
\end{equation}
The unique solution of \eqref{EV0} such that $\cV_{0}\rightarrow 0$ when $Y_{3}\rightarrow\infty$, is, with the choice \eqref{V0cd} for $\lambda$, the tangential field
$\cV_{0}(y_{\beta},Y_{3}) = \bsH_{0} (y_{\beta})\,\mathrm{e}^{-\lambda Y_{3}}$.
Combining with \eqref{Esv0}, we find that the first profile in the conductor region is exponential with the complex rate $\lambda$:
\begin{equation}
\VV_{0}(y_{\beta},Y_{3}) = \bsH_{0} (y_{\beta})\,\mathrm{e}^{-\lambda Y_{3}} \, .
\label{VOcdbis}
\end{equation}

The next term which is determined in the asymptotics is the normal component $\sv_1$ of the profile $\VV_1$ given by equation $(i)$ of \eqref{EH} for $n = 1$. We obtain 
\begin{equation}
\label{v1cd}
\sv_{1}(y_{\beta},Y_{3}) =\lambda^{-1} D_{\alpha}\sH^{\alpha}_{0}(y_{\beta}) \;\mathrm{e}^{-\lambda Y_{3}}\,.
\end{equation}
According to \eqref{EH} $(ii)$-$(iv)$ and \eqref{EH1}, the next term in the dielectric region solves:
\begin{equation}
\label{H1pbis}
 \left\{
   \begin{array}{lll}
    \rot\rot \HH^+_1 - \kappa^2\HH^+_1 = 0  \quad&\mbox{in}\quad \Omega_{+}
\\[0.5ex]
  \HH^+_1\cdot\nn= \sv_1
    \quad\mbox{and}\quad\!
   \rot\HH^+_1\times\nn= i \lambda \bsH_0 \quad &\mbox{on}\quad \Sigma
\\[0.5ex]
   \HH^+_1\times\nn=0 
   \quad &\mbox{on}\quad \Gamma.
   \end{array}
    \right.
\end{equation}

Recall that $\bsH_1$ is the trace of $\HH^+_1$ on the interface $\Sigma$. We denote by $\sH_{1,\alpha}$ its tangential components. According to equations $(v)$-$(vi)$ in \eqref{EH} for $n = 1$, $\cV_{1}$ satisfies the following ODE (for $Y_3\in I$)
\begin{equation}
\label{EVb1}
\left\{
   \begin{array}{lll}
   \partial_{3}^2\cV_{1,\alpha}(.\,,Y_{3}) -\lambda^2 \cV_{1,\alpha}(.\,,Y_{3}) = -2b_{\alpha}^{\sigma}\partial_{3}\cV_{0,\sigma}(.\,,Y_{3})+b_{\beta}^{\beta}\partial_{3}\cV_{0,\alpha}(.\,,Y_{3})
\\[2pt]
\cV_{1,\alpha}(.\,,0)=\sH_{1,\alpha}(.\,,0)\ .
   \end{array}
\right.
\end{equation} 
From \eqref{VOcdbis}, the unique solution of \eqref{EVb1} such that $\cV_{1}\rightarrow 0$ when $Y_{3}\rightarrow\infty$ is the profile
\begin{equation}
\label{cV1cd}
   {\cV_{1,\alpha}}(y_{\beta},Y_{3}) =
   \Big[ \sH_{1,\alpha}+Y_{3} \big(\cH\,\sH_{0,\alpha}
   -b^{\sigma}_{\alpha} \sH_{0,\sigma}\big)\Big](y_{\beta}) \;\mathrm{e}^{-\lambda Y_{3}},
   \quad\alpha=1,2\,. 
\end{equation}    
Using the relation \eqref{Emet} and  performing the scaling $h=\delta Y_{3}$ in the previous equation, we obtain for the contravariant components
\begin{equation*}
 {\cV^{\alpha}_{1}}= a^{\alpha\beta} {\cV_{1,\beta}} + 2 Y_{3}   b^{\alpha\beta}   {\cV_{0,\beta}}\ .
\end{equation*}
From \eqref{VOcdbis}, and \eqref{cV1cd}, the tangential components ${\cV}^{\alpha}_1$ are given by \eqref{V1cd}
\begin{equation*}
   {\cV^{\alpha}_{1}}(y_{\beta},Y_{3}) =
   \Big[ \sH^{\alpha}_{1}+Y_{3} \big(\cH\,\sH^{\alpha}_{0}
   +b_{\sigma}^{\alpha} \sH^{\sigma}_{0}\big)\Big](y_{\beta}) \;\mathrm{e}^{-\lambda Y_{3}},
   \quad\alpha=1,2\,. 
\end{equation*}    
 
\begin{rem}
Note that the boundary value problems \eqref{H0pbis} and \eqref{H1pbis}  are well-posed. It is a consequence of the spectral Assumption \ref{H1} on $\omega$. 
\end{rem}

\subsection{Equations for the coefficients of the electric field}
\label{A2}

The second order Maxwell operator for the electric field writes 
\begin{equation}
\label{EME}
 \left\{
   \begin{array}{lll}
   \rot\rot  \EE^+_{(\delta)} - \kappa^2 \EE^+_{(\delta)} =i \omega\mu_{0}\jj   \quad&\mbox{in}\quad \Omega_{+}
\\[0.5ex]
\rot\rot  \EE^-_{(\delta)} -\kappa^2({1+\frac{i}{\delta^2}})  \EE^-_{(\delta)} =0
 \quad&\mbox{in}\quad \Omega_{-}
\\[0.5ex]
  \rot \EE^+_{(\delta)}\times\nn=  \rot \EE^-_{(\delta)}\times\nn \quad &\mbox{on} \quad \Sigma
\\[0.5ex]
 \EE^+_{(\delta)}\times\nn= \EE^-_{(\delta)}\times\nn  \quad &\mbox{on}\quad \Sigma
\\[0.5ex]
  \EE^+_{(\delta)} \cdot\nn= 0 \quad\mbox{and}\quad\!  \rot \EE^+_{(\delta)}\times\nn = 0  \quad &\mbox{on}\quad \Gamma .
   \end{array}
    \right.
\end{equation}
We denote by $\BB(y_{\alpha}, h;D_{\alpha}, \partial_{3}^h)$ the tangent trace operator $\rot\cdot\times\nn$ on $\Sigma$ in a {\em normal coordinate system}. If 
$\WW=(\cW_{\alpha},\fke)$, then
\begin{equation*}
\left( \BB(y_{\alpha}, h;D_{\alpha} , \partial_{3}^h)\WW\right)_{\alpha}=\partial_{3}^h \cW_{\alpha} - D_{\alpha} \fke \ ,
\end{equation*}
see \cite[Ch. 3, Prop. 3.36]{Pe09}. We define $\BB[\delta]$ the operator  obtained from $\BB$ in $\cU_{-}$ after the scaling $Y_{3}=\delta^{-1} h$. This operator expands in power of $\delta$ :
\begin{equation*}
 \BB[\delta]={\delta}^{-1}\BB^{0}+\BB^{1} \ .
\end{equation*}
Recall that $\partial_{3}$ is the partial derivative with respect to $Y_{3}$. Thus, denoting by $B_{\alpha}^n$ the surface components of $\BB^{n}$,  we obtain
\begin{equation}
\label{EB0B1}
B_{\alpha}^0(\WW)= \partial_{3}\cW_{\alpha}   \quad \mbox{and}  \quad B_{\alpha}^1(\WW)= -D_{\alpha} \fke \ .
\end{equation}

According to the second and third equations in system \eqref{EME}, the profiles $\WW_{j}$ and the terms $ \EE_{j}^+$ of the electric field satisfy the following system
\begin{gather}
\label{ELW}
 \LL[\delta] \sum_{j\geqslant0} \delta^j\WW_j(y_{\alpha},Y_{3})=0  \ ,
   \quad \mbox{in }  \quad \Sigma \times I \ ,
   \\
  \label{EBW}
 \BB[\delta] \sum_{j\geqslant0} \delta^j\WW_j(y_{\alpha},0)= \sum_{j\geqslant0} \delta^j \rot \EE_{j}^+ \times\nn \quad \mbox{on} \quad \Sigma \ .
\end{gather}

We identify in \eqref{ELW}-\eqref{EBW} the terms with the same power in $\delta$. The components of the equation \eqref{ELW} are collections of equations, similar to the equations \eqref{EsLVj} set for the magnetic field. 
The surface components of the equation \eqref{EBW} write 
\begin{equation}
\label{EsBWj}
B_{\alpha}^0(\WW_{0})=0\ ,  \quad \mbox{and} \quad B_{\alpha}^0(\WW_{n+1})+B_{\alpha}^1(\WW_{n}) = \left(\rot \EE_{n}^+ \times\nn \right)_{\alpha}   \ ,
\end{equation}
for all $n\geqslant0$. 

According to the system \eqref{EME} and \eqref{EB0B1}, the profiles 
$\WW_n = (\cW_{n},\fke_n)$ and the terms $\EE^+_{n}$ have to satisfy, for all $n \geqslant 0$, 
(we recall $\lambda=\kappa\, \mathrm{e}^{-i\pi/4}$, cf \eqref{V0cd}, so that  $-\lambda^2=i\kappa^2$)
\begin{equation}
\label{EE}
 \left\{
   \begin{array}{clll}
   (i) & \quad &
   \partial_{3}^2 \cW_{n,\alpha}-\lambda^2 \cW_{n,\alpha} =  \di \sum_{j = 0}^{n-1} L_\alpha^{n-j}( \WW_j) 
 \quad&\mbox{in}\quad \Sigma \times I
  \\[0.5ex]
   (ii) & \quad &
   \partial_{3}\cW_{n,\alpha}  = D_{\alpha} \fke_{n-1}+ \left( \rot \EE_{n-1}^+ \times\nn \right)_{\alpha}
  \quad &\mbox{on} \quad \Sigma
\\[0.5ex]
   (iii) & \quad &
   -\lambda^2 \fke_n =  \di \sum_{j = 0}^{n-1} L_3^{n-j}( \WW_j) 
 \quad&\mbox{in}\quad \Sigma \times I
\\[0.5ex]
   (iv) & \quad &
   \rot\rot  \EE^+_{n} - \kappa^2 \EE^+_{n} = \delta_n^0 i \omega\mu_{0}\jj   \quad&\mbox{in}\quad \Omega_{+}
\\[0.5ex]
   (v) & \quad &
\EE^+_{n}\times\nn=\cW_{n}\times\nn \quad &\mbox{on}\quad \Sigma
\\[0.5ex]
   (vi) & \quad &
  \EE^+_{n} \cdot\nn= 0 \quad\mbox{and}\quad\!  \rot \EE^+_{n}\times\nn = 0 \quad &\mbox{on}\quad \Gamma .
   \end{array}
    \right.
\end{equation}
Hereafter, we determine the terms $\WW_{n}$ and $\EE^+_{n}$ by induction.

\subsection{First terms of the asymptotics for the electric field}
\label{A3}
According to equations $(i)$-$(ii)$ in system \eqref{EE} for $n=0$, $\cW_{0}$ satisfies the following ODE 
\begin{equation}
\left\{
   \begin{array}{ll}
\partial_{3}^2\cW_{0}(.,Y_{3}) -\lambda^2 \cW_{0}(.,Y_{3})&=0 \quad \mbox{for} \quad Y_{3}\in I \ ,
\\[0.2ex]
\partial_{3}\cW_{0}(.,0)&=0 \ .
\end{array}
  \right.
\label{EW0}
\end{equation}
The unique solution of \eqref{EW0} such that $\cW_{0}\rightarrow 0$ when $Y_{3}\rightarrow\infty$, is $\cW_{0}=0$. From equation $(iii)$ in system \eqref{EE}  for $n=0$, there holds $\fke_{0}=0$. We infer 
\begin{equation}
\WW_{0}(y_{\beta},Y_{3}) =0 \, .
\label{W0cd}
\end{equation}
From equations $(iv)$-$(vi)$ in \eqref{EE} for $n=0$, and from \eqref{W0cd}, the first asymptotic of the electric field in the dielectric part solves the following problem 
\begin{equation*}
 \left\{
   \begin{array}{lll}
    \rot\rot  \EE^+_{0} - \kappa^2 \EE^+_{0} =  i \omega\mu_{0}\jj   \quad&\mbox{in}\quad \Omega_{+}
\\[0.2ex]
\EE^+_{0}\times\nn=0 \quad &\mbox{on}\quad \Sigma
\\[0.2ex]
  \EE^+_{0} \cdot\nn= 0 \quad\mbox{and}\quad\!  \rot \EE^+_{0}\times\nn = 0 \quad &\mbox{on}\quad \Gamma .
   \end{array}
    \right.
\end{equation*}
According to the spectral Assumption  \ref{H1}, this boundary value problem is well-posed. 

The next term determined in the asymptotic expansion is $\cW_{1}$. From equations $(i)$-$(ii)$ in \eqref{EE}  for $n=1$, $\cW_{1}$ satisfies for $Y_3\in I$
\begin{equation*}
\left\{
   \begin{array}{lll}
\partial_{3}^2\cW_{1}(.,Y_{3}) -\lambda^2 \cW_{1}(.,Y_{3})&=0 
\\[0.2ex]
\partial_{3}\cW_{1}(.,0)&= (\rot{\EE}_{0}\is \times\nn)(.\,,0)    \ .
\end{array}
  \right.
\end{equation*}
We denote by $\jj_{k}(y_{\beta})=\lambda^{-1} (\rot{\EE}_{k}\is \times\nn)(y_{\beta},0)$ for $k=0,1$. Hence,
\begin{equation}
\label{W1cd}
\cW_{1}(y_{\beta},Y_{3}) = -\jj_{0}(y_{\beta})\;\mathrm{e}^{-\lambda Y_{3}}
\, .
\end{equation}
From equation $(iii)$ in \eqref{EE}  for $n=1$, and from \eqref{W0cd}, we obtain $\fke_{1}=0$. From equations $(iv)$-$(vi)$ in \eqref{EE} for $n=1$, and from \eqref{W1cd}, the asymptotic of order $1$ for the electric field in the dielectric part solves :
\begin{equation*}
 \left\{
   \begin{array}{lll}
    \rot\rot  \EE^+_{1} - \kappa^2 \EE^+_{1} =  0   \quad&\mbox{in}\quad \Omega_{+}
\\[0.2ex]
\EE^+_{1}\times\nn= -\jj_{0} \times\nn \quad &\mbox{on}\quad \Sigma
\\[0.2ex]
  \EE^+_{1} \cdot\nn= 0 \quad\mbox{and}\quad\!  \rot \EE^+_{1}\times\nn = 0 \quad &\mbox{on}\quad \Gamma .
   \end{array}
    \right.
\end{equation*}

 Then, from equations $(i)$-$(ii)$ in \eqref{EE} for $n=2$, $\cW_{2,\alpha}$ solves the ODE for $Y_3\in I$:
\begin{equation*}
\left\{
   \begin{array}{lll}
   \partial_{3}^2\cW_{2,\alpha}(.\,,Y_{3}) -\lambda^2 \cW_{2,\alpha}(.\,,Y_{3}) & = -2b_{\alpha}^{\sigma}\partial_{3}\cW_{1,\sigma}(.\,,Y_{3})+b_{\beta}^{\beta}\partial_{3}\cW_{1,\alpha}(.\,,Y_{3}) 
\\[0.2ex]
\partial_{3}\cW_{2,\alpha}(.\,,0) & = (\rot{\EE}_{1}\is \times\nn)_{\alpha}(.\,,0)\ .
   \end{array}
\right.
\end{equation*}
We denote by $\sj_{k,\alpha}$ the surface components of $\jj_{k}$, for $k=0,1$. We obtain
\begin{gather}
\label{W2cd}
   \cW_{2,\alpha}(y_{\beta},Y_{3})=
   \Big[ -\sj_{1,\alpha} + \big(\lambda^{-1} + Y_{3}\big) 
   \big( b_{\alpha}^{\sigma}\ \sj_{0,\sigma}-\cH\,\sj_{0,\alpha}\big) 
   \Big](y_{\beta}) \;\,\mathrm{e}^{-\lambda Y_{3}}\, .
 \end{gather}
From equation $(iii)$ in \eqref{EE} for $n=2$, we infer
$$
\fke_{2}(y_{\beta},Y_{3})= 
   -\lambda^{-1} D_{\alpha}\ \sj_{0}^{\alpha}(y_{\beta})\,\mathrm{e}^{-\lambda Y_{3}}  .
$$

\subsection{From E to H in the conducting part}
\label{A4}
An alternative way of calculating the magnetic profiles is to deduce them from the electric ones by means of a {\em normal parameterization} of the intrinsic curl operator, see \cite[Ch.\ 3]{Pe09}. Let $\WW=(\cW_{\alpha},\fke)$ be a 1-form fields in $\cU_-$. There holds: 
\begin{equation}
 \label{rot}
(\nabla\times\WW)^\alpha=\epsilon^{3\beta\alpha}(\partial_{3}^{h}\cW_{\beta}-
\partial_{\beta}\fke)
\quad \mbox{and} \quad 
(\nabla\times\WW)^3=\epsilon^{3\alpha\beta}D_{\alpha}^h \cW_{\beta} 
\quad \mbox{on} \quad \Sigma_{h} \, .
\end{equation}
Here, $D_{\alpha}^h$ is the {\em covariant derivative} on $\Sigma_{h}$, and $\epsilon$ is the Levi-Civita tensor, see \cite{Le26,Go80}. The contravariant components $\epsilon^{ijk}$ of $\epsilon$ depend on the normal coordinate $h$, and write in a \emph{normal coordinate} system 
$$
\epsilon^{ijk}=\big(\det a_{\alpha\beta}(h)\big)^{-1/2}\,\epsilon_{0}(i,j,k)\,.
$$
Here, $a_{\alpha\beta}(h)$ is the metric tensor of the manifold $\Sigma_{h}$, see \eqref{Emett}. The indices $i,j,k\in\{1,2,3\}$, and  $\epsilon_{0}(i,j,k)$ equals $1$ when $(i,j,k)$ is an even circular permutation, and equals $-1$ when $(i,j,k)$ is an odd circular permutation, and $\epsilon_{0}(i,j,k)=0$ otherwise.

\begin{rem}
Let $a=\det(a_{\alpha\beta})$, and recall that $b_{\nu}^{\nu}=2 \cH$. Then using \eqref{Emett} we obtain
$$
\epsilon^{ijk}= a^{-1/2}\big(1+2 \cH h+O(h^2)\big) \,\epsilon_{0}(i,j,k)
\,.
$$
\end{rem} 

We make the scaling $h=\delta Y_{3}$ and expand equations \eqref{rot} in power series of $\delta$: 
\begin{equation}
\label{6E27}
(\nabla\times\WW)^\alpha=
\delta^{-1} j^{\beta\alpha}\partial_{3} \cW_{\beta}+j^{\beta\alpha}\big(2 \cH Y_{3}\partial_{3} \cW_{\beta}-\partial_{\beta}\fke \big)+\mathcal{O}(\delta)\, ,
\end{equation}
\begin{equation}
\label{6E26}
 (\nabla\times\WW)^3 =j^{\alpha\beta}D_{\alpha}\cW_{\beta}
+\delta Y_{3} j^{\alpha\beta}(2 \cH  D_{\alpha}\cW_{\beta}+\cW_{\nu}D_{\alpha}b_{\beta}^{\nu})+\mathcal{O}(\delta^2)\,.
\end{equation}
Here, 
\begin{equation*}
j^{\alpha\beta}=a^{-1/2}\epsilon_{0}(\alpha,\beta,3)\,.
\end{equation*}

From the expansions \eqref{EAb} and \eqref{EAa} combined with Faraday's law \eqref{MS} written in normal coordinates, we obtain the \textit{profiles} $\VV_{j}$ in the expansion \eqref{EAc} of the magnetic field from the profiles $\cW_j$ of the electrical field. In particular, from \eqref{W0cd},  \eqref{W1cd}, and  \eqref{W2cd}, we obtain explicitly the first terms $\VV_{0}$ and $\VV_{1}$, see \eqref{VOcdbis}, \eqref{V1cd}, and \eqref{v1cd}.

\subsection{Convergence result}
\label{AConv}
The validation of the asymptotic expansion \eqref{EAa}, \eqref{EAb},\eqref{EAc},  and  \eqref{EAd}, consists in proving estimates for remainders defined as 
\begin{gather}
\label{EA1}
   \RR_{m;\,\delta}^{\EE} = \EE_{(\delta)} - \sum_{j=0}^m \delta^j\EE_j \quad\mbox{and}\quad   
  \RR_{m;\,\delta}^{\HH} = \HH_{(\delta)} - \sum_{j=0}^m \delta^j\HH_j \quad\mbox{in}\quad\Omega\ .
\end{gather}

By construction of the terms $\EE_j=(\EE_j^+,\EE_j^-)$ in the dielectric and conductor parts, the remainders $\RR_{m;\,\delta}^{\EE}$ satisfy the assumption of \cite[Th. 5.1]{CDP09}, which is an estimate on the right hand side when the Maxwell operator is applied to $\RR_{m;\,\delta}^{\EE}$.

Thus \cite[Th. 5.1]{CDP09} yields that for all $m\in\N$, and $\delta\in(0,\delta_{0})$, there holds the optimal estimate 
\begin{multline}
\label{6E8}
   \qquad \|\RR_{m;\,\delta}^{\EE,+} \|_{0,\Omega_+} 
   + \|\rot\RR_{m;\,\delta}^{\EE,+} \|_{0,\Omega_+} \\ + 
   \delta^{-\frac12} \|\RR_{m;\,\delta}^{\EE,-} \|_{0,\Omega_-} + 
   \delta^{\frac12} \|\rot\RR_{m;\,\delta}^{\EE,-} \|_{0,\Omega_-}
    \leqslant C_m\delta^{m+1}.\qquad
\end{multline}
Using Maxwell equations \eqref{MS}, we can deduce a similar estimate for $\RR_{m;\,\delta}^{\HH}$: We have
\begin{equation}
\label{ARH+}
 \RR_{m;\,\delta}^{\HH,+}= (i \omega\mu_0 )^{-1} \rot \RR_{m;\,\delta}^{\EE,+}  \quad \mbox{and} \quad \rot\RR_{m;\,\delta}^{\HH,+}=(i \omega\mu_0 )^{-1} \kappa^2 \RR_{m;\,\delta}^{\EE,+} \ .
\end{equation}
Hence, according to \eqref{6E8}, we infer
\begin{equation*}
   \|\RR_{m;\,\delta}^{\HH,+} \|_{0,\Omega_+} + \|\rot\RR_{m;\,\delta}^{\HH,+} \|_{0,\Omega_+}  \leqslant c_m\delta^{m+1}.
\end{equation*}
In $\Omega_{-}$, the presence of profiles in the expansion prevents to link $\RR_{m;\,\delta}^{\HH,-}$ and $ \RR_{m;\,\delta}^{\EE,-}$ via the Maxwell equations in a similar way as \eqref{ARH+}. Nevertheless, there holds also uniform estimates for $\RR_{m;\,\delta}^{\HH,-}$ :
\begin{equation*} 
   \delta^{-\frac12} \|\RR_{m;\,\delta}^{\HH,-} \|_{0,\Omega_-} + 
   \delta^{\frac12} \|\rot\RR_{m;\,\delta}^{\HH,-} \|_{0,\Omega_-} \leqslant C'_m\delta^{m+1}.
\end{equation*}

\section{The multiscale expansion of the orthoradial component}
\label{AppB}
We denote by $(\vec{e_{r}},\vec{e_{\theta}},\vec{e_{z}})$ the basis associated with the cylindric coordinates $(r,\theta,z)$. 

\subsection{Expansion of the operators}
\label{AppB1}
In the basis $(\vec{e_{r}},\vec{e_{z}})$, recall that $\big(r(\xi),z(\xi)\big)=\XX(\xi)$, $\xi \in (0,L)$ is an \textit{arc-length coordinate} on the interface $\Sigma^\m$, and $(\xi,h=y_{3})$ is the associate normal coordinate system, see \S \ref{AEh}. The normal vector $\nn(\xi)$ at the point $\XX(\xi)$ writes
$$
\nn(\xi)=\big( -z'(\xi),r'(\xi)\big)\ .
$$
Hence, the tubular neighborhood $\cU^\m_-$ of $\Sigma^\m$ inside $\Omega^\m_-$ is represented  thanks to the parameterization 
\begin{equation*}
\Psi :  (\xi,h)\longmapsto (r,z) \ ,
\end{equation*}
\begin{equation}
\label{EB1}
\mbox{where}\quad r=r(\xi)-h z'(\xi) \ , \quad \mbox{and}\quad z=z(\xi)+h  r'(\xi)\ . 
\end{equation}
\begin{rem}
The curvature $k(\xi)$ at the point $\XX(\xi)$ is  defined by 
\begin{equation}
\label{Ek}
k(\xi)=(r'z^{\prime\prime}-z'r^{\prime\prime})(\xi) \ .
\end{equation}
 For $h_{0}<1/\| k\|_{\infty}$ the change of coordinates $\Psi$ 
is a $\mathcal{C}^{\infty}$-diffeomorphism from the cylinder  $\T_L \times [0,h_0)$ into $\cU^\m_-$ :
$$
\cU^\m_-=\Psi \big(\T_L \times [0,h_0) \big)\ .
$$
\end{rem}
Thanks to \eqref{EB1}, we obtain 
$$
\left(\begin{array}{ccc}  \partial_{r} 
\\
 \partial_{z}
\end{array}\right)
=\left(1-h k(\xi)\right)^{-1}
\left(\begin{array}{ccc} r' & -(z'+h r^{\prime\prime} ) 
\\
z' & (r'-h z^{\prime\prime}) \end{array}\right)
\left(\begin{array}{ccc}  \partial_{\xi} 
\\
 \partial_{h}
\end{array}\right) \, .
$$
To perform the formal expansion of $\sH^-_\delta$ solution of \eqref{FVHtheta} in $\Omega^\m_-$, we first use
the change of variables $\Psi$ in  order to write the equations in the
cylinder $\T_L \times [0,h_0)$.  We then perform the rescaling
\begin{equation}
\label{EY}
Y= \delta ^{-1} h
\end{equation}
in the  equations set  in $\Omega^\m_{-}$ and  $\Sigma^\m$ in order  to make
appear the  small parameter $\delta$  in the equations. Actually
$\delta$  appears in  the  equations set  in  $\Omega_{-}$ through  the
expression  of  the operator  $(1+\frac{i}{\delta^2})^{-1} \D$, where $\D$ is defined by \eqref{EOpDB}. We obtain the formal expansion:
$$
\D=  \delta^{-2}\big[ \partial^2_{Y}+ \delta \D_{1} +\delta^2 {\mathsf R}_{\delta} \big] \ ,
$$
where $\D_{1}(\xi,Y;\partial_{\xi}, \partial_Y)=-\left(k(\xi)+\frac{z'}{r}(\xi)\right)\partial_{Y}$, and ${\mathsf R}_{\delta}$
is an  operator, which has smooth coefficients in $\xi$ and $Y$, bounded in $\delta$. Hence, we obtain
$$
i(1+\frac{i}{\delta^2})^{-1} \D+\kappa^2\Id= \A_{0}+ \delta \A_{1}+\delta^2 \Q_{\delta}\ ,
$$
with $\A_{0}=\partial^2_{Y} -\lambda^2\Id$ and $\A_{1}=-\left(k+\frac{z'}{r}\right)(\xi)\,\partial_{Y}$.
Similarly, there holds
$$
\B =\delta^{-1}\partial_{Y}- \frac{z'}{r}(\xi) 
$$
on the interface $\Sigma^{m}$, and 
$$
(1+\frac{i}{\delta^2})^{-1} \B= -i\delta\partial_{Y}+\delta^2 \P_{\delta}\ .
$$
Note that these expansions correspond to the expansions \eqref{ETTT} in orthoradial symmetry.

\subsection{From the 3D asymptotic expansion to the 1D expansion}
\label{AppB2}
In this subsection, we show that we can obtain the profiles $\svV_0$, and $\svV_{1}$ defined by \eqref{v0} and \eqref{v1} from the general profiles $\cV_0^\alpha$ and $\cV_1^\alpha$ defined above. 

We set $(y_\alpha):=(\xi,\theta)$ a coordinate system on $\Sigma$.
Thus, the {\it normal coordinate} system $(y_{\alpha}, h)$  on $\cU_{-}$ is induced by the {\it normal coordinate} system $(\xi,h)$ on $\cU_{-}^\m$. The tubular neighborhood $\cU_{-}$ is parameterized by 
\begin{equation*}
\Phi :  (y_{\alpha}, h)\longmapsto (r\cos\theta,r\sin\theta,z) \ ,
\end{equation*}
where $r$ and $z$ are defined by  \eqref{EB1}. The associated tangent coordinate vector fields are 
\begin{gather}
\label{Exx1}
\xx_{1}(h)=  \left(r'(\xi)-h z^{\prime\prime}(\xi)\right) \vec{e_{r}}
+ \left(z'(\xi)+h r^{\prime\prime}(\xi)\right) \vec{e_{z}}\ ,
\\
\label{Exx2}
\xx_{2}(h)= \left(r(\xi)-h z'(\xi)\right) \vec{e_{\theta}}\ .
\end{gather}
The normal coordinate vector field is $\xx_{3}(h)=\nn(\xi)= -z'(\xi) \vec{e_{r}}+r'(\xi) \vec{e_{z}}$\ .

A vector field $\VV: (y_\alpha,h)\mapsto\VV(y_\alpha,h)$ in $\cU_{-}$ writes (here $\alpha = 1,2$)  
$$
\VV= \cV^{\alpha}\xx_{\alpha} (h)+\svV \nn \ ,
$$
and corresponds to the magnetic field in {\it cylindric} components 
$$
\brH(r,\theta,z) = H_{r}(r,\theta,z)\vec{e_{r}} + H_{\theta}(r,\theta,z)\vec{e_{\theta}}+H_{z}(r,\theta,z)\vec{e_{z}}. 
$$
Using \eqref{Exx2} we easily obtain
\begin{equation*}
\left(r(\xi)-h z'(\xi)\right)  \cV^{2}=H_{\theta} \ ,
\end{equation*}
and introducing the stretched variable $Y=\delta^{-1} h$, it follows
\begin{equation}
\label{EcVH}
\left(r(\xi)-\delta Y z'(\xi)\right)  \cV^{2}=H_{\theta} \ .
\end{equation}
We insert the following expansions
$$
\cV^{2}=\cV^{2}_{0}(y_\alpha,Y)+ \delta \cV^{2}_{1}(y_\alpha,Y) + \mathcal{O}(\delta^2),
$$
$$
 H_{\theta}= \svV_0(\xi,Y) + \delta \svV_1(\xi,Y)
   + \mathcal{O}(\delta^2) \ 
 $$
 in equation \eqref{EcVH}. Then we perform the identification of terms with the same power in $\delta$. We obtain the equations  
\begin{gather}
\label{EcVsV0}
r(\xi) \cV^{2}_{0} = \svV_0\ ,
\\
\label{EcVsV1}
r(\xi) \cV^{2}_{1}-Y z'(\xi)\cV^{2}_{0} = \svV_1 \ .
\end{gather}
According to  \eqref{VOcdbis}, and \eqref{V1cd} for $\alpha=2$, there holds
\begin{equation}
\label{EcVH1}
r(\xi) \cV^{2}_{1}-Y z'(\xi)\cV^{2}_{0} =
r(\xi)   \Big[ \sH^{2}_{1}+Y \big(\cH\,\sH^{2}_{0}    +b_{\sigma}^{2} \sH^{\sigma}_{0}
-\frac{z'}{r}\sH^{2}_{0} 
\big)\Big] \;\mathrm{e}^{-\lambda Y}\ .
\end{equation}
\begin{lem}
\label{lAB}
The {\it main curvatures} $\kappa_{1}$, $\kappa_{2}$ of the interface $\Sigma$, and its {\it mean curvature} $\cH=\frac12(\kappa_1+\kappa_{2})$ are 
\begin{equation*}
\kappa_{1}=k \ ,\quad \kappa_{2}= \frac{z'}{r}\ , 
\quad \mbox{and} \quad \cH=\frac12\left(k+\frac{z'}{r}\right)
\ .
\end{equation*}
\end{lem}
This result together with  \eqref{EcVH1} yields
$$
r(\xi) \cV^{2}_{1}-Y z'(\xi)\cV^{2}_{0} =
r(\xi)   \Big[ \sH^{2}_{1}+  \frac{Y}{2} \, \left(k+\frac{z'}{r}\right)(\xi) \,\sH^{2}_{0}\Big] \;\mathrm{e}^{-\lambda Y}\ .
$$
We set $r(\xi) \sH^{2}_{j}(y_{\beta})=\sH^{+}_{j}\left(\XX(\xi)\right)$, for $j=0,1$. Then, we infer \eqref{v0}-\eqref{v1} from the equations  \eqref{VOcdbis}-\eqref{V1cd}.
 
 \subsection{Proof of the lemma \ref{lAB}}

The {\it main curvatures} are defined by  $\kappa_{1}=b_{1}^1$, $\kappa_{2}=b_{2}^2$, where $b_{\alpha\beta}$ is the curvature tensor.  In this coordinate system on $\cU_{-}$, the metric  is defined by
$$
g_{ij}(h)=\left<\xx_{i}(h), \xx_{j}(h) \right>_{\mathbb{R}^3}\ .
$$
From \eqref{Exx1}-\eqref{Exx2}, we obtain 
$$
g_{ij}(h)=
\left(\begin{array}{ccc}
\left(r'(\xi)-h z^{\prime\prime}(\xi)\right)^2 + \left(z'(\xi)+h r^{\prime\prime}(\xi)\right)^2 & 0 & 0 
\\
0 &  \left(r(\xi)-h z'(\xi)\right)^2 & 0 
\\
0 & 0 & 1
\end{array}\right)
$$
Recall that
\begin{equation*}
g_{\alpha\beta}(h)=a_{\alpha\beta}-2 b_{\alpha\beta} h + b_{\alpha}^{\gamma}b_{\gamma\beta} h^2
\, .
\end{equation*}
Hence, the curvature tensor on the interface $\Sigma$ is diagonal, and its diagonal components write
\begin{equation}
\label{Ect}
b_{11}=k
\quad \mbox{and} \quad b_{22}= r z' \ ,
\end{equation}
where $k=k(\xi)$ is defined by \eqref{Ek}. The inverse of the metric tensor $a_{\alpha\beta}$ in $\Sigma$ is diagonal and there holds
\begin{equation*}
\label{Emt}
a^{11}=1 \quad \mbox{and} \quad a^{22}= r^{-2} \ .
\end{equation*}
From \eqref{Ect}, we deduce the lemma.

\newpage
\bibliographystyle{plain}
\bibliography{biblio}

\end{document}